\date{} 
\title{Conformal invariance of planar loop-erased 
random walks and uniform spanning trees}
\author{Gregory F. Lawler\footnote {Duke University and
Cornell University;  %
partially supported by the  %
National Science Foundation and the Mittag-Leffler  %
Institute.}\and Oded Schramm\footnote {Microsoft Research.}\and  %
Wendelin Werner\footnote{Universit\'e Paris-Sud and IUF.}}  %
\documentclass[12pt,naturalnames]{article}
\usepackage{amsmath}
\usepackage{amssymb}
\usepackage{amsthm}
\usepackage{amsfonts}
\usepackage{graphicx}
\input labelfig.tex
\newif\ifhyper\IfFileExists{hyperref.sty}{\hypertrue}{\hyperfalse}

\ifhyper\usepackage{hyperref}\fi

\newif\ifdraft
\drafttrue
\def\note#1/{\ifdraft {\bf [#1]}\fi}
\long\def\comment#1{}
\numberwithin{equation}{section}
\numberwithin{figure}{section}

\newtheorem{theorem}{Theorem}
\numberwithin{theorem}{section}
\newtheorem{corollary}[theorem]{Corollary}
\newtheorem{lemma}[theorem]{Lemma}
\newtheorem{proposition}[theorem]{Proposition}

\theoremstyle{remark}
\theoremstyle{remark}\newtheorem{remark}[theorem]{Remark}
\def\eref#1{(\ref{#1})}
\def\QED{\qed\medskip}
\newcommand{\Prob} {{\bf P}}
\newcommand{\R}{\mathbb{R}}
\newcommand{\C}{\mathbb{C}}
\newcommand{\Z}{\mathbb{Z}}
\newcommand{\N}{\mathbb{N}}

\def\H{\mathbb{H}}
\def\U{\mathbb{U}}
\def\diam{\mathop{\mathrm{diam}}}

\def\dist{\mathop{\mathrm{dist}}}

\def\length{\mathop{\mathrm{length}}}
\def\Im{{\rm Im}\,}
\def\Re{{\rm Re}\,}
\def\SLEkk#1/{$\mathrm{SLE}_{#1}$}
\def\SLEk/{\SLEkk{\kappa}/}
\def\SLEtwo/{\SLEkk2/}
\def\SLE/{$\mathrm{SLE}$}
\def\Ito/{It\^o}
\def \eps {\epsilon}
\def \P {{\bf P}}
\def\md{\mid}
\def\Bb#1#2{{\def\md{\bigm| }#1\bigl[#2\bigr]}}
\def\BB#1#2{{\def\md{\Bigm| }#1\Bigl[#2\Bigr]}}
\def\Bs#1#2{{\def\md{\mid}#1[#2]}}
\def\Pb{\Bb\P}
\def\Eb{\Bb\E}
\def\PB{\BB\P}
\def\EB{\BB\E}
\def\Ps{\Bs\P}
\def\Es{\Bs\E}

\def \p {{\partial}}
\def \E {{\bf E}}

\def\closure{\overline}
\def\ev#1{{\mathcal{#1}}}
\def \proof {{ \medbreak \noindent {\bf Proof.} }}
\def\proofof#1{{ \medbreak \noindent {\bf Proof of #1.} }}

\def\bl{\bigl}\def\br{\bigr}\def\Bl{\Bigl}\def\Br{\Bigr}

\def\V#1{V(#1)}
\def\Vb#1{V_\p(#1)}

\def\capa{t}

\def\gauge{\Upsilon}

\def\XX{\mathcal X}

\def\ball{\mathcal{B}}

\def\inr#1{\mathrm{rad}_{#1}}
\def\inrad{\inr0}

\def\Doms{\mathfrak D}

\def\Dd{{D}}

\def\pdist{\rho}
\def\VV{Z}
\def\harm{{\mathfrak H}}
\def\W{W}
\def\TT{\bar t}  %
\def\GO{G_\to}
\def\V#1{V(#1)}
\def\Vd#1{V^\delta(#1)}
\def\Vbd#1{V^\delta_\p(#1)}
\def\WWzeta{W}
\def\KasteleynManhattan{MR28:2859}
\def\DuplantierManhattan{MR89b:82025}
\def\FukaiUchiyama{MR97m:60098}

\def\SmirnovPerc{MR1851632}

\def\LSWi{math.PR/9911084}
\def\LSWii{math.PR/0003156}

\def\LSWan{math.PR/0005295}
\def\SWpercexpo{math.PR/0109120}
\def\KenyonDD{MR99m:52026}
\def\SchSLE{MR1776084}
\def\Ahlfors{MR50:10211}

\def\CardyFormula{MR92m:82048}

\def\AizenmanBurchard{MR2000i:60012}
\def\ABNW{MR2001c:60151}
\def\Dudley{MR91g:60001}

\def\RSsle{math.PR/0106036}
\def\Lbook{MR92f:60122}

\def\Llerw{MR2000k:60092}
\def\LawlerLERWdef{MR81j:60081}
\def\LSWoneArm{math.PR/0108211}
\def\HaggstromUSTlim{MR97b:60170}
\def\PemantleUST{MR92g:60014}
\def\SchPercForm{math.PR/0107096}
\def\WilsonAlg{MR1427525}
\def\LyonsUSFsurvey{MR99e:60029}
\def\KenyonLERW{MR1819995}
\def\KenyonUST{MR1757962}
\def\PommeBDRY{MR95b:30008}
\def\Dubins{MR38:2837}
\def\Strassen{MR35:4969}
\def\RevuzYor{MR92d:60053}
\def\BSsqharm{MR98d:60134}

\def\KaratsasShreve{MR89c:60096}
\def\Collatz{MR22:322}

\def\noopsort#1{}

\begin{document}
\maketitle

\begin{abstract}
This paper proves that the scaling limit of loop-erased 
random walk in a simply connected domain $D\subsetneqq\C$
is equal to the radial \SLEkk2/ path.  In particular,
the limit exists and is conformally invariant.  It follows
that the scaling limit of the uniform spanning tree
in a Jordan domain exists and is conformally invariant.
Assuming that $\p D$ is a $C^1$ simple closed curve,
the same method is applied to show that the
scaling limit of the uniform spanning tree Peano curve,
where the tree is wired along a proper arc $A\subset\p D$,
is the chordal \SLEkk8/ path in $\closure D$ joining the endpoints
of $A$.   A by-product of this result is that 
 \SLEkk8/ is almost surely generated
by a continuous path.
The results and proofs are not restricted to a particular choice
of lattice.
\end{abstract}

\newpage
\tableofcontents
\newpage

\section {Introduction}
\subsection{Motivation from statistical physics}
One of the main goals of both
probability theory and statistical physics
is to understand the 
asymptotic behavior of random systems when the number of 
microscopic random inputs goes to infinity. 
These random inputs can be independent, such as 
a sequence of independent random variables,
or dependent, as in the Ising model.
Often, one wishes to 
understand these systems
via some relevant ``observables'' that can 
be of geometric or analytic nature.
In order to understand this asymptotic behavior, one 
can attempt to prove convergence towards a suitable continuous model.
 The simplest and  most important example of such random continuous
models is Brownian motion, which 
is the scaling limit of random walks. 
In particular, simple 
random walk on any lattice in $\R^d$ converges to
(a linear image of) Brownian motion in the scaling limit.

Physicists and chemists have observed that critical systems (i.e.,
systems at their phase transition point) can exhibit macroscopic
randomness. 
Hence, various quantities related to the corresponding lattice models
should converge as the mesh refines.
In fact, one of the important starting points for theoretical physicists
working on two-dimensional critical models
is the assumption that the continuous limit is independent of the 
lattice and, furthermore, displays conformal invariance.
This  assumption has
 enabled them to develop and use techniques from conformal field theory 
to predict  exact values of  certain critical 
exponents. Until very recently, the existence of the limit,
its conformal invariance, and the derivation of the exponents
assuming conformal invariance remained beyond mathematical justification
for the basic lattice models in critical phenomena, such as percolation,
the Ising model, and random-cluster measures.
Although there are many interesting questions about higher dimensional
systems, we will limit our discussion to two dimensions
where conformal invariance plays an essential role.  

\subsection{Recent progress}
In~\cite{\SchSLE}, a one-parameter family of random growth processes 
(loosely speaking, random curves) in two dimensions was introduced.
The growth process is based on Loewner's differential equation where
the driving term is time-scaled one-dimensional Brownian motion,
and is therefore called sto\-chastic Loewner evolution, or \SLEk/.
The parameter $\kappa\ge 0$ of \SLE/ is the time scaling constant
for the driving Brownian motion.  
It was conjectured that the scaling 
limit of the loop-erased random walk (LERW) 
is \SLEkk2/, and this conjecture was proved to be 
equivalent to the 
conformal invariance of the LERW scaling limit~\cite{\SchSLE}. 
The argument given was quite general and shows that
a conformally invariant random path satisfying a mild
Markovian property, which will be described below, 
must be \SLE/.
On this basis, it was also conjectured there
that the scaling limits of the critical
percolation interface and the uniform spanning tree Peano curve
are the paths of \SLEkk6/, and \SLEkk8/, respectively,
and it was claimed that conformal invariance is sufficient
to establish these conjectures.
(For additional conjectures regarding curves tending to \SLE/,
including the interfaces in critical random cluster
models --- also called FK percolation models --- 
for $q\in[0,4]$, see~\cite{\RSsle}.)

At some values of the parameter 
$\kappa$, \SLE/ has some remarkable properties.
For instance, \SLEkk6/ has a locality
property \cite {\LSWi} that makes it possible
to relate its outer boundary with that of 
planar Brownian motion. This has led
 to the proof of conjectures 
concerning planar Brownian motion and simple random walks
\cite {\LSWi,\LSWii,\LSWan}.
 
Smirnov \cite {\SmirnovPerc,Smirnov} recently
proved the existence and conformal invariance of the
scaling limit of critical site percolation on the
two-dimensional triangular lattice:  he managed 
to prove Cardy's formula~\cite{\CardyFormula}
 which is a formula
for the limit of the probability of a percolation crossing
between two arcs on the boundary of the domain.
Combining this information
with independence properties of percolation, Smirnov then 
showed that the
 scaling limit of 
the percolation interface is \SLEkk6/.
This has led to the rigorous determination of critical 
exponents for this percolation model~\cite{\LSWoneArm,\SWpercexpo}.

\subsection{LERW and UST defined}
The
uniform spanning tree (UST), which can be interpreted as the
$q=0$ critical random cluster  
model~\cite {\HaggstromUSTlim},
is a dependent model that
has many remarkable features.
In particular, it is very closely related to the loop-erased 
random walk, whose definition~\cite{\LawlerLERWdef} we now
briefly recall.

Consider any finite or recurrent connected graph $G$, a vertex $a$ and a set 
of vertices $V$.
{\em Loop-erased random walk} (LERW) from $a$ to $V$ is a random simple
curve joining $a$ to $V$ obtained by erasing the loops in
chronological order from
a simple random walk started at $a$ and stopped upon hitting $V$.
In other words, if $(\Gamma({n}), 0 \le n \le T)$ is a simple 
random walk on $G$ started from $a$ and stopped at its first hitting
time $T$ of $V$, the loop-erasure $\beta=(\beta_0,\dots,\beta_\ell)$ is defined 
inductively as follows:
 $\beta_0 = a$;  if $\beta_n\in V$, then $n=\ell$;
and otherwise
$\beta_{n+1} = \Gamma(k)$
where $k = {  1+ \max \{ m \le T \ : \ \Gamma(m) = \beta_n \}}$.

\begin{figure}
\centerline{\includegraphics*[height=2.3in]{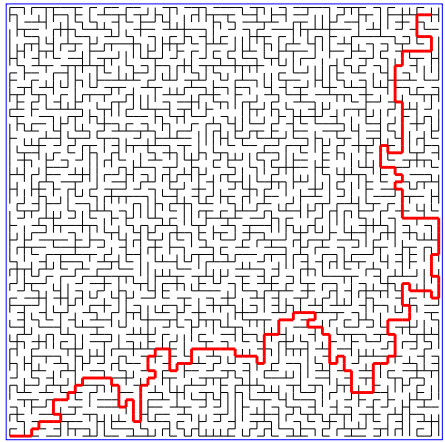}}
\caption{\label{f.ustlerw}The LERW in the UST.}
\end{figure}

A spanning tree $T$ of a connected 
graph $G$ is a subgraph of $G$ such that 
for every pair of vertices $v,u$ in $G$ there is 
a unique simple path (that is, self-avoiding) 
in $T$ with these vertices as endpoints. 
A {\em uniform spanning tree} (UST) in a finite, connected graph $G$ 
  is a sample from
the uniform probability measure on spanning trees of $G$.
It has been shown \cite {\PemantleUST} that the law of the 
self-avoiding path
with endpoints $a$ and $b$ in the UST is the same
as that of LERW from $a$ to $\{b\}$.  See Figure~\ref{f.ustlerw}.

David Wilson~\cite{\WilsonAlg} established
an even
 stronger connection  between LERW and UST by
giving an algorithm to generate USTs using LERW.
Wilson's algorithm runs as follows.
Pick an arbitrary ordering $v_0,v_1,\dots,v_m$ for the vertices
in $G$.  Let $T_0=\{v_0\}$.  Inductively, for $n=1,2,\dots,m$
define $T_n$ to be the union of $T_{n-1}$ and a (conditionally independent) 
LERW path from $v_n$ to $T_{n-1}$.  (If $v_n\in T_{n-1}$, then 
$T_n=T_{n-1}$.) Then, regardless of the chosen order of
the vertices,  $T_m$ 
is a UST on $G$.

Wilson's algorithm gives a natural extension of the definition of UST 
to infinite recurrent graphs.  In fact, for transient 
graphs, there are two natural definitions
 which often coincide, but this interesting
theory is somewhat removed from the topic of this paper.
Many striking properties of UST and LERW have been discovered. 
See~\cite{\LyonsUSFsurvey} for a survey of UST's and
\cite{\Llerw} for a survey of properties of LERW in $\Z^d$, $d>2$.

Exploiting a link with domino tilings and deriving discrete 
analogs of Cauchy-Riemann equations,
Richard Kenyon~\cite{\KenyonLERW,\KenyonUST}
rigorously established
the values of various critical exponents predicted for
the LERW \cite{GuttmannBursill, MajumdarLERW, DuplantierLERW}
in two dimensions.
In particular, he showed that the expected number of steps 
of a LERW joining two
corners of the $N \times N$
square in the square grid $\Z^2$ 
is of the order of magnitude of 
$N^{5/4}$.
He also 
 showed 
conformal invariance for the leading term 
in the asymptotics of the  
probability that the LERW contains a given edge. 
This was the first 
mathematical evidence for full conformal invariance of the LERW
scaling limit. 

In~\cite{\AizenmanBurchard,\ABNW} subsequential scaling limits of the
UST measures in $\Z^d$ were shown to exist,
using a compactness argument.
Moreover, these papers prove that all the paths  in the scaling limit
that intersect a fixed bounded region are  uniformly H\"older continuous.  
In~\cite{\SchSLE} the topology of subsequential scaling limits of the UST on
$\Z^2$ was determined.  In particular, it was shown that every subsequential
scaling limit of LERW is a simple path.

\subsection{A short description of \SLE/}
We now briefly describe   \SLE/; precise definitions
are deferred to Section~\ref{s.backgroundSLE}.
{\em Chordal \SLE/} is a random growing family
of compact sets $K_t,\: t \in [0,\infty)$, in the closure $\closure{\H}$
of the upper half plane $\H$.
The evolution of $K_t$ is given by the
Loewner differential equation with ``driving
function'' Brownian motion.
It is known~\cite{\RSsle} that when $\kappa\ne8$
the process is described by a random curve $\gamma:[0,\infty)\to\closure\H$,
in the sense that for every $t\ge 0$,
$\H \setminus K_t$ is the unbounded component of
$\H \setminus \gamma[0,t]$.    A corollary of our results
is that this holds for $\kappa=8$ as well. 
The curve $\gamma$   satisfies $\gamma(0)=0$ and
$\lim_{t\to\infty}\gamma(t)=\infty$. If $\kappa\le 4$, then
$\gamma$ is a simple curve and
 $K_t=\gamma[0,t]$.

There is another version of \SLE/ called {\em radial \SLE/}.
Radial \SLE/ also satisfies the description above,
except that the upper half plane $\H$ 
is replaced by the unit disk $\U$, $\gamma(0)$ is  
  on the unit circle $\p\U$  and $\lim_{t\to\infty}\gamma(t)=0$.

Both radial and chordal versions of \SLE/ may be defined
in an arbitrary simply connected domain $D\subsetneqq\C$
by mapping over to $D$ using a fixed conformal map $\phi$ from $\H$ or $\U$
to $D$.  

\subsection{The main results of the paper}

\begin{figure}
\centerline{\includegraphics*[height=2.3in]{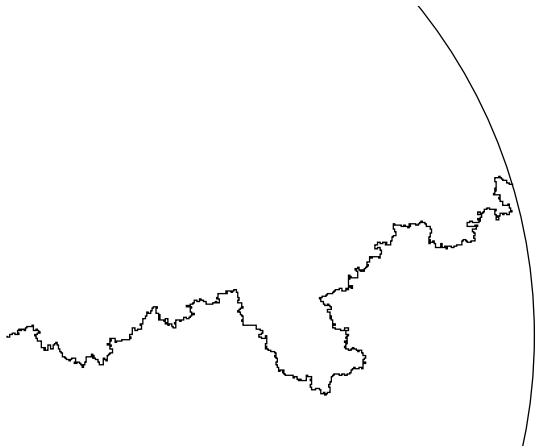}}
\caption{\label{f.lerw}A sample of the loop-erased
random walk; proved to converge to radial \SLEkk2/.}
\end{figure}

Let $D\subsetneqq\C$ be a simply connected
domain  with $0\in D$.  For $\delta>0$, let
$\mu_\delta$ be the law of the loop-erasure of simple random
walk   on the grid $\delta\Z^2$, started at $0$
and stopped when it hits $\p D$. 
See Figure~\ref{f.lerw}. 
Let $\nu$ be the law of the image of the radial \SLEkk2/ path
under a conformal map from the unit disk $\U$ to $D$ fixing $0$.
When the boundary of $D$ is very rough, the conformal
map from $\U$ to $D$ might not extend continuously to the boundary,
but the proof of the following theorem in fact
shows that
even in this case the image of the \SLEkk2/ path
has a unique endpoint on $\p D$.

On the space of unparameterized paths in $\C$, consider the
metric $\pdist(\beta,\gamma)=
\inf\,\sup_{t\in[0,1]}|\hat\beta(t)-\hat\gamma(t)|$,
where the infimum is over all choices of parameterizations
$\hat\beta$ and $\hat \gamma$
in $[0,1]$ of $\beta$ and $\gamma$.

\begin{theorem}[LERW scaling limit]\label{t.lerw}
The measures $\mu_\delta$ converge 
weakly to $\nu$  as $\delta\to0$ 
with respect to the metric $\pdist$
on the space of curves.
\end{theorem}

Since \SLE/ is conformally invariant by definition, 
 this theorem implies conformal invariance
of the LERW.  The theorem and proof apply 
also to some other walks on lattices in the plane
where the scaling limit of the walk is isotropic Brownian motion. 
It even applies in the non-reversible setting.
See Section~\ref{s.other} for further details.

\medskip

There are two distinct definitions for the UST corresponding to
a domain $D\subsetneqq\C$, as follows.
Let $G_F(D)$ denote the subgraph of $\Z^2$ consisting
of all the edges and vertices which are contained
in $D$.
If $G_F(D)$ is connected, then we refer to the UST on $G_F(D)$
as the UST on $D$ with {\em free} boundary conditions. 
Let $G_W(D)$ denote the graph obtained from $\Z^2$
by contracting all the vertices outside of $D$ to
a single vertex (and removing edges which become loops).
Then the UST on $ G_W(D)$ is the UST on $D$ with {\em wired}
boundary conditions.

Since the UST is built from LERW via Wilson's Algorithm,
it is not surprising that conformal invariance of 
the UST scaling limit 
should follow from that of the LERW scaling limit.
  In fact,~\cite[Thm.~11.3]{\SchSLE}
says just that.

\begin{corollary}[UST scaling limit]\label{c.ust}
The wired and free UST
scaling limits (as defined in~\cite{\SchSLE})
in a simply connected domain $D\subset\C$
whose boundary is a $C^1$ smooth simple closed curve
exist, and are conformally invariant.
\QED
\end{corollary}

One can easily show, using~\cite[Thm.~11.1.(i)]{\SchSLE},
that the wired tree depends continuously on the domain,
and hence for that case $D$ may be an arbitrary
simply connected domain.
However, some regularity assumption is needed
for the free UST scaling limit: conformal invariance
fails for the domain whose boundary contains the topologist's
sine curve (the closure of $\{x+i\sin(1/x):x\in(0,1]\}$). 

\medskip

The UST Peano curve
is an entirely different curve derived from the UST 
in two dimensions.  The curve is rather remarkable,
as it is a natural random path visiting every vertex in
an appropriate graph or lattice.
We now roughly describe two natural definitions of this curve;
further details appear in Section~\ref{s.peano}.

Let $G$ be a finite planar graph, with a particular embedding
in the plane, and $G^\dagger$ denote its planar dual,
again with a particular embedding.
Then there is a bijection $e\leftrightarrow e^\dagger$
between the edges of $G$ and those of $G^\dagger$, such
that for every edge $e$ in $G$,
$e\cap e^\dagger$ is a single point, and $e$ does
not intersect any other edge of $G^\dagger$.
Given a spanning tree $T$ of $G$, let $T^\dagger$
denote the graph whose vertices are the vertices
of $G^\dagger$ and whose edges are those edges
$e^\dagger$ such that $e\notin T$.  It is then easy
to verify that $T^\dagger$ is a spanning tree for
$G^\dagger$.
Therefore, if $T$ is a UST on $G$, then $T^\dagger$
is a UST on $G^\dagger$. 

\begin{figure} 
\centerline{\includegraphics*[height=2.0in]{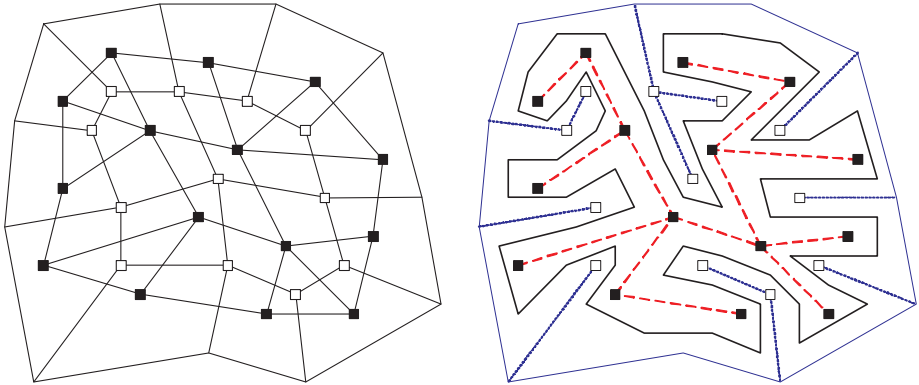}} 
\caption{\label{f.Peano}The graph, dual graph, tree, dual tree, and the 
Peano curve.  The vertex of the dual graph corresponding to the 
unbounded face is drawn as a cycle.} 
\end{figure}

The UST Peano curve is a curve that winds between $T$ and $T^\dagger$
and separates them.  More precisely, consider the graph $\hat G$
drawn in the plane by taking the union of $G$ and $G^\dagger$,
where each edge $e$ or $e^\dagger$ is subdivided into two edges
by introducing a vertex at $e\cap e^\dagger$.
The subgraph of the planar dual ${\hat G}^\dagger$
of $\hat G$
containing all edges which do not intersect $T\cup T^\dagger$
is a simple closed path --- the {\em UST Peano} path.
See Figure~\ref{f.Peano}. 

Some properties of the UST Peano path on $\Z^2$ have been
studied in the physics literature;
e.g.,~\cite{\KasteleynManhattan,\DuplantierManhattan}.  
There, it has been called the
Hamiltonian path on the Manhattan lattice.
The reason for this name is as follows.
On $\Z^2$, say,
orient each horizontal edge whose $y$-coordinate
is even to the right and each horizontal edge whose
$y$-coordinate is odd to the left.  Similarly,
orient down each vertical edge whose $x$-coordinate is even,
and orient up each vertical edge whose $x$-coordinate is odd.
Now rescale the resulting oriented graph by $1/2$ and
translate it by $(1/4,1/4)$.  It is easy to check that
a Hamiltonian path (a path visiting every
vertex exactly once) respecting the orientation on the resulting
oriented graph is the same as the UST Peano path of $Z^2$.
It should be expected that the uniform measure on Hamiltonian
paths in $\Z^2$ has the same scaling limit as that of
the UST Peano path.

Given a domain $D$, one can consider the UST Peano curve
for the wired or for the free UST (which is essentially the same as 
the wired, by duality).  
However, the conjecture from~\cite{\SchSLE} regarding the convergence to 
chordal \SLE/ pertains to the UST Peano curve associated with the tree with
mixed wired and free conditions.

Let $D\subset\C$ be a domain whose boundary is
a $C^1$-smooth simple closed curve,
and let $a,b\in\p D$ be distinct boundary points.
Let $\alpha$ and $\beta$ denote the two complementary arcs
of $\p D$ whose endpoints are $a$ and $b$.
For all $\delta>0$, consider an approximation $G_\delta$ of
the domain $D$ in the grid $\delta\Z^2$.
(A precise statement of what it means for $G_\delta$ to
be an approximation of $D$ will be given in Section~\ref{s.peano}.)
Let $\gamma_\delta$ denote the Peano curve associated to
the UST on $G_\delta$ with wired boundary near $\alpha$ and
free boundary near $\beta$.  
Then $\gamma_\delta$ may be considered
as a path in $D$ from a point near $a$ to a point near $b$.

\begin{figure} 
\centerline{\includegraphics*[height=2.0in]{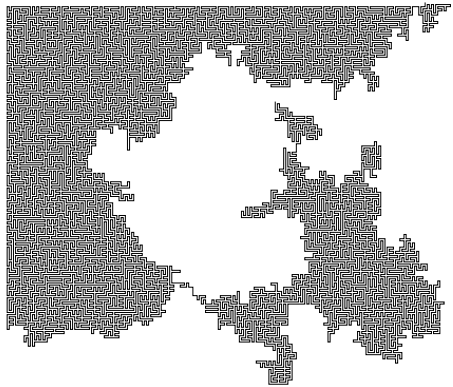}} 
\caption{\label{f.finepeano}An arc from a sample of the UST Peano path;
proved to converge to chordal \SLEkk8/.} 
\end{figure} 

\begin{theorem}[UST Peano path scaling limit]\label{t.intropeano} 
The UST Peano curve scaling limit in $D$
with wired boundary on $\alpha$ and free boundary on
$\beta$ exists, and is equal
to the image of the chordal \SLEkk8/ path under
any   conformal map from $\H$ to $D$ mapping
$0$ to $a$ and $\infty$ to $b$.
\end {theorem}

Again, the convergence is weak convergence of measures
with respect to the metric $\rho$.
Figure~\ref{f.finepeano} shows a sample of the UST Peano path on
a fine grid.

As explained above, it was proved in~\cite{\RSsle} that each
\SLEk/ is generated by a path, except for $\kappa=8$.  
In Section~\ref{s.peanoconclusions},
the remaining case $\kappa=8$ is proved,
 using the convergence of the Peano curve.

Corollary~\ref{c.ust} 
and Theorem~\ref{t.intropeano} (and their proofs)
 apply to other reversible walks on planar lattices
(the self-duality of $\Z^2$ does not play an important role);
see Section~\ref{s.other}.

\medskip 

To add perspective, we note that the convergence to SLE of the LERW and the
UST Peano curve are two boundary cases of the conjectured
convergence~\cite{\RSsle} of the critical FK random cluster
measures with parameter $q\in(0,4)$.
For these parameter values,
the scaling limit of the interface of a critical cluster with mixed
boundary values is conjectured to converge to chordal \SLEkk\kappa(q)/,
where $ \kappa(q)={4\pi}/{\cos^{-1}\bigl(-\sqrt{q}/2\bigr)} $.
The boundary case $\kappa(0)=8$ corresponds to the convergence
of the UST Peano path to \SLEkk8/.

The outer boundary of the scaling limit of a macroscopic critical cluster
is not the same as the scaling limit of a critical cluster outer boundary,
because of \lq\lq fjords\rq\rq\ which are pinched off in the limit. 
The former
is conjectured to \lq\lq look like\rq\rq\ \SLEkk{16/\kappa(q)}/, but a precise
form of this conjecture is not yet known.
In the case $q=0$, however, such a correspondence is easy to explain.
In $\Z^2$, an arc of the Peano curve is surrounded on one side by
a simple path in the tree, and on the other side by a simple
path in the dual tree.
Both these paths are LERW's.  Similar correspondences exist for
the UST in a subdomain of $\R^2$, but one has to set appropriate
boundary conditions.  Thus, the convergence of LERW to \SLEkk2/
also corresponds to the case $q=0$, as $16/\kappa(0)=2$.

\medskip

Suppose that $0\in D$, and $\alpha,\beta\subset \p D$,
as before.  Consider the simple random walk on $\delta\Z^2$
which is reflected off $\beta$ and stopped when it
hits $\alpha$.
Using an analogous method to the one of the 
present paper, one could handle the scaling limit
of the loop-erasure of this walk. 
It is described
by a variant of \SLEtwo/ where the driving term
is Brownian motion with time scaled by $2$, but
having an additional drift.  The drift is not
constant, but can be explicitly computed.

\medskip 

The identification of the scaling limit as one of the \SLE/'s 
should 
facilitate the derivation of critical exponents and also the
asymptotic probabilities of various events, including
some results which have not been predicted by arguments from physics.
This was the case for critical site percolation on the
triangular grid~\cite{\SmirnovPerc,\SchPercForm,\LSWoneArm,\SWpercexpo}. 

\subsection{Some comments about the proof}

Since loop-erased random walk is obtained in a deterministic
 way from
simple random walk (by erasing its loops) and since 
simple random walk converges to Brownian motion in the scaling limit, 
it is natural to think that 
the scaling limit of LERW should simply be the process
obtained by erasing the loops from a planar Brownian motion.
The problem with this approach is that planar
Brownian motion has loops at every scale, so that there is 
no simple algorithm to erase loops.  In particular,
 there
is no ``first'' loop.
Our proof does use the relation
between LERW and simple random walks, combined with the 
fact that quantities related to simple random walks,
such as hitting probabilities, converge to their 
continuous conformally invariant counterparts.

The proof of each of our main theorems is naturally divided into two parts.
The first part establishes the convergence to \SLE/ with
respect to a weaker topology than the topology induced by
the metric $\pdist$ of paths, namely, we show that the Loewner
driving process for the discrete random path converges to
a Brownian motion.  
This part of the proof, which we consider to be the more important one,
is essentially self-contained.  The second part uses some  regularity
properties of the discrete processes from~\cite{\SchSLE} to 
prove convergence with respect to the stronger topology.

The method for the first part  can be considered
as a rather general method for identifying the
scaling limit of a dependent system that is
conjectured to be conformally invariant.  It
requires having some  ``observable'' quantity 
that can be estimated
well and
  a mild Markovian property, which we now describe.
Suppose that to every simply connected domain $D$ containing $0$ there is
associated a random path $\gamma$ from $\p D$ to $0$ (e.g., the orientation
reversal of LERW).  The required property is
that if $\beta$ is an arc with one endpoint in $\p D$ and
we condition on $\beta\subset\gamma$ (assuming this has positive
probability, say), then the conditioned
distribution of $\gamma\setminus\beta$ is the same as the
random path in the domain $D\setminus\beta$ conditioned to
start at the other endpoint $q$ of $\beta$. 
(Thus, $(D\setminus\beta,q)$ is the state of a Markov chain 
whose transitions correspond to adding edges from $\gamma$ 
to $\beta$ and modifying $q$ appropriately.) 
Interestingly, among the discrete processes conjectured to converge
to \SLE/, the LERW is the only one where the verification of
this property is not completely trivial.
(For LERW it is not trivial, but not difficult; 
see part 3 of Lemma~\ref{l.mark}.) 
The statement of this property for the UST Peano curve
is in Lemma~\ref{l.pmark}.  The fact that \SLE/ satisfies this
property follows from the Markovian property of its
driving Brownian motion.

The particular choice of observable is not so important. 
What is essential is that one can conveniently calculate the asymptotics of the
observable for appropriate large-scale configurations. 
The particular observable that we have chosen for the LERW 
convergence is the
expected number of visits to a
vertex $v$ by the simple random walk generating the LERW.
Conformal invariance is not assumed but comes
out of the calculation --- 
hitting probabilities for random walks are discrete
harmonic functions, which converge to continuous harmonic functions.
One technical issue is to establish this convergence without
any boundary smoothness assumption.
Once the observable has been approximated, 
  the conditional expectation and variance of increments of 
the Loewner driving function for the discrete process
can be estimated,
and standard techniques (the Skorokhod embedding) can be used to show that
this random function
approaches the appropriate Brownian motion.

Although Theorem~\ref{t.intropeano} can probably be
derived with some work from Corollary~\ref{c.ust}, instead, to illustrate
our method we prove it by applying again the same general
strategy of the proof of Theorem~\ref{t.lerw}, with the choice
of a different observable. 

Actually, it is easier to explain the main ideas
behind the proof of Theorem~\ref{t.intropeano}.
Fix some vertex $v$ in $D$ and a subarc $\alpha_1\subset\alpha$.
Let $\ev A$ be the event that the UST path (not the Peano-path, but
the path contained in the UST) from $v$ to $\alpha$
hits $\alpha_1$.    By Wilson's algorithm the
probability of $\ev A$ 
is the same as the probability that simple random walk
started at  $v$ reflected off $\beta$  first hits 
$\alpha$ in $\alpha_1$. The latter probability
 can be estimated
directly.
If $\gamma[0,n]$ denotes the restriction of the Peano path
to its first $n$ steps, then $\Pb{\ev A\md\gamma[0,n]}$, the
probability of $\ev A$ conditioned on $\gamma[0,n]$, is clearly a martingale
with respect to $n$.  But,
by the Markovian property discussed above, 
 the value of $\Pb{\ev A\md\gamma[0,n]}$
may be estimated in precisely the same way that $\Ps{\ev A}$ is
estimated.  The estimate turns out to 
be a function of the conformal 
geometry of the configuration
$(v,D\setminus\gamma[0,n],\gamma(n),\alpha_1,\beta)$.
Knowing that this is a martingale for two appropriately chosen
vertices $v$ is sufficient to characterize the large
scale behavior of $\gamma$.

\medskip

As mentioned above,
in the case of LERW, the observable we chose to look at
is the expected number of visits to a fixed vertex $v$
by the simple random walk $\Gamma$ generating the LERW $\gamma$.
The walk $\Gamma$ can be considered as the union of
$\gamma$ with a sequence of
loops $\Gamma^j$ based at vertices of $\gamma$.
We look at the conditioned expectation of the number of
visits of $\Gamma$ to $v$ given an arc $\tilde\gamma$ of $\gamma$
adjacent to the boundary of the domain.  This is clearly
a martingale with respect to the filtration obtained
by taking larger and larger arcs $\tilde\gamma\subset\gamma$.
This quantity falls into two parts: the visits to
$v$ in the loops $\Gamma^j$ based at $\tilde\gamma$, and those that are not.
Each of these two parts can be estimated well by random-walk calculations.
Translating the fact that this is a martingale to information
about the Loewner driving process for $\gamma$ inevitably leads to
the identification of this driving process as appropriately
scaled Brownian motion.

\medskip

Actually, we first had a longer proof of convergence 
of LERW to \SLEtwo/, based on the fact that it is 
possible to construct the hull of a Brownian motion
by adding Brownian loops to \SLEtwo/.  This
can be viewed as 
a particular case of the restriction properties of 
\SLEkk\kappa/ with Brownian loops added, which we  study 
in the subsequent paper \cite {LSWinprep}.  
Let us also mention the following related open 
 question.
Consider a sequence of simple random walks $S^k(n)$ on 
a lattice with lattice spacing $\delta_k\to 0$, from $S^k(0)=0$ to $\p\U$,
and let $\gamma^k$ denote the corresponding loop-erased paths.
Theorem~\ref {t.lerw} shows that one can find a subsequence 
such that the law of the pair $(\gamma^k,S^k)$
converges to a coupling of \SLEtwo/ with
Brownian motion.  (That is, a law
for a pair $(X,Y)$, where $X$ has the
same distribution as the \SLEtwo/ path and $Y$ has
the same distribution as Brownian motion.)
The question is whether in this coupling,
the \SLEtwo/ is a deterministic 
function of the Brownian motion.
In other words, is it possible to show that this
is not a deterministic procedure to erase loops from a Brownian motion?

\section{Preliminaries}\label{s.prel}

The reading of this paper requires some background knowledge in
several different
fields.  Some background about Loewner's equation and \SLE/ is reviewed
in the next subsection.  It is assumed that the reader is familiar
with some of the basic properties of Brownian motion (definition, 
strong Markov property, etc.).  Some of the basic
properties of conformal maps (Riemann's mapping theorem,
compactness, Koebe distortion)
are also needed for the proof.  This material may be learned from the first
two chapters of~\cite{\PommeBDRY}, for example.  
In terms of the theory of conformal mappings, this suffices
for understanding the argument showing that the driving process of
the LERW converges to Brownian motion.  For improving the topology
of convergence, some familiarity with the notion of extremal length
(a.k.a.\ extremal distance) is also required.  A possible source
for that is~\cite{\Ahlfors}.  The reader would also need to know
some of the very basic properties of harmonic measure.

\subsection {Loewner's equation and \SLE/}\label{s.backgroundSLE}

We now review some facts concerning Loewner's equations 
and stochastic Loewner evolutions. 
For more details, see e.g., \cite {\SchSLE, \RSsle, \LSWi, \LSWii}.

Suppose that $D\subsetneqq\C$ is a simply connected domain
with $0\in D$.  Then there is a unique conformal
homeomorphism $\psi=\psi_D:D\to\U$ which is onto the unit
disk $\U=\{z\in\C:|z|<1\}$ such that 
$\psi_D(0)=0$ and $\psi_D'(0)$ is a positive real.
If $D\subset\U$, then $\psi_D'(0)\ge 1$,
and $\log\psi_D'(0)$ is called the {\em capacity}
of $\overline\U\setminus D$ from $0$.

Now suppose that $\eta : [0, \infty] \to \closure\U$ is a continuous 
simple curve in the unit disk with  
$\eta (0) \in \partial \U$, $\eta (\infty)=0$ and
$\eta (0,  \infty] \subset \U$.
For each $t \ge 0$, set $K_t:= \eta[0,t]$,
$U_t:=\U\setminus K_t$ and $g_t:=\psi_{U_t}$.
Since $t \mapsto g_t'(0)$ is increasing (by the Schwarz Lemma, say),
one can reparametrize the
path in such a way that $g_t' (0) = \exp (t)$.
If that is the case,  we say that $\eta$ 
is parametrized by {\em capacity} from $0$.
By standard properties of conformal maps
(\cite[Proposition 2.5]{\PommeBDRY}), for each $t\in[0,\infty)$ 
the limit 
$$
\WWzeta(t):=\lim_{z\to \eta(t)} g_t(z)
\,,
$$
where $z$ tends to $\eta(t)$ from within $\U\setminus\eta[0,t]$
exists.
One can also verify that
$$
\WWzeta:[0,\infty)\to\p\U
$$
is continuous.
Assuming the parameterization by capacity,
Loewner's theorem states that $g_t$ satisfies the differential
equation
\begin {equation}
\label {e.lode}
\partial_t g_t (z) = - g_t (z) \,
\frac { g_t (z) + \WWzeta (t) }{g_t (z) - \WWzeta (t)}
\,.
\end {equation}
It is also clear that 
\begin{equation}\label{e.initial}
\forall z\in\U\qquad g_0(z)=z\,.
\end{equation}
We call $(\WWzeta (t), t \ge 0)$ the {\em driving function}
of the curve $\eta$.

The driving function $\WWzeta$ is sufficient to recover
the two-dimensional path $\eta$, because the
procedure may be reversed, as follows. 
Suppose that $\WWzeta:[0,\infty)\to \p\U$ is continuous.
Then for every $z\in\closure\U$
there is a solution $g_t(z)$ of
the ODE~\eref{e.lode} with initial value $g_0(z)=z$
up to some time $\tau(z)\in(0,\infty]$, beyond which the solution
does not exist.  
In fact, if $\tau(z)<\infty$ and $z\ne \WWzeta(0)$, then we have
$\lim_{t\uparrow\tau(z) } g_t(z)-\WWzeta(t)=0$, since this is the
only possible reason why the ODE cannot be solved beyond time $\tau(z)$.
Then one defines $K_t:=\{z\in\closure\U:\tau(z)\le t\}$
and $D_t:=\U\setminus K_t$ is the domain of definition of $g_t$.
The set $K_t$ is called the {\em hull} at time $t$.
If $\WWzeta$ arises from a simple path $\eta$ as described in
the previous paragraph, then we can recover $\eta$ from $\WWzeta$
by using $\eta(t)=g_t^{-1}\bl(\WWzeta(t)\br)$.
However, if $\WWzeta:[0,\infty)\to\U$ is an arbitrary continuous
driving function, then in general $K_t$ need not be a path,
and even if it is a path, it does not have to be a simple
path.

{\em Radial \SLEk/} is the process $(K_t,t\ge 0)$, where
the driving function $\WWzeta(t)$ is set to be
$\WWzeta(t) := \exp \bl(i B_{\kappa t}\br)$,
where $B:[0,\infty)\to\R$ is Brownian motion.
Often, one takes the starting point $B_0$
to be random uniform in $[0,2\pi]$.
It has been shown~\cite {\RSsle} that the hull
$K_t$ is a.s.\ a simple curve for every $t>0$ if $\kappa\le 4$
and that a.s.\ for every $t>0$ $K_t$ is not a simple curve
if $\kappa>4$.
For every $\kappa\ge 0$, there is a.s.\ some random continuous path 
$\eta:[0,\infty)\to\closure\U$ such that for all $t>0$,
$D_t$ is the component of $\U\setminus\eta[0,t]$ containing
$0$.  When $\kappa\ne8$, this was proved in~\cite{\RSsle},
while for $\kappa=8$ this will be proven in the current paper.
This path is called the radial \SLE/ path.

Suppose that $D$ is a simply connected domain containing $0$.
If $\gamma$ is a continuous simple
curve joining $\partial D$ to $0$ with only an endpoint in
$\p D$, one can reparametrize
the path $\eta:=\psi\circ\gamma$ according to capacity and find
its driving function
$\WWzeta$, as before.
The conformal map 
$$
\hat g_t=\psi_{ D \setminus \gamma [0,t]}:D \setminus \gamma [0,t] \to \U
$$
still satisfies~\eref {e.lode}, but this time, $\hat g_0 = \psi_D$.
(Here, the parameterization chosen for $\gamma$ is
according to the capacity of $\psi\circ\gamma[0,t]$.)
Radial \SLE/ in $D$ is then simply the image under $\psi_D^{-1}$ 
of radial \SLE/ in the unit disk.

\medskip

Similarly, one can encode continuous 
simple curves $\eta$  from $0$ to $\infty$ in the 
closed upper half-plane $\closure\H$ via a variant
of Loewner's equation.
For each time $t \ge 0$, there is a unique conformal map $g_t$
from $H_t := \H \setminus \eta  [0,t]$ onto $\H$ satisfying the
so-called {\em hydrodynamic normalization}
\begin{equation}
\label{e.hydro}
\lim_{z\to\infty} g_t(z)-z = 0\,,
\end{equation}
where $z \to \infty$ in $\H$.
If we write $g_t(z)=z+a(t)\,z^{-1}+ o(z^{-1})$ near $\infty$,
it turns out that $a(t)$ is monotone. 
Consequently, one can reparametrize $\eta$ in such a way that
$a(t)=2\,t$, that is $g_t(z) = z + 2\,t\,z^{-1} + o(z^{-1})$ when
 $z \to \infty$.  This parameterization of $\eta$ is
called the parameterization by {\em capacity} from infinity.
(This notion of capacity is analogous to the notion of capacity
in the radial setting, however, these are two distinct notions and
should not be confused.)
If $g:\H\setminus K\to\H$ is the conformal homeomorphism
satisfying the hydrodynamic normalization,
then $\lim_{z\to\infty} (g_t(z)-z)\,z/2$ is
called the {\em capacity} of $K$ from $\infty$.
Assuming that $\eta$ is parameterized by capacity,
the following analogue of Loewner's equation holds:
\begin {equation}
\label {e.clode}
\forall t>0\ 
\forall z\in H_t\qquad
\partial_t g_t (z) = \frac {2}{g_t (z) - \WWzeta (t)}\,,
\end {equation} 
where the driving function $\WWzeta$ is again defined by
$\WWzeta(t):=g_t\bl(\eta(t)\br)$.
As above, $\eta$ is determined by $\WWzeta$.

Conversely, suppose that $\WWzeta$ is a real-valued continuous
function.  For $z\in\closure\H$, one can solve the
differential equation~\eref {e.clode} starting with
$g_0(z)=z$, up to the first time $\tau(z)$ where $g_t(z)$ and $\WWzeta(z)$
collide (possibly, $\tau(z)=\infty$).
Let the hull be defined by $K_t:=\{z\in\closure\H: \tau(z)\le t\}$.
Then $g_t:\H\setminus K_t\to\H$ is a conformal map onto $\H$,
and $g_0 (z)=z$.
In general, $K_t$ is not necessarily a simple curve.
If $\WWzeta(t)= B_{\kappa t}$, then $( K_t , t \ge 0)$
is called {\em chordal} \SLEk/.

It turns out~\cite[\S 4.1]{\LSWii} that the local properties of
chordal \SLEk/ and of radial \SLEk/ are essentially the same.
(That is the reason why the normalization $a(t)=2\,t$ was
chosen over the seemingly more natural $a(t)=t$.)
In particular, for every $\kappa$ chordal \SLEk/
is generated by a random continuous path, called the
chordal \SLEk/ path.

\medskip

At some points in our proofs, we will need the following 
simple observation:

\begin{lemma}[Diameter bounds on $K_t$]\label{l.diambound}
There is a constant $C>0$ such that the following always holds.
Let $\WWzeta:[0,\infty)\to\R$ be continuous and let $(K_t, t\ge 0)$ 
be the corresponding hull for Loewner's chordal equation~\eref{e.clode} 
with driving function $\WWzeta$.
Set 
$$
k(t):= \sqrt t+ \max\bl\{|\WWzeta (s)-\WWzeta(0)|:s\in[0,t]\br\}\,.
$$
Then
$$
\forall t\ge 0\qquad
C^{-1}k(t) \le \diam K_t  \le C\,k(t)\,.
$$
Similarly,
when $K_t\subset\closure\U$ is the radial hull
for a continuous driving function $\WWzeta:[0,\infty)\to\p\U$,
then
$$
\forall t\ge 0\qquad
C^{-1}\min\bl\{k(t),1\br\} \le \diam K_t  \le C\,k(t)\,.
$$
\end{lemma}

\proof
This lemma can be derived by various means.
We will only give a detailed argument in the 
radial case. The chordal case is actually easier and 
can be derived using the same methods. It can also be 
seen as a consequence of the result in the radial setting
(because chordal Loewner equations can be interpreted as 
 scaling limits of radial Loewner equations).

 We start by proving the upper bound on $\diam K_t$.
Let $\delta\ge \max \bl\{ |\WWzeta (s)-\WWzeta(0)|:s\in[0,t]\br\}$.
Then, as long as $\bl|g_t (z)-\WWzeta(0)\br| \ge 3 \delta$, 
we have $|\partial_t g_t (z)| \le 1/ \delta$.
Hence, if $\bl|z-\WWzeta(0)\br| \ge 4 \delta$, then for all $t \le \delta^2$, 
$|g_t (z) - z| \le \delta$, and therefore $z \notin K_t$.  
Hence, $\diam K_t \le 8 k(t)$. 

In order to derive the lower bound, we will compare capacity with harmonic 
measure. 
It is sufficient to 
consider the case 
where $\diam K_t<1/10$.
Let $\mu$ denote the harmonic measure on $K_t\cup\p\U$ from $0$.
Because $K_t$ is contained in the disk of radius $\diam K_t$
with center $\WWzeta (0)\in K_t\cap\p\U$, there is a universal
constant $c$ such that $\mu(K_t)\le c\,\diam K_t$.
Hence, it suffices to give a lower bound for $\mu (K_t)$.

Since $g_t(z)/z$ is analytic and nonzero in a neighborhood
of $0$, the function $h(z)=\log|g_t(z)|-\log |z|$
is harmonic in $U_t:=\U\setminus K_t$.
Note that $h(0)=t$.  Because $|g_t(z)|\to 1$
as $z$ tends to the boundary of $U_t$, 
the mean value property of $h\circ g_t^{-1}$ implies the following
relation between harmonic measure and capacity:
$
t=h(0) = \int \log(1/|z|)\,d\mu(z)\,
$.
Since $K_t$ contains points in $\p\U$ and $\diam K_t
\le 1/10$, we have
$ \log(1/|z|)\le  c' \diam K_t$ for all $z \in K_t$.
Therefore,
$ t\le  c'\, \mu (K_t) \, \diam K_t \, \le
 c''\,(\diam K_t)^2 $.

It now remains to compare  $\mu (K_t)$ 
and $|\WWzeta (t)-\WWzeta (0)|$. 
We still assume that  
$\diam K_t <1/10$.
Let $A_t:= \p\U\setminus g_t(\closure\U\setminus K_t)$.
If $z\in \p\U\setminus K_s$
and $s \le t$
 then~\eref{e.lode} shows that $\p_u |g_{s+u}(z)-\WWzeta (s)|\ge 0$
at $u=0$.  This implies that $(A_s, s \le t)$ is 
non-decreasing.
Hence, for all $s \le t$,
we have $\WWzeta (0)\in A_0\subset A_t$ and $\WWzeta (s) 
\in A_s \subset A_t$ so that
$
|\WWzeta (s)-\WWzeta (0)|$ is bounded by the 
length of $A_t$, which is equal to $ 2\,\pi\,\mu(K_t)$.
This completes the proof of the lemma.
\QED

\subsection {A discrete harmonic measure estimate} \label{s.RWbg}

In this section we introduce some notation 
and state an  estimate relating discrete harmonic
measure and continuous harmonic measure in domains in the plane.
In order to get more quickly to the core of our method
in Section~\ref{s.core},
we postpone the proof of the harmonic measure estimate to Section~\ref{s.RW}.

A {\em grid domain} $D$ is a domain whose boundary
consists of edges of the grid $\Z^2$.
For an arbitrary domain $D\subset\C$, and $p\in D$ define
the inner radius of $D$ with respect to $p$,
$$
\inr p (D) := \inf \{ |z-p| \ : \ z \notin D \}\,.
$$
Let $\Doms$ denote the set of all simply connected grid domains 
such
that 
$0 < \inrad (D) < \infty$ (i.e., $D \not= \C$ and $0 \in D$).

Points in $\R^2=\C$  with integer coordinates will be called
{\em vertices}, 
or {\em lattice points}.
Let $\V D:=D\cap\Z^2$ denote the lattice points in $D$.

Let $D\in\Doms$ and $v$ a vertex in $\p D$.
If $\p D$ contains more than one edge
incident with $v$, then it may happen that the intersection
of $D$ with a small disk centered at $v$
will not be connected.  Hence, as viewed from $D$,
$v$ appears as more than one vertex.  In particular,
$\psi=\psi_D$ does not extend continuously to $v$.
This is a standard issue in conformal mapping theory,
which is often resolved by introducing the notion of prime
ends. 
But in the present case, there is a simpler solution which suffices for
our purposes.  Suppose $v\in \Z^2\cap\p D$,
and $e$ is an edge incident with $v$ that intersects $D$.
The set of such pairs $w=(v,e)$ will be denoted $\Vb D$.
If $\psi:D\to\U$ is conformal, then $\psi(w)$ will be shorthand
for the limit of $\psi(z)$ as $z\to v$ along $e$
(which always exists, by~\cite[Proposition 2.14]{\PommeBDRY}).  
Similarly, if a random walk first exists $D$ at $v$,
we say that it exited $D$ at $w$ if the edge $e$ was used
when first hitting $v$.
A reader of this paper who chooses to be sloppy and not distinguish
between $v$ and $w$
will not loose anything in the way of substance.
We will not always be so careful to make this distinction.

If $a\in \V D$ and $b\in\V D\cup\Vb D$,
define $H(a, b)= H_D (a,b)$ as the 
probability that  simple random walk started from $a$ and
stopped at its first exit time of $D$ visits $b$.

For any $w  \in D$ and $u \in \Vb D$, we define
\begin{equation}
\label{e.lambdadef}
\lambda = \lambda (w, u ; D) 
:= \frac {1- | \psi (w) |^2 }{| \psi(w) - \psi(u)|^2}
= \Re \left( \frac { \psi(u) + \psi (w)}{\psi(u) - \psi (w)} 
\right).
\end{equation}
Note that $\lambda$ is also equal to the imaginary part of the 
image of $w$ by the conformal map from $D$ onto the
upper
half-plane that maps $0$ onto $i$ and $u$ to $\infty$.
It is also the limit when $\eps \to 0$ of the ratio between the 
harmonic measure in $D$ of the $\eps$ neighborhood of $u$ in $\partial D$, 
taken respectively at $w$ and at $0$ (that is,
it corresponds to the Poisson kernel).
Therefore, $\lambda$ can be viewed as the continuous analog of
$H(w,u)/ H(0, u)$.  Note that the function $h(w)=H(w,u)/H(0,u)$
is {\em discrete harmonic}, on $\V D$, which means
that $h(w)$ is equal to the average of $h$ on the neighbors
of $w$ when $w\in\V D$.

\begin{proposition}[Hitting probability]\label{p.phit}
For every $\eps>0$ there is some $r_0>0$ such that the following
holds.  Let $D\in\Doms$ satisfy $\inrad(D)>r_0$, let
$u\in \Vb D$ and $w\in\V D$.
Suppose $|\psi_D(w)|\le1-\eps$ and $H(0,u)\ne0$.  Then
\begin{equation}\label{e.phit}
\left|\frac{H(w,u)}{H(0,u)} -
\lambda (w,u; D)
\right|<\eps
\,.
\end{equation}
\end{proposition}

The proof is given in Section~\ref{s.RW}.

\section {Conformal invariance of LERW}
\label{s.lerw}

\subsection {Loop-erased random walk background}
\label{s.lerwbackground}

We now recall some well-known facts concerning loop-erased 
random walks.  

\begin{lemma}[LERW reversal]\label{l.reverse}
Let $D\in\Doms$
and let $\Gamma$ be simple random walk
from $0$ stopped when it hits $\p D$.
Let $\beta$ be the loop-erasure of $\Gamma$,
and let $\gamma$ be the loop-erasure of the
time reversal of $\Gamma$.
Then $\gamma$ has the same distribution as
the time-reversal of $\beta$.
\end{lemma}

See~\cite{\Lbook}.  A simpler proof follows immediately
from the symmetry of equation (12.2.3) in~\cite{\Llerw}.
This result (and the proofs) also holds if we condition $\Gamma$
to exit $\partial D$
at a prescribed $u \in \Vb D$, which 
correspond to the event $\bl \{ \gamma \cap \partial D = \{u\}\br \}
= \bl\{ \beta \cap \partial D = \{ u \}\br\} $ (assuming this has positive
probability).

\medskip

Throughout our proof we will use the simple random walk $\Gamma$ and
the loop-erasure 
 $\gamma = ( \gamma_0, \gamma_1, \ldots, \gamma_\ell)$
of its time-reversal (so that
 $\gamma_0 \in \partial D$ and $\gamma_\ell = 0$).
We use $D_j$ to denote the grid domains
$D_j := D \setminus  \bigcup_{i=0}^{j-1}  
[\gamma_i , \gamma_{i+1}]$.
Define for $j \in \{ 0,1, \ldots, \ell\}$,
$$
n_j := \min \{ n \ge 0 \ : \ \Gamma (n) = \gamma_j \},
$$
and note that $n_{j+1}<n_j$ for $j=0,1,\dots,\ell-1$, by
the definition of $\gamma$.
Also set
$$
\Gamma^{j+1} := \Gamma [n_{j+1}, n_{j}]
\,.$$
More precisely, consider $\Gamma^j$ as the grid-path
given by
$$
\Gamma^j(m):=\Gamma(m+n_j),\qquad m=0,1,\dots,n_{j-1}-n_j\,.
$$

\begin{lemma}[Markovian property]\label{l.mark}
Let $j\in\N$ and let $u_0,\ldots, u_j\in\Z^2$.
Suppose that the probability of the event
${(\gamma_0,\ldots,\gamma_j) = (u_0, \ldots, u_j)}$ 
is positive.
Conditioned on this event, the following holds.
\begin {enumerate}
\item
The paths
$\Gamma^1,\dots,\Gamma^j$ and $\Gamma [0, n_j]$ are conditionally independent.
\item
For $k\in\{1,\dots,j\}$, 
the conditional law of $\Gamma^{k}$ is that of a simple random 
walk in $D_{k-1}$ started from $u_k$ and conditioned to leave 
$D_{k-1}$ through the edge $[u_k, u_{k-1}]$.
\item
The conditional law of $\Gamma [ 0, n_j]$ is that of a simple 
random walk started from $0$ conditioned to leave $D_j$ at $u_j$,
and $\gamma[j,\ell]$ is the loop erasure of the time reversal
of $\Gamma[0,n_j]$. \qed
\end {enumerate}   
\end{lemma}

\proof
Since $\gamma$ is the loop-erasure of the reversal of
$\Gamma$, the event $(\gamma_0,\dots,\gamma_j)=(u_0,\dots,u_j)$
is equivalent to the statement that for
each $k=0,1,\dots,j-1$ the first hit of
$\Gamma$ to $\{u_0,\dots,u_k\}\cup \p D$
is through the edge $[u_{k+1},u_k]$.
Let $\tau_k:=\min\bigl\{n:\Gamma(n)\in \{u_0,\dots,u_k\}\cup\p D\bigr\}$,
$k=0,\dots,j$.
The strong Markov property of $\Gamma$ with 
the stopping times $\tau_k$ now implies the lemma. 
\QED

The following simple lemma will also be needed.

\begin{lemma}[Expected visits]\label{l.G1}
Suppose that $v\in\V D$ and that $u_0$ and $u_1$ are two
vertices
satisfying $\Pb{\gamma_0=u_0,\,\gamma_1=u_1}>0$.  Conditioned on
$\gamma_0=u_0$ and $\gamma_1=u_1$,
the expected number of visits to $v$ by $\Gamma^1$ is
$ G(u_1,v)\, H(v,u_1) $. 
\end{lemma} 

Here, $G(u,v)$ denotes the discrete Green's function;
that is, the expected number of visits to $v$ by
a simple random walk started at $u$, which is stopped on
exiting $D$.

\proof Let $X$ be simple random walk from $u_1$
stopped on exiting $D$ and let $k$ be the last time 
such that $X(k)=u_1$.  Then 
$\Gamma^1$ conditioned on $\gamma_0=u_0$ and $\gamma_1=u_1$
has the same distribution as $X$ conditioned
on $X(k+1)=u_0$.  But the path $j\mapsto X(k+j)$
is independent from $X[0,k]$.   Consequently,
the expected number of visits of $X$ to $v$ conditioned
on $X(k+1)=u_0$ is equal to the expected number of visits
to $v$ of $X[0,k]$.  The lemma follows.
\QED

\subsection{The core argument}\label{s.core}

We keep the previous notation and also use
the conformal maps
$
\psi_j:D_j\to\U
$ 
satisfying $\psi_j(0)=0$ and
$\psi_j'(0)>0$.
Set $U_j:= \psi_j ( \gamma_j)$ and $U:=U_0$.
Note that $\gamma$ can also be viewed as a continuously growing
simple curve from $\p D$ to $0$, and therefore can be represented
by Loewner's equation.
Let $W:[0,\infty)\to\p\U$ denote the (unique) continuous function such that
solving the radial Loewner equation with driving function $W(t)$ 
gives the path $\gamma$.
Note that $U_j = W ( \capa_j)$, where
$\capa_j$ is the continuous capacity
of $\gamma [0,j]$ from $0$ in $D$
(that is, the capacity of $\psi\bl(\gamma[0,j]\br)$ from $0$ in $\U$).
We denote by $\bl(\vartheta(t), t \ge 0\br)$ the continuous real-valued
function with $\vartheta(0) = 0$ such that 
$W(t) = W(0)\, \exp \bl( i \vartheta(t)\br)$. 
We also define $\Delta_j = \vartheta({t_j})$, so that 
$U_j = U\, \exp ( i \Delta_j)$.

\begin{proposition}[The key estimate]\label{p.bas}
There exists a positive constant $C$ such that for 
all small positive $\delta$, there exists $r_0= r_0 (\delta)$
such that the following holds.  
Let $D\in\Doms$ satisfy $\inrad(D)>r_0$.
For every  $u_0 \in \Vb D$ with $\Pb{\gamma_0=u_0}>0$,
let $\gamma$ denote the random path from $u_0$ to $0$
obtained by loop-erasure of the time reversal
of a simple random walk from  $0$ to $\p D$ conditioned
to hit $\p D$ in $u_0$.
Let
$$
m:=\min\bl\{j\ge 1:
\capa_j\ge\delta^2\hbox{ or }|\Delta_j|\ge\delta\br\}\,,
$$
where $\Delta_j$ and $\capa_j$ are as described above.
Then
\begin{equation}
\label{e.mart1}
\bl|\Eb{\Delta_m}\br| \le C \delta^3\,,
\end{equation}
and
\begin{equation}
\label{e.time1}
\bl|\Eb{\Delta_m^2}-2\,\Eb{\capa_m}\br| \le C\delta^3\,.
\end{equation}
\end{proposition}

Recall that Lemma~\ref{l.reverse} says that 
 $\gamma$ has the same
distribution as the {\em chronological} 
loop-erasure of random 
walk from $0$ to $\p D$ 
conditioned to hit $\p D$ at $u_0$.

Here is a rough sketch of the proof.
Let $v\in \V D$ 
satisfy
\begin{equation}\label{e.vrange}
\inrad(D)/200 <|v|<\inrad(D)/5\,.
\end{equation}
Let $h_0^+$ denote the number of visits to $v$ by $\Gamma$.
(This is the quantity which we referred to in the introduction
as the ``observable''.)
The proof is based on estimating 
the two sides of the equality
\begin{equation}\label{e.mot}
\Eb{h_0^+} = \EB{\Eb{h_0^+\md\gamma[0,m]}}\,.
\end{equation}
The estimate for the right-hand side will involve
the distribution of $t_m$ and $\Delta_m$.
We get the two relations~\eref{e.mart1} and~\eref{e.time1}
by considering two different choices for such a $v$.

The estimates for the two sides of~\eref{e.mot}
are rather straightforward.
Basically, each side is translated into expressions
involving the Green's functions $G_j$ and
the hitting probabilities $H_j$.  These are then
translated into analytic quantities using~\eref{e.phit}.
Earlier versions of the proof required other estimates, somewhat more delicate,
in addition to~\eref{e.phit}.  Fortunately, it turned out
that~\eref{e.phit} is sufficient.
Since we came across several different variants for
the proof, based on choosing different observables, 
it may be said that the proof is inevitable, rather than
accidental (and this also applies to 
Theorem~\ref{t.intropeano}).
Basically,  the reason the proof works is that
the expected number of visits to $v$ in
$\bigcup_{j=1}^m\Gamma^j$ given $\gamma[0,m]$ can be estimated
rather well given the rough geometry of $\gamma[0,m]$ in a scale
much coarser than the scale of the grid.
Similarly, it is important that $\Eb{h_0^+}$ can be estimated given
the rough-geometry of $D$, but this fact is not surprising.

\medskip
In the following,
we abbreviate the Green's function and hitting
probabilities in $D_j$ by $G_j:=G_{D_j}$ and $H_j:=H_{D_j}$.
The following lemma will be needed.

\begin{lemma}[Green's function bounds]\label{l.dG}
There is a constant $C>0$ such that for every $D\in\Doms$ and
$v\in\V D$ satisfying~\eref{e.vrange}
\begin{equation}
\label{e.lowG}
1/C \le G_D(0,v)\le C
\end{equation}
holds.  Also, given $\delta>0$ there is
an $r=r(\delta)$ such that if $\inrad(D)>r$,
then with the notations of Proposition~\ref{p.bas}
\begin{equation}
\label{e.dG}
G_0(0,v)-G_m(0,v)\le C\,\delta^2\,.
\end{equation}
\end{lemma}

\proofof{Proposition~\ref{p.bas}}
Since $\capa_{m-1}<\delta^2$ it follows from
the Koebe 1/4 theorem 
that $\inrad(D_{m-1})>\inrad(D)/5\ge r_0/5$ if $\delta$ is
small.
(Apply~\cite[Cor.~1.4]{\PommeBDRY} with
$z=0$ to $\psi_0^{-1}$ and $\psi_{m-1}^{-1}$.)
Moreover, the continuous harmonic measure in $D$ at   
$0$ of any edge $e$ with a 
vertex on $\p D$ can be made arbitrarily small
by requiring $\inrad(D)$ to be large.
(A Brownian motion started at $0$ has probability going
to $1$ to surround the disk $\inrad(D)\U$ before hitting $e$,
as $\inrad(D)\to\infty$.)
By conformal invariance of harmonic measure,
this implies that the diameter of $\psi(e)$ can be made
arbitrarily small.
 Applying this to the domains $D_j$ and 
using Lemma~\ref {l.diambound}, we see that 
we may take
$r_0$ large enough so that for all 
$j < m$, for all $t \in [t_j, t_{j+1}]$,  
$|\vartheta(t) - \vartheta({t_j})| \le \delta^{3}$ and
$|t_{j+1} - t_j| \le \delta^{3}$.
In particular,
$ t_m \le \delta^2 + \delta^{3}$ and 
$ |\Delta_m| \le \delta+ \delta^{3}$.
We also require $r_0/8$ to be larger than the $r(\delta)$
of Lemma~\ref{l.dG}.

Suppose $v\in \V D$ satisfies~\eref{e.vrange}.
Set  $\VV_j := \psi_j (v)$ and $\VV:=\VV_0$. 
For each $j \in \{1 , \ldots , \ell \}$,  let 
$h_j$ denote the number of visits to $v$ by 
$\Gamma^j$. Also let
$$
h_j^+:=\sum_{k=j+1}^\ell h_k\,,
$$
which is the number of visits of $v$ by $\Gamma [0, n_j]$.
Let
 $\lambda_j:= \lambda (v,\gamma_j; D_j)$,
where $\lambda(v,v'; D_j)$ is as in~\eref{e.lambdadef}. 
Since, conditionally on $\gamma[0,j]$,
 $\Gamma[0, n_j]$ is a random walk in $D_j$ conditioned to leave
$D_j$ at $\gamma_j$, 
$$
\Eb{  h_j^+ \md \gamma[0,j]  } =
\frac { G_j ( 0,v)\, H_j (v, \gamma_j) }{H_j(0, \gamma_j)}
$$
and Proposition \ref {p.phit} (together with~\eref{e.lowG}) implies  that if
 $r_0(\delta)$ is sufficiently large, for every
$j\in\{0,1,\dots,m\}$
$$ %
\Eb{  h_j^+ \md  \gamma[0,j]}
= G_j(0,v) \lambda_j + O (\delta^3)\,,
$$ %
(This $O(\cdot)$ notation is shorthand for the statement that there is an 
absolute constant $C$ such that
$\bl|\Eb{  h_j^+ \md  \gamma[0,j]} - G_j(0,v) \lambda_j\br|\le C\,\delta^3$.
We freely use this shorthand below.)
In particular
\begin {equation}
\label {e.diff}
\EB{ \sum_{j=1}^{m} h_j}
= \Eb{ h_0^+ - h_m^+}
=\Eb{  G_0 (0,v)\, \lambda_0 - G_m (0,v)\, \lambda_m }
+ O (\delta^3).
\end {equation}

We will now get a different approximation for the left-hand side.
Applying Lemma~\ref{l.G1} to the domain $D_{j-1}$ gives
$$
\Eb{h_j\md\gamma[0,j]} =
 G_{j-1}(\gamma_j,v)\, H_{j-1}(v,\gamma_j)\,.
$$
Proposition~\ref{p.phit} implies that for $r_0(\delta)$ large
enough,
\begin{equation}
\label{e.hj}
\Eb{h_j\md\gamma[0,j]} =
\bl(\lambda_{j-1}+O(\delta)\br)\,
 G_{j-1}(\gamma_j,v)\, H_{j-1}(0,\gamma_j)\,.
\end{equation}
Considering the same simple random walk starting at zero
and stopped when it exits $D_j$ or $D_{j-1}$ shows that
\begin{equation}
\label{e.tele}
G_{j-1}(0,v)-G_j(0,v)= H_{j-1}(0,\gamma_j)\, G_{j-1}(\gamma_j,v)
\,.
\end{equation}
  
We now derive an a priori 
bound on  $\max\bl\{\bl|\lambda_j - \lambda_m\br|:j\le m\br\}$.
Recall that 
\begin {equation}
\label {defl}
\lambda_j - \lambda_0
= \Re \left( \frac {U_j + \VV_j}{U_j - \VV_j} - \frac {U+\VV}{U-\VV} \right)
.\end {equation}
But $|U_j - U| \le O (\delta )$ for $j\le m$ and 
 Loewner's equation shows that  
\begin {equation}
\label {taylorV}
\forall j\le m\qquad
\VV_j  = \VV + t_j\,  \VV\,  \frac {U+\VV } { U-\VV} 
+ t_j\,O (\delta)  
 = \VV + t_j \,\VV \,\frac {U+\VV}{U-\VV}
  + O (\delta^3)
\,,
\end{equation}
and, in particular $\VV_j = \VV+ O (\delta^2)$.
(The equation blows up when $|U-\VV|$ is small, and such estimates
would not be valid in such a situation.  However, this is not
a problem here.  First, $\psi_0'(0)\le 1/\inrad(D)$, by
the Schwarz Lemma applied to the restriction of $\psi_0$ to
$\inrad(D)\U$. Now, the Koebe 1/4 theorem
(the case $z=0$ in the left hand inequality
in~\cite[Cor.~1.4]{\PommeBDRY}) gives  
$\psi_0^{-1}\bl((4/5)\U\br)\supset (1/4)\,|\psi_0'(0)|^{-1}\,(4/5) \U
\supset (\inrad(D)/5)\,\U$.
In particular, $|Z|=\bl|\psi_0(v)\br|\le 4/5$, by~\eref{e.vrange}.
Since $t_m=O(\delta^2)$, it is clear that if $\delta$ is small
and one starts flowing from $Z$ according to Loewner's equation,
it is impossible for $Z$ to get close to $\p \U$ up to time $t_m$.)
Thus, we get our bound,
$$ %
\forall j\le m\qquad
| \lambda_j - \lambda_m | \le O (\delta).
$$ %

Using~\eref{e.hj}, this implies
$$
\Eb{h_j\md\gamma[0,j]} =
\bl(\lambda_{m}+O(\delta)\br)\,
 G_{j-1}(\gamma_j,v)\, H_{j-1}(0,\gamma_j)\,.
$$
Now applying~\eref{e.tele} yields
$$
\EB{\sum_{j=1}^{m} h_j} = \EB{\bl(\lambda_{m}+O(\delta)\br)\,\bl(G_0(0,v)-G_{m}(0,v)\br)}
\,,
$$
and hence~\eref{e.dG} implies
$$
\EB{\sum_{j=1}^{m} h_j} = \EB{\lambda_m \,\bl(G_0(0,v)-G_m(0,v)\br)} + O(\delta^3)
\,.
$$
Comparing with~\eref{e.diff} gives
$ G_0(0,v)\,\Eb{\lambda_m-\lambda_0} = O(\delta^3) $, and hence~\eref{e.lowG}
implies
\begin{equation}
\label{e.core}
\Eb{\lambda_m-\lambda_0} = O(\delta^3)\,.
\end{equation}
(The reader may wonder about the apparent miracle happening here;
that $\lambda_j$ turns out to be ``almost'' a martingale.
In fact, this is not important for identifying the scaling limit.
If the right hand side in~\eref{e.core} turned out to be any other explicit
quantity, up to $\delta^3$ error terms,
the proof would still work, but give a different limiting process.
In Remark~\ref{r.lam} below, we give a short proof of~\eref{e.core} and 
further comments.)
  
Recall that this equation is valid uniformly over all choices
of $v$.
We now Taylor-expand $\lambda_m - \lambda_0$
with respect to $U_m - U$ and $\VV_m-\VV$,
 up to $O(\delta^3)$ error terms.
As we have seen, $U_m-U=O(\delta)$ and $\VV_m-\VV=O(\delta^2)$,
and hence only the 
first order derivative with respect to $\VV_m-\VV $ and 
the first two derivatives with respect to 
$$
U_m - U =  \bl(e^{ i\Delta_m} -1\br)\,U = i\,U\, \Delta_m - U\, \Delta_m^2 /2 
+ O(\delta^3)
$$
 come into play (the mixed derivatives can be ignored). 
Using \eref {defl} and \eref {taylorV} we get 
$$
\lambda_m - \lambda_0 
=  \Delta_m\, \Im\Bl(  \frac {2\,\VV\, U}{(U-\VV)^2} \Br)
+ (2\, t_m -\Delta_m^2)\, \Re\Bl( \frac { \VV\,U\, (U+\VV)}{(U-\VV)^3}\Br)
+ O (\delta^3)\,,
$$
and therefore~\eref{e.core} gives
\begin{equation}
\label{e.dlam}
\Im\Bl(  \frac {2\,\VV\, U}{(U-\VV)^2} \Br)\,
\Eb{  \Delta_m}
+
\Re\Bl( \frac { \VV\,U\, (U+\VV)}{(U-\VV)^3}\Br)
\,
\Eb{ 2\, t_m -\Delta_m^2}
= O (\delta^3)\,.
\end{equation}

We claim that when $r_0(\delta)$ is large enough, we may find
$v_1,v_2\in \V D$ in the range~\eref{e.vrange}
satisfying $\bl|\psi(v_1)-U/30\br|<\delta^3$
and $\bl|\psi(v_2)- i\,U/30\br|<\delta^3$.
Indeed, by Theorem~1.3 and Corollary~1.4 from \cite{\PommeBDRY},  
for every $R\in(0,1)$ there is a $C=C(R)<\infty$
such that $|\psi'(z)|\le C/\inr z(D)$ and
$\inr z(D)\ge C^{-1}\,\inrad(D)$  hold for all $z\in\psi^{-1}(R\,\U)$.
Let $v_1$ be a vertex closest to $\psi^{-1}(U/30)$.
By integrating the above bound on $\psi'$ along the line segment
from $\psi^{-1}(U/30)$ to $v_1$ (whose length is less than $1$),
we get $\bl|\psi(v_1)-U/30\br|<\delta^3$, 
if $\inrad(D)$ is large enough.
Another application of
Theorem~1.3 and Corollary~1.4 from \cite{\PommeBDRY} now show  
that $v_1$ satisfies~\eref{e.vrange}.  An entirely similar
argument produces $v_2$.

Consequently,~\eref{e.dlam} holds with $\VV\in\{U/30,i\,U/30\}$.
Plugging in these two values for $\VV$ produces two linearly
independent equations in the variables 
$ \Eb{ 2\, t_m -\Delta_m^2}$ and $\Eb{\Delta_m}$, 
and thereby proves~\eref{e.mart1}
and~\eref{e.time1}.
\QED

\begin{remark}\label{r.lam}
Here is another proof of~\eref{e.core}.
Given a vertex $v\in \V D$, let
$\beta=(\beta_0,\beta_1,\dots)$ denote the loop-erasure
of the reversal of the simple random walk $\Gamma_v$ started from $v$ and
stopped on exiting $D$ (i.e., the analogue of $\gamma$,
but starting from $v$ instead of $0$).
Abbreviate $\gamma^n:=(\gamma_0,\dots,\gamma_n)$,
and similarly $\beta^n:=(\beta_0,\dots,\beta_n)$.
For a sequence of vertices $u=(u_0,u_1,\dots,u_n)$,
let $a_n(u):=\Ps{\gamma^n=u}$ and $b_n(u):=\Ps{\beta^n=u}$.
Set $ M_n:={b_n(\gamma^n)}/{a_n(\gamma^n)}$.
(In other words, $M_n$ is the Radon-Nikodym derivative of
the law of $\beta^n$ with respect to the law of 
$\gamma^n$.)
It is easy to verify that $M_n$ is a martingale:
$$
\Eb{M_{n+1}\md \gamma^n} =
\sum_w
\frac{ b_{n+1}(\gamma^nw)}
{ a_{n+1}(\gamma^nw)}
\frac
{ a_{n+1}(\gamma^nw)}
{ a_{n}(\gamma^n)}
=
\frac{ \sum_w b_{n+1}(\gamma^nw)}
{ a_{n}(\gamma^n)}
=M_n\,.
$$
Lemma~\ref{l.mark} implies that $M_n=H_n(v,\gamma_n)/H_n(0,\gamma_n)$,
since, on the event that $\Gamma_v$ and $\Gamma$ first hit
$\{u_0,\dots,u_n\}\cup\p D$ at $u_n$, we may couple them to
agree after that first visit to $u_n$.
Now~\eref{e.phit} implies~\eref{e.core}.

Although this proof is shorter than the first proof of~\eref{e.core},
it is harder to motivate and less natural.  For this reason,
we chose to stress the first proof.

Let us finally note that (as opposed to the martingale that shows up in 
the analysis of the UST Peano curve),
the quantity corresponding to this martingale 
in the scaling limit is unbounded and converges almost surely
to zero (it is not uniformly integrable), so that it can not 
be interpreted as  
a  conditional probability or a conditional expectation.
Correspondingly, in the discrete setting, $M_n$ is very large 
when the path hits $v$ (if it does) and $M_n$ is very small
when the path hits $0$.
\end{remark}

\subsection{Recognizing the driving process}\label{s.recog}

The objective in this subsection is to show that
$W$ of the previous section
 is close to a time-scaled Brownian motion on the unit circle.

\begin{theorem}[Driving process convergence]
\label{t.sleconv}
For every $T>0$ and $\eps>0$ there is an
$r_1=r_1(\eps,T)>0$ such that for all $D\in\Doms$ with
$\inrad(D)>r_1$ there is a coupling of $\gamma$ with Brownian motion
$B(t)$ starting at a random uniform point in $[0,2\pi]$
such that
$$
\PB{\sup\bl\{|\vartheta (t) - B(2\,t)|:t\in[0,T]\br\}>\eps}<\eps\,.
$$
\end{theorem}

Recall that a coupling of two random variables
(or random processes) $A$ and $B$
is a probability space with two random variables
$A'$ and $B'$, where $A'$ has the same distribution
as $A$ and $B'$ has the same distribution as $B$.
In the above statement (as is customary) we don't distinguish
between $A$ and $A'$ and between $B$ and $B'$.

In order to deduce this theorem from Proposition \ref {p.bas},
we will use
the Skorokhod Embedding Theorem, which is one of
the standard tools for proving convergence to Brownian motion
(one could work out a more direct proof but the following proof
seems cleaner).

\begin{lemma}[Skorokhod embedding]
\label {l.sko}
If $(M_n)_{n \le N}$ is an $({\cal F}_n)_{n \le N}$ martingale,
with  $\| M_n-M_{n-1}\|_\infty\le 2\, \delta$ and $M_0=0$ a.s.,
 then there are stopping times
$0=\tau_0\le\tau_1\le \cdots \le\tau_N$  for standard Brownian motion 
$(B_t, t \ge 0)$, such that
$(M_0,M_1, \ldots, M_N)$ and $(B_{\tau_0},B_{\tau_1}, \ldots, B_{\tau_N})$ 
have the same law.
Moreover, one can impose for $n=0,1,\dots,N-1$
\begin{equation}\label{e.etau}
\Eb{  \tau_{n+1} - \tau_n  \md B[0, \tau_n] }
= \Eb { (B_{\tau_{n+1}} - B_{\tau_n} )^2 \md B [0, \tau_n]}
\end{equation}
and
\begin {equation}
\label {e.bdt}
\tau_{n+1} \le \inf \bl\{ t \ge \tau_n \ : \ | B_t - B_{\tau_n}| \ge 
2\, \delta\br \}.
\end {equation}
\end {lemma}
The proof can be found in many 
probability textbooks including \cite {\Dudley, \RevuzYor}. Often, it is stated for just one random 
variable $M_1$; for a statement in terms of martingales see,
for instance, \cite {\Strassen,\Dubins}.
The relation (\ref {e.bdt}) is not stated explicitly in these
references (since the assumption that the increments of $M_n$
are bounded is weakened), but is a consequence of the proof.
It can also be derived a posteriori
from $\Es{\tau_{n+1}-\tau_n}=\Eb{(M_{n+1}-M_n)^2}<\infty$,
since the expected time for Brownian motion started outside
an interval to hit the interval is infinite.

\proofof{Theorem~\ref{t.sleconv}}
Since the hitting measure of simple random walk from zero is
close to the hitting measure for Brownian motion when 
$\inrad(D)$ is large
(see, e.g., Section~\ref{s.RW}), it is clear that $W(0)$ is nearly
uniform in $\p\U$.  It is therefore
enough to show that $\vartheta (t/2)$ is close to
standard Brownian motion.

Assume, with no loss of generality, that $T\ge 1$.
Pick $\delta=\delta(\eps,T)>0$ small.
Let $r_0$ be as in Proposition~\ref{p.bas}
and take $r_1:=8\exp(20\, T)\,r_0$.
Let $\gamma^t$ denote the initial segment of $\gamma$ 
such that $\psi_D(\gamma^t)$ has capacity $t$ from $0$.
By the Schwarz Lemma $\psi_D'(0)\le \inrad(D)^{-1}$.
Therefore, the Koebe 1/4 Theorem implies 
$\inrad(D\setminus \gamma^t)\ge \exp(-t)\,\inrad(D)/4$.
Hence, if $\inrad (D) \ge r_1$, 
Proposition~\ref{p.bas} is valid not only for the initial
domain $D$, but also for the domain $D$ slitted by
subarcs of $\gamma$, up to capacity $20\,T$.

As in Proposition~\ref{p.bas},
define $m$ to be the first $j=1,2,\dots$ such that
$|\Delta_j|\ge\delta$ or $t_j\ge\delta^2$.
Set $m_0:=0$,  $m_1:=m$, and inductively let $m_{n+1}$ be the first
$j \ge m_n+1$ such that
$|\Delta_j-\Delta_{m_n}|\ge\delta$ or
$t_j-t_{m_n}\ge\delta^2$, whichever happens first.
Let $\ev F_n$ denote the $\sigma$-field generated by $\gamma[0,m_n]$.
Set
$$N := \lceil 10\, T \,\delta^{-2}\rceil.$$
Our choice of $r_1$ ensures that 
$t_{j+1} - t_j \le 2\, \delta^2$ for all $j<N$,
and that $t_N \le 20\, T$. Hence, 
Proposition~\ref{p.bas} holds for all domains $D_{m_n}$ with $n <N$.
Applying clause 3 of Lemma~\ref{l.mark} therefore gives
\begin{equation}\label{e.mart}
\Eb{\Delta_{m_{n+1}}-\Delta_{m_n}\md \ev F_n}=O(\delta^3)\,,
\end{equation}
and
\begin{equation}
\label{e.time}
\Eb{(\Delta_{m_{n+1}}-\Delta_{m_n})^2\md\ev F_n}=
2\,\Eb{\capa_{m_{n+1}}-\capa_{m_n}
\md\ev F_n}+O(\delta^3)\,.
\end{equation}

For $n\le N$, set
\begin {equation}
\label {e.M}
M_n 
:=
\sum_{j=0}^{n-1}
\left( \Delta_{m_{j+1}}-\Delta_{m_{j}}
-\Eb{\Delta_{m_{j+1}}-\Delta_{m_{j}}\md \ev F_{j}} \right) \,. 
\end {equation}
Clearly, $M_0, \ldots, M_N$ is a martingale for ${\cal F}_0, \ldots,
{\cal F}_N$. The definition of $m_n$ and the choice of $r_1$ imply that 
$\|M_{n+1}-M_n\|_\infty\le 2\,\delta$.

By Lemma~\ref{l.sko}, we may couple $(M_0,\dots,M_N)$ with
a standard Brownian motion with stopping times
$\tau_0\le\tau_1\le\cdots\le\tau_N$
such that $B_{\tau_n}=M_n$ and~\eref{e.etau} hold.
Extend the coupling to include $\gamma$ (this clearly
can be done).

Note that the definition of $t_{m_n}$ and~\eref{e.bdt} ensure that for 
all $n < N$,
\begin{equation}\label{e.intermed}
\begin{split}
&
\sup \{ |B_t - B_{\tau_n}| \ : \  t \in [\tau_{n}, \tau_{n+1}] \}   
\le 2\, \delta\,, 
\\
& 
\sup \{ |\vartheta (t) -  \Delta_{t_{m_n}} | \ : \ t \in [t_{m_n}, t_{m_{n+1}}] \}
\le 2\, \delta
\end{split}
\end{equation}
and (\ref {e.mart}) shows that 
\begin{equation}
\label{e.MD}
\sup \{ | \Delta_{t_{m_n}}  - M_n | \ : \ n \le N \} 
= O( \delta^3N)  = O (\delta\,T)\,.
\end{equation}
Hence, as $M_n = B_{\tau_n}$ and $B_t$ is a.s.\ continuous,
it remains to relate the 
capacities $t_{m_n}$ with the stopping times $\tau_n$ and
verify that $t_{m_N}>T$ with high probability. 
For this purpose, define 
$$ Y_n = \sum_{j=0}^{n-1} ( M_{j+1} - M_j )^2 .
$$
We first show that $Y_n$ is close to $2\,t_{m_n}$.
Let $Z_n:=Y_n-2\,t_{m_n}$.
By~\eref{e.M} and~\eref{e.mart}, we have for $n<N$,
$|M_{n+1}-M_n-\Delta_{t_{m_{n+1}}}+\Delta_{t_{m_{n}}}|=O(\delta^3)$.
This implies $|M_{n+1}-M_n| = O(\delta)$ and hence also
$$
Y_{n+1}-Y_n=(M_{n+1}-M_n)^2= (\Delta_{t_{m_{n+1}}}-\Delta_{t_{m_{n}}})^2
+ O(\delta^4)\,.
$$
Consequently,~\eref{e.time} gives
$$
\Eb { Z_{n+1}-Z_n \md {\cal F}_n } \le  O (\delta^3)  \,.
$$
{}From the fact that the
increments of $t_{m_n}$ and those of $Y_n$ are bounded by $O(\delta^2)$,
we also have $
\Eb {  \bl(Z_{n+1}-Z_n\br)^2 \md {\cal F}_n }  \le  O (\delta^4)$.
Set
$
Z_n':=
Z_n-\sum_{j=1}^n\Eb{Z_{j}-Z_{j-1}\md \ev F_{j-1}}
$.
Since this is an $\ev F_n$-martingale, we
have
$ \Eb{{Z_N'}^2}=\sum_{j=1}^N \Eb{(Z_j'-Z'_{j-1})^2}$
and the above estimates give
$ \Eb{{Z_N'}^2} =O(N\,\delta^4) $.
Assuming $N\,\delta^3<\delta^{1/2}/2$, without loss of generality,
and applying Doob's maximal inequality~\cite[II.1.7]{\RevuzYor}
for $L^2$ martingales to $Z_n'$, we get
\begin {equation}
\label {e.doob1}
\PB{  \max_{n \le N} \bl|Y_n-2\,t_{m_n} \br|> \delta^{1/2} }
 = O (N\,\delta^{3}) = O (T\,\delta).
\end {equation}
By the definition of the $t_{m_n}$, we have
$Y_{n+1}-Y_n+t_{m_{n+1}}-t_{m_n}\ge \delta^2$.
Summing gives $Y_N+t_{m_N}\ge N\,\delta^2\ge 10\,T$.
Therefore,~\eref{e.doob1} implies
\begin{equation}\label{e.Nok}
\Pb{t_{m_N}<2\,T}=O(T\,\delta)
\,.
\end{equation}

We now show that with high probability $\tau_n$ is also close to $Y_n$
for every $n\le N$.
By~\eref{e.bdt}, it is clear that
$\Eb{(\tau_{n+1} - \tau_n)^2\md B[0,\tau_n]}=O(\delta^4)$, and
therefore
$$
\Eb { ((\tau_{n+1} - Y_{n+1}) -
( \tau_n - Y_ n ))^2 \md B[0, \tau_n] }  =  O (\delta^4)
\,.
$$
Also,~\eref{e.etau} gives
$$
\Eb { (\tau_{n+1} - Y_{n+1}) -
( \tau_n - Y_ n ) \mid B[0, \tau_n] } =0\,.
$$
Doob's inequality therefore implies
$$ %
\PB{ \max_{n \le N} \bl|\tau_n - Y_n \br|>\delta^{1/2} } =O(T\,\delta)\,.
$$ %
Combining this with (\ref {e.doob1}) leads to 
$$
\PB{ \max_{n \le N} \bl|\tau_n - 2\,t_{m_n} \br|>\delta^{1/2} } =O(T\,\delta)\,.
$$
Since $B_t$ is a.s.\ continuous,
together with~\eref{e.Nok},~\eref{e.intermed}
and~\eref{e.MD}, this completes the proof.
 \QED

\subsection{Convergence with respect to a stronger topology}\label{s.improve}

Theorem~\ref{t.sleconv} provides a kind of convergence of loop-erased
random walk to \SLEtwo/. As we will see in
the present subsection, this kind of convergence suffices,
for example, to show that the scaling limit with respect
to the Hausdorff metric of the union of $\p\U$ and LERW in $\U$ 
is the union of $\p\U$ and the \SLEtwo/ path.

Let $\alpha:[0,1]\to \C$ and $\beta:[0,1]\to \C$
be two continuous paths.  Define
$$
\pdist(\alpha,\beta):=
\inf_{\phi\in\Phi}
\sup_{t\in[0,1]} |\alpha(t)-\beta\circ\phi(t)|
\,,
$$
where $\Phi$ is the collection of all monotone non-decreasing
continuous maps from
$[0,1]$ onto $[0,1]$.  It is an easy well-known fact that
$\rho$ is a metric on equivalence classes of paths, where
two paths $\alpha$ and $\beta$ are equivalent if
$\alpha\circ \phi_1=\beta\circ \phi_2$, where
$\phi_1,\phi_2\in\Phi$.
Since $\pdist(\alpha,\beta)$ does not depend on the
particular parameterization of $\alpha$ or $\beta$,
the metric $\pdist$ is also defined for paths
on intervals other than $[0,1]$.

To explain our present goal, let us point out
that there is a sequence of paths $\alpha_n$ from
$1$ to $0$ in $\closure\U$ such that 
their Loewner driving functions $W_n(t)$
converge uniformly to the constant $1$ but
$\alpha_n$ does not converge to the path
$\alpha(t)=1-t$, $t\in[0,1]$, in the metric $\pdist$,
although the driving function for $\alpha$
(reparameterized by capacity) is the constant $1$.
For example, we may take $\alpha_n$ as 
the polygonal path through the points
$a_1,b_1+i\,n^{-2},a_2,b_2-i\,n^{-2},a_3,b_3+i\,n^{-2},
\dots,a_{\lfloor n/2\rfloor},0$,
where $a_j:= 1-n^{-1}+(jn)^{-1}$ and
$b_j:=1-j/n$.

\begin{theorem}[LERW image in $\closure\U$ converges]\label{t.str}
For any sequence $D_n \in \Doms$
with $\inrad (D_n) \to \infty$, if $\mu_n$ 
denotes the law of $\tilde\gamma^n:=\psi_{D_n} \circ \gamma^n$,
where $\gamma^n$ is the time-reversal of LERW from
$0$ to $\partial D_n$, then $\mu_n$ converges
weakly (with respect to the metric $\rho$)
to the law of the radial \SLEtwo/ path started uniformly on the 
unit circle.
\end{theorem}

The outline of the proof goes as follows. We define a 
suitable family of compact subsets of the space of simple 
paths from $\partial \U$ to $0$ in $\closure\U$, which we can
use to show that the sequence $\mu_n$ is tight.
(See, e.g.,~\cite{\Dudley} for background on weak convergence
and the notion of tightness.)
This implies that a subsequence of $\mu_n$ converges weakly to some
probability measure.
Theorem~\ref {t.str}
then shows that the law of \SLEtwo/ is the unique possible subsequential limit.

In order to prove tightness, we will
use properties of loop-erased random walk proved in~\cite{\SchSLE}.
The actual details will require some background in the geometric theory
of conformal maps.  In particular,
some properties of extremal distance (a.k.a.\ extremal
length) will be used.  See, for example,~\cite{\Ahlfors} for background.
The basic ideas that are 
used in the proof are taken from
\cite {\AizenmanBurchard} and \cite {\SchSLE}.

\medskip

For a simply connected $D\subsetneqq\C$
containing $0$, let $\XX_0(D)$ denote the space of
all simple paths
$\gamma:[0,\infty]\to\closure{D}$
from $\p D$ to $0$ in $\closure{D}$, which intersect
$\p D$ only at the starting point.
Given a monotone nondecreasing
function $\gauge:(0,\infty)\to (0,1]$,
let $\XX_\gauge(D)\subset\XX_0(D)$
denote the space of all simple paths
$\gamma\in\XX_0(D)$ such that
for every $0\le s_1<s_2$,
$$
\dist(\gamma[0,s_1]\cup\p D,\gamma[s_2,\infty])/\inrad(D)\ge
\gauge\bl(\diam(\gamma[s_1,s_2])/\inrad(D)\br)\,.
$$
Note that whether or not $\gamma\in\XX_\gauge(D)$ does not depend
on the parameterization of $\gamma$, and is scaling invariant.

\begin{lemma}[Compactness]\label{l.iscomp}
Let $\gauge:(0,\infty)\to(0,1]$ be monotone nondecreasing.
Then $\XX_\gauge(\U)$ is compact in the topology of
convergence with respect to $\pdist$.
\end{lemma}

\proof
We use an idea from~\cite{\AizenmanBurchard}.
For all $n\in \N$ let $Z_n$ be a finite collection
of points such that the open balls $\ball(z,2^{-n})$,
$z\in Z_n$, cover $\U$.  Given a set $K\subset\closure\U$
and a point $z\in Z_n$, let
$s(K,z,n)$ denote the diameter
of $K\cap \ball(z,2^{1-n})$ in the
metric obtained from the Euclidean metric
on the disk $\closure{\ball(z,2^{1-n})}$ by collapsing
the boundary $\p \ball(z,2^{1-n})$ to a single point.
(In other words, this metric $d(x,y)$ is defined
as $d(x,y)=\min\{|x-y|,\dist(x,\p \ball)+\dist(y,\p \ball)\}$,
where $\ball=\ball(z,2^{1-n})$.)

Fix $\gamma\in\XX_\gauge(\U)$.  Given $t\ge 0$, let
$$
s(t)=s_\gamma(t):=
\sum_{n\in\N}\sum_{z\in Z_n}
\frac
{s(\gamma[0,t],z,n)}
{|Z_n|}
\,.
$$
Clearly, 
$s_\gamma(t) \le \sum_{n \ge 0} 2^{2-n} = 8$,
 and $s:[0,\infty]\to [0,\infty)$ is continuous
and strictly monotone increasing.
(To verify that $s$ is strictly monotone increasing,
note that if $t_2>t_1\ge 0$, then there 
is some $n\in\N$ such that
$\dist(\gamma(t_2),\gamma[0,t_1])\ge 2^{2-n}$,
and so $s(\gamma[0,t_2],z,n)\ge s(\gamma[0,t_1],z,n)+2^{-n}$
if $z\in Z_n$ satisfies $\gamma(t_2)\in \ball(z,2^{-n})$.)
Let $\hat \gamma(s)$ be $\gamma$ parameterized by
$s$; that is, $\hat\gamma=\gamma\circ s^{-1}$.
Let $s_1<s_2$ and set $\eps:=\diam\hat\gamma[s_1,s_2]>0$.
Then
$\dist\bl(\hat\gamma(s_2),\hat\gamma[0,s_1]\cup\p\U\br)\ge \gauge(\eps)$.
By the argument for strict monotonicity given above,
this shows that $s_2-s_1\ge 2^{-n}/|Z_n|$, where
$n:=\min\{k\in\N: 2^{2-k}\le\gauge(\eps)\}$.
  Therefore, $\hat\gamma$ satisfies an equicontinuity
estimate.  By the Arzela-Ascoli Theorem, it follows that
the closure of $\XX_\gauge(\U)$ is compact in the $\pdist$ metric.
It is also clear that $\XX_\gauge(\U)$ is closed.
\QED

Our next goal is to use these compact sets to prove tightness,
and start by observing that the diameter is tight.

\begin{lemma}[Diameter is tight]\label{l.tightdiam}
There are constants $c,C>0$ such that for every
$D\in\Doms$  and every $r\ge 1$ the simple random walk
$\Gamma$ starting from $0$ and stopped on hitting $\p D$
satisfies
$$
\Pb{\diam(\Gamma)\ge r\,\inrad(D)}\le C\,r^{-c}\,.
$$
Consequently, the same estimate holds for the loop-erasure $\gamma$.
\end{lemma}

The first statement is an easy well-known fact.  Since the complement
of $D$ is connected and unbounded, if the random walk makes a loop
separating the circle $\inrad(D)\,\p\U$ from the circle
$(r/2)\,\inrad(D)\,\p\U$ before hitting the latter circle, then it must hit
$\p D$ before $(r/2)\,\inrad(D)\,\p\U$. 
Thus, the lemma is easily proved directly,
and also follows from the convergence of simple random walk to Brownian motion.
A rather precise form of this estimate for the random walk, where $c=1/2$,
is known as the discrete Beurling theorem~\cite[Theorem 2.5.2]{\Lbook}.

\begin{lemma}[Tameness]\label{l.gtight}
For every $\epsilon>0$ there is some monotone nondecreasing
$\gauge:(0,\infty)\to(0,1]$
and some $r_0>0$ such that for every
$D\in\Doms$ with $\inrad(D)\ge r_0$ its time-reversed loop-erased
walk $\gamma=\gamma_D$ satisfies
$$
\Pb{\gamma\in\XX_\gauge(\Dd)}\ge 1-\epsilon\,.
$$
\end{lemma}

\proof The proof is essentially contained in
the proof of~\cite[Thm.~1.1]{\SchSLE},
where it is established that every subsequential scaling limit of
LERW is a.s.\ a simple path.
We will not repeat the complete proof from~\cite{\SchSLE} here, but
indicate how it may be adapted to yield the statement of the lemma.

Let $\eps>0$.
Clearly, $\gamma\in \XX_0(D)$.  If $\gamma\notin\XX_\gauge(D)$,
then there are $0\le s_1 < s_2<\infty$
such that the distance between
$\gamma[0,s_1]\cup\p D$ and $\gamma[s_2,\infty]$ is
smaller than $\inrad(D)\,\gauge\bl(\diam\gamma[s_1,s_2]/\inrad(D)\br)$.
Let us first deal with the case where the distance between
$\gamma[s_2,\infty]$ and $\p D$ is small.
Let $\Gamma$ be the walk generating the time-reversal of $\gamma$,
and let $t_n$ be the first time $t$ where the distance
from $\Gamma(t)$ to $\p D$ is smaller than $2^{-n}\,\inrad(D)$,
and let $\tau=\inf\{t:\Gamma(t)\in\p D\}$.
By the Markov property of $\Gamma$ at time
$t_n$ and Lemma~\ref{l.tightdiam}, 
$$
\Pb{\diam\Gamma[t_n,\tau]> 2^{-n/2}\,\inrad(D)}\le C\, 2^{-c\,n/2}\,.
$$
Consequently, there is an $N=N(\eps)$
such that with probability $1-\epsilon/2$
for every integer $n\ge N$ we have
$\diam\Gamma[t_n,\tau]\le 2^{-n/2}\,\inrad(D)$.
In this case, if $\diam\gamma[0,s_2]> 2^{-n/2}\,\inrad(D)$,
where $n>N$, then $\gamma[0,s_2]$ is not contained
in $\Gamma[t_n,\tau]$, which implies that
$\gamma[s_2,\infty]\subset\Gamma[0,t_n]$, and gives
$\dist(\gamma[s_2,\infty],\p D)\ge 2^{-n}\,\inrad(D)$.
In other words, if $\gauge$ satisfies
\begin{equation}\label{e.gau1}
\gauge(t)<\min\{t^2,2^{-2\,N}\}/4\,,
\end{equation}
then with probability at least $1-\eps/2$,
for every $s_1,s_2\in[0,\infty]$,
\begin{equation}\label{e.nearbd}
\dist(\p D,\gamma[s_2,\infty])\ge
\inrad(D)\,\gauge (\diam\gamma[s_1,s_2]/\inrad(D))\,.
\end{equation}

We now focus on the case where the distance between
$\gamma[0, s_1]$ and $\gamma [s_2, \infty)$ is small.
We shall say that $\gamma$ has a $(\beta,\alpha)$-quasi-loop
if there are $0<s_1<s_2<\infty$ such that
$|\gamma(s_1)-\gamma(s_2)|\le  \alpha\,\inrad(D)$ but
$\diam\gamma[s_1,s_2]\ge \beta\,\inrad(D)$.
Note that if there are $0<s_1<s_2<\infty$ such that
$\dist(\gamma[0,s_1],\gamma[s_2,\infty])< \alpha\,\inrad(D)$
and $\diam\gamma[s_1,s_2]\ge\beta\,\inrad(D)$, then
$\gamma$ has a $(\beta,\alpha)$-quasi-loop.
Let $\ev A(\beta,\alpha)$ denote the event that
$\gamma$ has a $(\beta,\alpha)$-quasi-loop.
Assume, for the moment, that for all $n\ge 0$,
\begin{equation}\label{e.Alim}
\lim_{\alpha\searrow0}\Ps{\ev A(2^{-n},\alpha)}=0\,,
\end{equation}
uniformly in $D$.
Then we may take a decreasing sequence $\alpha_n\searrow 0$
such that $\sum_{n=1}^\infty \Ps{\ev A(2^{-n},\alpha_n)}<\eps/2$
holds for every $D\in\Doms$.
Then with probability at least $1-\eps/2$, $\gamma$
has no $(2^{-n},\alpha_n)$-quasi-loop for any $n=1,2,\dots$.
Assuming that $\gauge(t)<\alpha_n$ holds whenever
$t\le 2^{1-n}$, $n\in\N$, and $\gauge(t)<\alpha_1$ for all $t$,
on this event we also have
$$
\dist(\gamma[0,s_1],\gamma[s_2,\infty])\le
\inrad(D)\,\gauge(\diam\gamma[s_1,s_2]/\inrad(D))
$$
for all $0<s_1<s_2<\infty$.
If we also assume~\eref{e.gau1},
then together with~\eref{e.nearbd}
we get
$ \Pb{\gamma\in\XX_\gauge(\Dd)}\ge 1-\epsilon$, completing the proof of
the lemma.  Thus, it remains to verify~\eref{e.Alim}.

Let $\ev A(z_0,\beta,\alpha)$ denote the event that
there are $0<s_1<s_2<\infty$ such that
$|\gamma(s_1)-\gamma(s_2)|\le \alpha\,\inrad(D)$,
$\gamma(s_1),\gamma(s_2)\in \ball(z_0,\beta\,\inrad(D)/4)$
and $\diam(\gamma[s_1,s_2])\ge \beta\,\inrad(D)$.
In particular, this implies that $\gamma[s_1,s_2]$ is
not contained in the interior of $\ball(z_0,\beta\,\inrad(D)/2)$.
Assume that $8\,\alpha<\beta$.
By Lemma~\ref{l.tightdiam}, there is an $R=R(\eps)>0$ such that
with probability at least $1-\eps/2$ we have
$\gamma[0,\infty]\subset \ball(0,R\,\inrad(D))$.
There is a collection $\{z_1,z_2,\dots,z_k\}$ of points
such that every disk of radius $2\,\alpha\,\inrad(D)$
with center in $\ball(0,R\,\inrad(D))$ is contained in one of the $k$ balls
$\ball(z_j,\beta\,\inrad(D)/2)$, $j=1,2,\dots,k$,
and we may take $k< c((R/\beta)^2+1)$, where $c$ is an
absolute constant.
On the event
$\gamma[0,\infty]\subset \ball(0,R\,\inrad(D))$ 
we have $\ev A(\beta,\alpha)\subset \bigcup_{j=1}^k\ev A(z_j,\beta,\alpha)$.
Since $\eps>0$ was arbitrary and
$\Pb{\gamma[0,\infty]\subset \ball(0,R\,\inrad(D))}\ge 1-\eps/2$,
it is therefore sufficient to show that
$\Pb{ \ev A(z_j,\beta,\alpha)}\to 0$ as $\alpha\to 0$,
uniformly in $D$.  The proof of this statement
is given (with minor changes in the setup)
in~\cite[Theorem 1.1]{\SchSLE}.
\QED

Let $\XX_\gauge^r(D)$ denote the set of paths
$\gamma\in\XX_\gauge(D)$ that are contained in the ball
of radius $r\,\inrad(D)$ about $0$.  Given $\gamma\in \XX_0(D)$,
let $\gamma^*:[0,\infty)\to\closure{\U}$
denote the path $\psi_D\circ\gamma$, parameterized by capacity.

\begin{lemma}[Tameness invariance]\label{l.transfer}
For every monotone nondecreasing $\gauge:(0,\infty)\to(0,1]$
and every $r>1$
there is a monotone nondecreasing $\gauge^*:(0,\infty)\to(0,1]$
such that for all $D\in\Doms$ and $\gamma \in \XX_\gauge^r (D)$,
$
\gamma^* \in \XX_{\gauge^*}(\U)
$.
\end{lemma}

\proof
Let $D\in\Doms$, $\gamma\in\XX^r_\gauge(D)$ and $0\le s_1'<s_2' \le\infty$.
Note that there exist $s_1$ and $s_2$ satisfying
$s_1' \le s_1 \le s_2 \le s_2'$ such that 
$$
\diam \bl( \gamma^* [s_1, s_2] \br) \ge \diam \bl( \gamma^* [s_1', s_2']\br) / 4
$$
and
\begin {equation}
\label {e.s0}
\dist \bl( 0 , \gamma^* [ s_1, s_2] \br) \ge
\diam \bl( \gamma^* [s_1, s_2] \br).
\end {equation}
Since
\begin {equation}
\label {e.s1}
\dist \bl( \gamma^* [0, s_1'] \cup \partial \U , \gamma^* [s_2', \infty)\br)
\ge 
\dist \bl( \gamma^* [0, s_1] \cup \partial \U , \gamma^* [s_2, \infty)\br)
,\end {equation}
it is sufficient to give a lower bound of the right-hand side of 
(\ref {e.s1}) in terms of $\eps:= \diam (\gamma^* [s_1, s_2])$. 

The Schwarz Lemma gives $\psi_D'(0)\le 1/\inrad(D)$.
Therefore,
by the Koebe 1/4 theorem 
(applied to the restriction of $\psi_D^{-1}$ to $\eps\U$)
and~\eref{e.s0}, 
$\dist\bl(0,\gamma[s_1,s_2]\br)>c_1\,\inrad(D)$, where
$c_1=\eps/4$.
On the other hand, the harmonic measure in $\U$ from $0$ of
$\gamma^*[s_1,s_2]$ is at least $c_2$, where $c_2=c_2(\eps)>0$; so 
that the harmonic measure in $D$  from $0$ of
$\gamma[s_1,s_2]$ is at least $c_2$.  
Hence, 
\begin{equation}
\label{e.dia}
\diam\gamma[s_1,s_2]\ge c_3\,\inrad(D)
\,,
\end{equation}
where $c_3=c_3(\eps)$.

Also set
$\delta:=\dist\bl(\gamma^*[0,s_1]\cup\p\U,\gamma^*[s_2,\infty]\br)$.
Since
$$
\diam\gamma^*[s_2,\infty]\ge\dist\bl(0,\gamma^*[s_1,s_2]\br)\ge
\eps\,,
$$
the extremal distance between
$\gamma^*[0,s_1]\cup\p\U$ and $\gamma^*[s_2,\infty]$
is at most $\phi_1(\delta,\eps)>0$, where $\phi_1$
is some function satisfying $\phi_1(\delta,\eps)\to 0$ as
$\delta\downarrow 0$.
By conformal invariance of extremal distance,
this implies that the extremal distance
between $\gamma[0,s_1]\cup\p D$ and
$\gamma[s_2,\infty]$ is at most $\phi_1(\delta,\eps)$.
Because $\gamma$ is contained in the disk of radius
$r\,\inrad(D)$ about $0$, this implies that
$$
\dist\bl( \gamma[0,s_1]\cup\p D, \gamma[s_2,\infty]\br)\le
\phi_2(\delta,\eps)\,r\,\inrad(D)\,,
$$
where $\phi_2\to 0$ as $\delta\downarrow 0$.
Because $\gamma\in\XX_\gauge(D)$,~\eref{e.dia} and this together imply
$$
\phi_2(\delta,\eps)\,r
\ge \gauge\bl( c_3(\eps)\br)\,,
$$
which gives a positive lower bound for
$\delta=\dist\bl(\gamma^*[0,s_1]\cup\p\U,\gamma^*[s_2,\infty]\br)$
in terms of $\gauge$, $r$ and
$\eps=\diam\bl(\gamma^*[s_1,s_2]\br)$.
This completes the proof.
\QED

\begin{lemma}[Convergence relations]\label{l.lconv} 
Suppose $W^n,W$  are continuous functions
from $[0,\infty)$ to $\p\U$ such that 
$W^n \rightarrow W$ locally uniformly.  Let
$g^n_t,g_t$ be the corresponding solutions to
Loewner's radial equation and set $f_t^n = (g_t^n)^{-1}$,
$f_t = g_t^{-1}$. 
Then $f_t^n \rightarrow f_t$ locally uniformly on $[0,\infty) \times\U$.
If there are continuous curves $\gamma^n:[0,\infty) \to \closure\U$
such that for all $t \geq 0$, the image of $f_t^n$
is the component of $0$ in $\U \setminus \gamma^n[0, t]$
and there is a $\gamma:[0,\infty)\to\closure\U$
such that $\gamma^n\to\gamma$ locally uniformly on $[0,\infty)$,
then for all $t\ge 0$ the image of $f_t$ is the component of $0$
in $\U \setminus \gamma[0,t]$.
\end{lemma}
\proof
Since $g_t$ is obtained by flowing along a vector field depending on $W$,
the inverse $f_t$ is obtained by flowing along the opposite field, with
the time reversed.  Hence, the first statement is an immediate
consequence of the principle that solutions of ODE depend continuously
on the parameters of the ODE.
The second statement is an immediate consequence of the 
Carath\'eodory kernel theorem~\cite[Theorem 1.8]{\PommeBDRY}.  
\QED

\proofof{Theorem \ref{t.str}}
Let $W^n$ denote the Loewner parameter of $\tilde\gamma^n$
and let $\hat \mu_n$ denote the law of the pair $(\tilde\gamma^n,W^n)$.
By Theorem~\ref{t.sleconv}, we know that the law of $W^n$ tends
weakly to the law of Brownian motion.
The lemmas~\ref{l.iscomp}--\ref{l.transfer}
show that the set of measures $\{\mu_n\}$
is tight with respect to the metric $\pdist$.  
Consequently, the sequence $\hat\mu_n$ is also tight.
Prokhorov's theorem (e.g., \cite {\Dudley,\RevuzYor}) implies that
there is a subsequence
such that $\hat\mu_{n}$ converges weakly along the subsequence.
Let $\hat\mu$ be any subsequential weak limit,
and let $(\tilde \gamma,W)$ be a sample from $\hat\mu$.
The lemmas show that $\tilde \gamma$ is a.s.\ a simple path
and Theorem~\ref{t.sleconv} shows that $W$ is Brownian
motion (with time scaled).
By properties of weak convergence,
we may couple the subsequence of pairs $(\tilde\gamma^n,W^n)$ and $(\tilde \gamma,W)$
so that a.s.\ $\pdist(\tilde\gamma^n,\tilde \gamma)\to 0$
and $W^n\to W$ locally uniformly.

Recall that the capacity is continuous with respect to the metric
$\pdist$; that is, if $\beta,\beta_n:[0,1]\to\overline\U\setminus\{0\}$
and $\pdist(\beta_n,\beta)\to 0$, then the capacity of $\beta_n[0,1]$
tends to the capacity of $\beta[0,1]$.  (In fact, it is enough that
$\beta_n[0,1]$ tends to $\beta[0,1]$ in the Hausdorff metric.)
Indeed, this follows immediately from
Carath\'eodory's kernel theorem~\cite[Theorem 1.8]{\PommeBDRY},  
and the fact that local uniform convergence of conformal maps
implies the convergence of the derivatives (by Cauchy's formula for
the derivative).

Since $\tilde \gamma$ is almost surely a simple path, the capacity
of $\tilde\gamma$ increases strictly, and one can parametrize the 
path continuously by its capacity. 
We also parameterize the paths $\tilde\gamma^n$ by capacity.
The next goal is to show that $\tilde\gamma^n\to\tilde\gamma$,
locally uniformly on $[0,\infty)$.
Since $\pdist(\tilde\gamma,\tilde\gamma^n)\to 0$, there are strictly
monotone continuous onto maps $\varepsilon_n:[0,\infty)\to[0,\infty)$
so that $\tilde\gamma^n\circ\varepsilon_n \to\tilde\gamma$
locally uniformly.  If $t_n\in[0,\infty)$ and $t_n\to t\in[0,\infty)$,
then it follows from the continuity of capacity with respect to
$\rho$ that $\varepsilon_n(t_n)\to t$ (because if
$s$ is a subsequential limit of $\varepsilon_n(t_n)$,
then the capacity of $\tilde\gamma(s)$ must be $t$; that
is $s=t$).  This implies that
$\varepsilon_n$ converges to the identity map $t\mapsto t$,
locally uniformly.  By continuity of $\tilde\gamma$,
it follows that $\tilde\gamma\circ \varepsilon_n^{-1}\to\tilde\gamma$
locally uniformly.  This gives $\tilde\gamma^n\to\tilde\gamma$
locally uniformly.

We can now finally apply Lemma~\ref{l.lconv} to 
see that $\tilde \gamma$ is the \SLEtwo/ path.
As the law of the limit $\tilde \gamma$ does not depend on the subsequence,
the theorem follows.
\QED

In the following proof of Theorem~\ref{t.lerw}, the main
technical point is that we do not make any smoothness assumptions
on $\p D$.  If $\p D$ is a simple closed path, the Theorem
follows easily from Theorem~\ref{t.str}, because
the suitably normalized conformal maps from $\U$ to the
discrete approximations of $D$ converge uniformly to
the conformal map onto $D$.

\proofof{Theorem~\ref{t.lerw}}
Let $D_\delta$ be the component of $0$ in the
complement of all the closed square faces of the
grid $\delta\Z^2$ intersecting $\p D$.
Let $\gamma_\delta$ be the time reversal of the loop-erased random
walk from $0$ to $\p D_\delta$,
and let $\beta$ be the radial \SLEtwo/ path in $\closure\U$.
Let $\phi_\delta:\U\to D_\delta$ be the conformal map
satisfying $\phi_\delta(0)=0$ and $\phi'_\delta(0)>0$,
and let $\phi:\U\to D$ be the conformal
map satisfying $\phi(0)=0$, $\phi'(0)>0$.
Theorem~\ref{t.str} tells us that
we may couple $\beta$ with each of the paths
$\gamma_\delta$ such that
$\pdist(\phi_\delta^{-1}\circ \gamma_\delta,\beta)\to 0$
in probability as $\delta\downarrow 0$. 
Moreover, the proof shows that if we use the
capacity parameterization for both, then in probability
$$
\sup\bl\{|\phi_\delta^{-1}\circ\gamma_\delta(t)-
\beta(t)|: t\ge 0\br\}\to 0\,.
$$
(There is no problem with convergence in a neighborhood
of $t=\infty$, because we know that the weak limit of
$\phi_\delta^{-1}\circ\gamma_\delta$ with respect to
$\pdist$ is a simple path tending to $0$ as $t\to\infty$.)

The Carath\'eodory kernel theorem~\cite[Theorem 1.8]{\PommeBDRY}.  
implies that 
$\phi_\delta\to\phi$ uniformly  on
compact subsets of $\U$ as $\delta\searrow0$.  Consequently, the above gives 
\begin{equation}
\label{e.gloop}
\forall t_0>0\qquad
\sup\bl\{|\gamma_\delta(t)-
\phi\circ\beta(t)|: t\ge t_0\br\}\to 0\,,
\end{equation}
in probability.  Let $\eps>0$ be small.  
Then, by Lemma~\ref{l.tightdiam}, there is an $\eps'>0$ such that for
every $D'\in\Doms$ the probability
that simple random walk from $0$ gets to
distance $\inrad(D')/(2\eps')$ before hitting $\p D'$ is
less than $\eps/2$.
Let $A$ be the connected component of $0$
in the set of points in $D\cap (\inrad(D)/\eps')\U$
having distance at least $\eps\,\eps'\,\inrad(D)$ from $\p D$.
By considering the first point where the random walk
generating $\gamma_\delta$ exits $A$, it follows
that with probability at least $1-\eps$, the
diameter of $\gamma_\delta[0,\infty]\setminus A$
is at most $\eps\,\inrad(D)+\delta$.
Now note that there is a compact $A'\subset\U$ such
that $\phi_\delta^{-1}(A)\subset A'$ for all
sufficiently small $\delta$, since $\phi_\delta\to\phi$ uniformly on
compacts.  Therefore, there is
some $t_1>0$ such that  
$\gamma_\delta[0,t_1]\cap A=\emptyset$ a.s.\ for all
sufficiently small $\delta>0$.
In particular, 
\begin{equation}
\label{e.gs}
  \Pb{\diam \gamma_\delta[0,t_1]>\eps\,\inrad(D)+\delta}<\eps\,.
\end{equation}
If we take $t_2\in(0,t_1)$, then taking $\delta\searrow0$ in~\eref{e.gloop}
implies 
$$
\Pb{\diam\phi\circ \beta[t_2,t_1]>2\,\eps\,\inrad(D)}<\eps\,.
$$
Since this holds for every $t_2$, it follows that
$$
\Pb{\diam\phi\circ \beta(0,t_1]>2\,\eps\,\inrad(D)}<\eps\,.
$$
Using this with~\eref{e.gs} and choosing $t_0=t_1$ in~\eref{e.gloop}
gives
$$
\PB{
\sup\bl\{|\gamma_\delta(t)-
\phi\circ\beta(t)|: t>0\br\}<3\,\eps\,\inrad(D)}\to 1\,.
$$
Since this holds for every $\eps>0$, the theorem follows.
\QED

\section{The UST Peano curve}
\label{s.peano}
\subsection{Setup}
\label{s.peanosetup}

\begin{figure}
\SetLabels
\T\R(0.02*0.02)$a$\\
\B\L(.83*.95)$b$\\
\L(.87*.295)Peano\\
\L(.87*.52)tree\\
\L(.87*.75)dual tree\\
\endSetLabels
\centerline{\AffixLabels{\includegraphics*[height=2.3in]{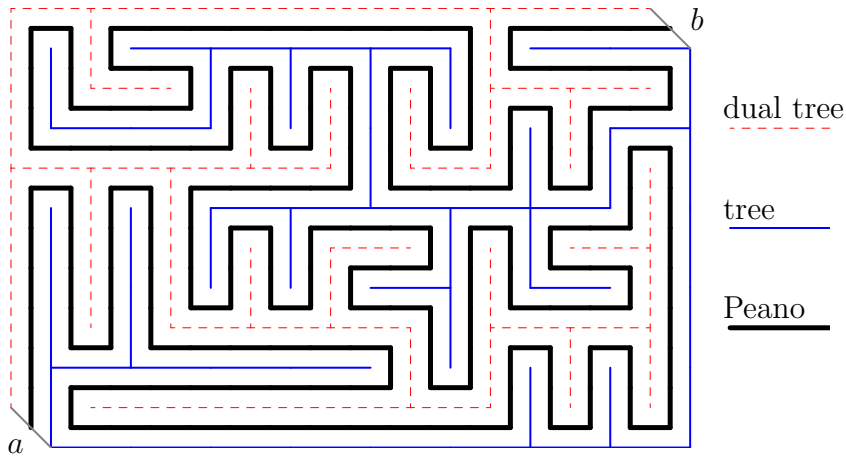}}}
\caption{\label{f.peanoc} The tree, dual tree, and Peano UST path $\gamma$.} 
\end{figure}

The UST Peano curve is obtained as the interface between the
UST and the dual UST.  The setup which corresponds 
 to chordal \SLEkk8/
is where there is symmetry between the UST and the dual UST.
 Loosely speaking, the UST is the uniform 
spanning tree on the grid inside a
domain $D$ but with an entire arc $\alpha\subset\p D$ on the boundary
identified (wired) as a single vertex,
and the dual UST also has an arc $\beta\subset\p D$ on the boundary
which is identified.  The arcs $\alpha$ and $\beta$ are essentially
complementary arcs.  See Figure~\ref{f.peanoc}, where $D$ is approximately
a rectangle.
As mentioned in the introduction, it was conjectured \cite {\RSsle}
that for an analogous setup, the 
interface defined for the critical random cluster models 
with $q \in (0,4]$ converges to \SLEkk\kappa/, where
$\kappa=\kappa(q) \in [4,8)$.
  
A combinatorial framework is necessary in order to be more precise.
There are several different possible setups that would work,
and the following is somewhat arbitrary.

If a tree $T$ lies in the
grid $\Z^2$, then its dual tree $T^\dagger$ will lie
in the dual grid $(\Z+1/2)^2$,
and the Peano path $\gamma$ will lie
in the graph $G$ whose vertices are $(1/4+\Z/2)^2$
and where $v,u$ neighbor iff $|v-u|=1/2$.
We have three kinds of vertices: elements of $\Z^2$ are
the {\em primal} vertices, elements of $(1/2+\Z)^2$ are
the {\em dual} vertices, and elements of 
$(1/4+\Z/2)^2$ are the {\em Peano} vertices.
If $w\ne v$ are vertices of any kind, not necessarily
the same, we say that they are {\em adjacent} if
the distance between them is as small as it
can be for distinct vertices of these particular kinds.
In other words, 
if they are of the same kind, this means that they are neighbors,
if $v\in (1/4+\Z/2)^2$ and $w\in \Z^2\cup (1/2+\Z)^2$,
this means $\|v-w\|_\infty=1/4$, while if
$v\in \Z^2$ and $w\in (1/2+\Z)^2$, this means
$\|v-w\|_\infty=1/2$.

\begin{figure}
\SetLabels
(0.69*0.29)$a$\\
(.31*.89)$b$\\
\L(.86*.6)$\alpha$\\
\R(.14*.6)$\beta$\\
\endSetLabels
\centerline{\AffixLabels{\includegraphics*[height=2.3in]{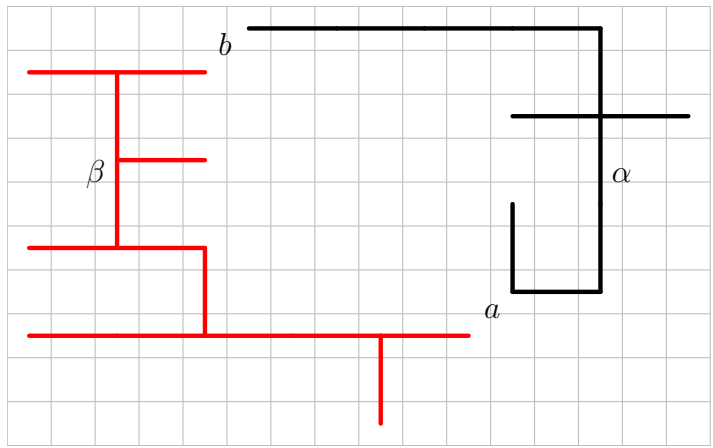}}}
\caption{\label{f.peanotreebd}The boundary data and the Peano grid.}
\end{figure}

Since there is no added complication, we consider a more
general case where $\alpha$ and $\beta$ are 
trees, rather than arcs.
Let $\alpha$ be some finite tree in the primal grid
$\Z^2$ and let $\beta$ be a finite tree in the dual
grid $(1/2+\Z)^2$.  Suppose that no edge of $\alpha$
intersects an edge of $\beta$.
Further suppose that there are two Peano vertices
$a,b\in(1/4+\Z/2)^2$ such that
$a$ is adjacent to both a primal vertex
$\alpha_a\in\alpha$ and a dual vertex
$\beta_a\in\beta$, and $b$ is adjacent to both a primal
vertex $\alpha_b\in\alpha$ and a dual vertex
$\beta_b\in\beta$.  See Fig.~\ref{f.peanotreebd}.
Note that the line segment $[\alpha_a,\beta_a]$
has $a$ as its midpoint, and the line segment
$[\alpha_b,\beta_b]$ has $b$ as its midpoint.
Let $D=D(\alpha,\beta,a,b)$ be the (unique) bounded connected
component of
$\C\setminus\bl( %
\alpha\cup [\alpha_b,\beta_b]\cup\beta\cup[\beta_a,\alpha_a]\br)$. %
Let $V_P=V_P(D)$ denote the collection of all Peano vertices
in $\closure D$, and, as before, $\V D$ denotes the collection
of all primal vertices in $D$.
Let $\ell=\ell(D)$ denote the cardinality of $V_P\setminus \{a,b\}$.
By switching the role of $a$ and $b$, if necessary,
assume that $D$ lies to the immediate right
of the oriented segment $[\alpha_a,\beta_a]$.
Let $\Doms^*$ denote the collection of all domains obtained in this
way.

Let $H=H(D)$ denote the subgraph of $\Z^2$
whose vertices are the vertices of $\alpha$ and $\V D$,
and whose edges are those edges
on this set of vertices which do not intersect
$\beta$.
Since $\beta$ is a tree, $H$ is connected.

Since $H$ is connected, there is at least one spanning
tree $T$ of $H$ which contains $\alpha$.
If we replace $\alpha$ by $T$ and apply the dual argument,
it follows that there is also a tree $T^\dagger$
in the dual grid $(1/2+\Z)^2$, which is disjoint
from $T$, contains $\beta$, and whose vertices
are the dual vertices in $\beta$ and the dual
vertices in $D$.  In fact, $T^\dagger$ contains %
every dual edge lying in $\closure D\setminus T$. %

We now need to give
an orientation to the Peano grid $G$.
Every edge in $G$ is either on the boundary
of a square face of $G$ centered on a primal
vertex, or is on the boundary of a square
face of $G$ centered on a dual vertex,
and these two possibilities are exclusive.
We orient the edges of $G$ by specifying that the square
faces of $G$ containing a primal vertex are oriented clockwise,
while those containing a dual vertex are oriented counterclockwise.
When we want to emphasize the orientation of the edges,
we write $\GO$ instead of $G$.
Note that the edges of $G$ contained in a horizontal 
or vertical line all get the same direction in $\GO$,
and consecutive parallel lines get opposite orientations.
For this reason, $\GO$ is often called the {\em Manhattan lattice}.

Let $\gamma=\gamma(T)$ denote the set of all edges
of $\GO$ which
do not intersect $T\cup T^\dagger$ and
which have at least one endpoint in $D$.
Let $v\in V_P\setminus \{a,b\}$ be some  Peano vertex in $D$.
Note that there are precisely two oriented edges of $\GO$
with initial point $v$, say $e_1$ and $e_2$, where
one of these, say $e_1$, intersects an edge $f_1$
of the primal grid $\Z^2$, and the other intersects
an edge $f_2$ of the dual grid $(1/2+\Z)^2$.
Note also that $f_1\cap f_2\ne\emptyset$.
It therefore follows that exactly one of the edges
$f_1,f_2$ is in $T\cup T^\dagger$.
Consequently, exactly one of the edges
$e_1,e_2$ is in $\gamma$.  This shows that
$\gamma$ has out-degree $1$ at every $v\in V_P\setminus\{a,b\}$.
  An entirely similar argument
shows that $\gamma$ has in-degree $1$ at every such $v$.
In particular, this shows that $\gamma$
does not contain the entire boundary of a square face
of $G$ that does not contain a primal or dual vertex.
If $\gamma$ had a cycle, the cycle therefore
would have to surround some primal or dual vertex.
But as $T$ and $T^\dagger$ are connected
and disjoint from $\gamma$, this is impossible.
It therefore follows that $\gamma$ is an
oriented simple path (i.e., self avoiding path),
and the endpoints of $\gamma$ are $a$
and $b$.
Since we are assuming that $D$ lies to
the right of $[\alpha_a,\beta_a]$,
the initial point of $\gamma$ is
$a$ and the terminal point is $b$.

Conversely, suppose that
$\gamma^*=(\gamma^*_0,\dots,\gamma^*_{\ell+1})$ is
any oriented simple path in $\GO$,
respecting the orientation of $\GO$, from $a$
to $b$, whose vertices are $V_P$.
For $n\in\{0,\dots,\ell+1\}$, 
let $v_n$ be the (unique) primal vertex
adjacent to $\gamma^*_n$, and let
$v_n^\dagger$ be the dual vertex
adjacent to $\gamma^*_n$.
Note that $v_n$ and $v_{n+1}$
are either the same vertex, or are adjacent
vertices when $n=\{0,\dots,\ell\}$.
Let $\alpha_n=\alpha_n(\gamma^*)$ denote the union
of $\alpha$ with the collection of
all edges $[v_k,v_{k+1}]$ for
$k<n$ such that $v_k\ne v_{k+1}$,
and similarly, let $\beta=\beta_n(\gamma^*)$
denote the union of $\beta$ with the
collection of all dual edges
$[v^\dagger_k,v^\dagger_{k+1}]$ for
$k<n$ such that $v^\dagger_k\ne v^\dagger_{k+1}$.
Then $T(\gamma^*):=\alpha_{\ell+1}(\gamma^*)$ and
$T^\dagger(\gamma^*):=\beta_{\ell+1}(\gamma^*)$ are obviously connected,
and there are no edges in $T(\gamma^*)$ intersecting edges in
$T^\dagger(\gamma^*)$.  Now, $T(\gamma^*)$ cannot
contain a cycle, for such a cycle would have to separate
$T^\dagger(\gamma^*)$.  Hence, $T(\gamma^*)$
is a spanning tree of $H$ containing $\alpha$.
It is also clear that $\gamma^*=\gamma\bl(T(\gamma^*)\br)$.
That is, $T\mapsto \gamma(T)$ is a bijection
between the set of spanning trees of $H$ containing
$\alpha$ and the set of oriented paths in $\GO\cap\closure{D}$ from $a$ to
$b$ containing $V_P$.
Hence, when $T$ is the UST on $H$ conditioned to
contain $\alpha$, $\gamma$ is uniformly distributed
among such Peano paths; it is the {\em UST Peano path}
associated with $(\alpha,\beta,a,b)$.

\medskip

Let $(a=w_0,w_1,\dots,w_{\ell+1}=b)$ be the order of
the vertices in the UST Peano path $\gamma$.
For $n\in\{0,1,\dots,\ell\}$ let
$\gamma[0,n]$ denote the initial arc of $\gamma$ from $w_0$ to $w_n$.
Since $\gamma$ is uniformly distributed among simple 
oriented paths in $\GO$ from $a$ to $b$ which contain $V_P$, %
 we immediately get the following Markov %
property.

\begin{lemma}[Markovian property]\label{l.pmark}
Fix any $n\in\{1,2,\dots,\ell\}$.
Conditioned on $\gamma[0,n]$, the distribution of
$(\gamma\setminus\gamma[0,n])\cup\{w_n\}$ is the same as that
of the UST Peano curve associated with
$\bl(\alpha_n(\gamma),\beta_n(\gamma),w_n,b\br)$.
 \QED
\end{lemma}
 
This lemma will play the same role
 in the proof as Lemma \ref {l.mark} in the case of LERW.
We will also use the convergence of  certain 
discrete harmonic functions towards their continuous 
counterparts.  To facilitate this, we have to 
set the combinatorial notation for the discrete Dirichlet-Neumann  %
problem.

Let $H$ be a finite nonempty  %
connected subgraph of $\Z^2$
with vertices $V_H$,
and let $E_\p$ denote the set of oriented edges in $\Z^2$
whose initial endpoint is in $V_H$, but whose unoriented
version is not in $H$.  Suppose $E_\p=E_0\cup E_1\cup E_2$
is a disjoint union, where $E_0\cup E_1\ne\emptyset$.
Suppose also that $\hat h:V_H\to [0,1]$ is some function.
For $v\in V_H$ set $\Delta_{H,E_0,E_1,E_2} \hat h(v):= \sum d\hat h[v,u]$,
where the sum is over all neighbors $u$ of $v$
in $\Z^2$, and $d\hat h[v,u]:=\hat h(u)-\hat h(v)$ when $[v,u]\notin E_\p$,
$d\hat h[v,u]:=0-\hat h(v)$ when $[v,u]\in E_0$,
$d\hat h[v,u]:=1-\hat h(v)$ when $[v,u]\in E_1$ and
$d\hat h[v,u]:= 0$ when $[v,u]\in E_2$.
Note that there is a unique $\hat h: V_H\to [0,1]$
such that $\Delta_{H,E_0,E_1,E_2} \hat h(v)=0$
in $V_H$: $\hat h(v)$ is the
probability that a simple random walk on $H\cup E_0\cup E_1$
started from $v$ will use an edge in $E_1$ before using
an edge of $E_0$. 
This $\hat h$ will be called the $\Delta_{H,E_0,E_1,E_2}$-harmonic function.

\begin{proposition}[Dirichlet-Neumann approximation]\label{p.free}
For every $\eps>0$ there is an $r_0=r_0(\eps)$ such
that the following holds.
Let $D\subset\C$ be a simply connected domain
satisfying $\inrad(D)\ge r_0$.
Let $A_0,A_1\subset\p\U$ be two disjoint arcs, each of length at least $\eps$,
and set $A_2:=\p\U\setminus(A_0\cup A_1)$.
Let $\eta\subset D$ be a simple closed path
which surrounds $0$, such that each point of $\eta$ is within
distance $5$ from $\p D$.
Suppose that $A_0',A_1'\subset\eta$ are two disjoint arcs,
$A_2':=\eta\setminus (A_0'\cup A_1')$, and the triple $(A_0',A_1',A_2')$
corresponds to $(A_0,A_1,A_2)$ under $\psi_D$, in the sense that
for each $j=0,1,2$ and each $p\in A_j'$ there is a continuous
path $\sigma:[0,1)\to D$ satisfying
$\diam \sigma[0,1)\le 5$,
$\sigma(0)=p$, and $\lim_{s\uparrow 1} \psi_D\circ\sigma(s)$
exists and is in $A_j$.

Let $H$ be the component of $0$ in the set of edges
of $\Z^2$ that do not intersect $\eta$.
For $j=0,1,2$, let $E_j$ denote the set of oriented edges $[v,u]$
intersecting $\eta$,
where $v$ is in $H$, and the first point of intersection from
the direction of $v$ is in $A'_j$.
Let
$\hat h$ denote the $\Delta_{H,E_0,E_1,E_2}$-harmonic
function.
Let $h:\U\to[0,1]$ be the continuous harmonic function
which has boundary value $0$ on $A_0$,
$1$ on $A_1$, and satisfies the
Neumann boundary condition on $A_2$.
Then $\bl|\hat h(0)-h(0)\br|<\eps$.
\end{proposition}

The proof will be given in Section~\ref{s.pfreepf}. %

\subsection{Driving process convergence}\label{s.peanodrive}

\begin{figure}
\SetLabels
(.975*.04)$a$\\
(.04*.96)$b$\\
\R\T(.60*.51)$\gamma(n)$\\
\R(.66*.15)$\beta_n$\\
\B\R(.77*.7)$\alpha_n$\\
\endSetLabels
\centerline{\AffixLabels{\includegraphics*[height=2.3in]{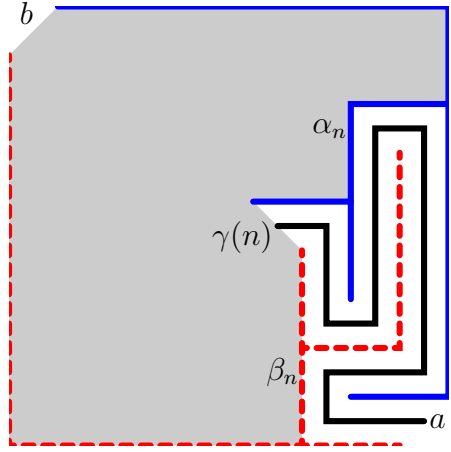}}}
\caption{\label{f.domains}The domain $\tilde D_n$ is shaded.}
\end{figure}

Let $\alpha,\beta,a,b$ and $D=D(\alpha,\beta,a,b)$ be as
above, and suppose now that $0\in D$.
As before, let $\ell$ denote the number of Peano vertices
in $D$, and let $\gamma=\bl(\gamma(0),\dots,\gamma({\ell+1})\br)$
be the UST Peano path from $a$ to $b$ in $\GO\cap\closure D$.
For each $n \le \ell$, there are two domains 
that are naturally associated to $\gamma [0,n]$.
The first one (as in Lemma~\ref {l.pmark}) is 
$\tilde D_n := D( \alpha_n, \beta_n , \gamma(n), b)$
(see figure~\ref {f.domains}).
But $\tilde D_n$ is not so useful if we   %
want to make estimates using Loewner's equation.  %
We therefore also define
$D_n:=D\setminus\gamma[0,n]$.
Let $\phi_0=\phi:D\to\H$ be the conformal map
which takes $D$ to $\H$, takes $b$  to
$\infty$, takes $a$  to $0$, and 
satisfies $|\phi(0)|=1$.
Let $\phi_n:D_n\to\H$ be the conformal maps
satisfying $\phi_n(z)-\phi_0(z)\to 0$
as $z\to b$ within $D_n$.
Define $W_n:=\phi_n\bl(\gamma(n)\br)\in\R$. 
Also let $t_n$ denote the capacity from infinity in $\H$
of $\phi_0 \circ \gamma[0,n]$, so that
$\phi_n\circ \phi^{-1}_0(z)= z + 2\, t_n / z + o(1/z)$
when $z \to \infty$ in $\H$.

We now prove the analog of Proposition~\ref{p.bas} for the
UST Peano curve.  Let $\harm_D(z,A)$ denote the
continuous harmonic measure of $A$ from $z$ in
the domain $D\setminus A$.

\begin{proposition}[The key estimate]\label{p.peanoLoc}
For every sufficiently small $\delta,\eps>0$ there is some
$r_0=r_0(\delta,\eps)$ such that the following holds.  %
Let $\gamma,D_n,\phi_n,W_n$ and $t_n$ be as above,
let $k\in\N$, and let $m$ be the first $n\ge k$
such that
$|W_n-W_k|\ge \delta$ or $t_n-t_k\ge\delta^2$.
Then
\begin{equation}
\label{e.martp}
\Eb{W_m-W_k\md D_k}=O(\delta^3)\,,
\end{equation}
and
\begin{equation}
\label{e.timep}
\Eb{(W_m-W_k)^2 \md D_k}=8\,\Eb{t_m-t_k\md D_k}+O(\delta^3)\,,
\end{equation}
provided that $\inrad(D_k)\ge r_0$
and $\harm_{D_k}\bl(0,\alpha_k(\gamma)\br)\in[\eps,1-\eps]$.
\end{proposition}

\proof
Assume first $k=0$.
Let $v_0\in \V D$ be some vertex 
such that $|\phi(v_0)|\le 2$
and $\Im \phi(v_0)\ge 1/2$, say.
(As we have seen in Section~\ref{s.core}, there is such a
$v_0$ when $\inrad(D)$ is large.)
If $Q=[q,q']$ is a line segment where
$q\in D$ is a dual vertex and $q'\in\alpha$ is the
midpoint of a dual edge containing $q$, %
then let $\phi^*(Q)\in \R_+$ denote the limit
of $\phi(z)$ as $z$ tends to $\p D$ along $Q$
(which always exists, by~\cite[Proposition 2.14]{\PommeBDRY}).  
Fix such a $Q_0$ satisfying
$U:=\phi^*(Q_0)\in [1/2,2]$; there clearly is such $Q_0$
when $r_0$ is large, because the harmonic
measure from $0$ of any square of the dual grid adjacent to
the boundary of $D$ is small.
Let $\eta\subset D$ be the set of points within distance $1/10$ from $\p D$.
Then $\eta$ is a simple closed path.
Consider it as oriented counterclockwise around the bounded domain
of $\C\setminus\eta$.
Let $p_0$ be the point of $\eta$ closest to $a$,
$p_1$ the point in $\eta\cap Q_0$ and $p_2$ the point of
$\eta$ closest to $b$. 
Let $A_0'$ be the positively oriented subarc of
$\eta$ from $p_0$ and $p_1$, $A_1'$ the positively oriented
arc from $p_1$ to $p_2$, and $A_2'$ the positively oriented
arc from $p_2$ to $p_0$.

Let $\ev A$ be the event that the path in the tree
$T(\gamma)=\alpha_{\ell+1}(\gamma)$ from $v_0$
to $\alpha$ hits $A'_1$.
We will now estimate both sides of the identity 
\begin{equation}\label{e.pea}
\Pb{\ev A}=\Eb{\Ps{\ev A\md  D_m}}
\end{equation}
using Proposition~\ref{p.free}.
By Wilson's algorithm, $\Pb{\ev A}$ is the
probability that
a simple random walk on the graph $H( D)$ started at $v_0$ 
stopped on hitting $\alpha$
will cross $A_1'$.  This is exactly $\hat h (v_0)$,
where the function $\hat h$ is
as defined in Proposition~\ref {p.free}.
Set 
$$
 h(z):=
\frac 1{\pi}
\cot^{-1}\Bl(\frac{1-|z|}{2\,\Im\sqrt z}\Br)
=
\frac 1{\pi}
\cot^{-1}\Bl(\frac{1-r}{2\,\sqrt r\,\sin(\theta/2)}\Br)\,,
$$
where $z=r\,e^{i\theta}$ and we take the value of $\cot^{-1}$
between $0$ and $\pi$.
Note that $h$ is harmonic in $\H$, is equal to $0$
on $(0,1)$, is equal to $1$ on $(1,\infty)$, and
$\p_y h=0$ on $(-\infty,0)$.  
(Of course, we found the map $h$ satisfying these
boundary conditions by reflecting the
domain along the negative real axis, mapping this larger domain
to $\H$ with $z\mapsto\sqrt z$, and then using a conformal
map from $\H$ to $\U$ to calculate the hitting probabilities.)
Consequently, Proposition~\ref{p.free} shows that
if $r_0$ is sufficiently large, then
\begin{equation}\label{e.pevA}
\Ps{\ev A} =  h\bl(\phi(v_0)/U\br) + O(\delta^3)\,.
\end{equation}
Set $V_j:=\phi_j(v_0)$ and $U_j:=\phi_j\circ \phi_0^{-1}(U)$.
By the chordal version of Loewner's equation and the definition of %
$m$, we have %
\begin{equation}\label{e.apx}
V_m=V_0+\frac {2\,t_m}{V_0}+O(\delta^3)\,,
\qquad
U_m=U_0+\frac {2\,t_m}{U_0}+O(\delta^3)\,.
\end{equation} 
Note also that $\inrad(D_m)>\inrad(D)/2$, $v_0\in D_m$,
and $U_m\in [1/4,4]$ provided that $\delta$ is small enough.

We now employ a similar argument to estimate $\Pb{\ev A\md D_m}$.
Recall that $\tilde D_n =D\bl(\alpha_n,\beta_n,\gamma(n),b\br)$.
Assume that $Q_0$ intersects $\tilde D_n$, which will
be the case if $U_n>1/4$, say.
Let $\eta_n$ be the set of points in $\tilde D_n$ at distance
$1/10$ from $\p\tilde D_n$.
Again $\eta_n$ is a simple closed path, and we 
write $\eta_n=A_0'(n)\cup A_1'(n)\cup A_2'(n)$,
where $A_0'(n)$ is the arc of $\eta_n$ from the closest
point to $\gamma(n)$ to the point of intersection
of $Q_0$ with $\eta_n$, $A_1'(n)$ is the arc of $\eta_n$
from the point in $Q_0\cap\eta_n$ to the point of 
$\eta_n$ closest to $b$, and let $A_2'(n)$ be the remaining
part of $\eta_n$.
By Lemma~\ref{l.pmark}, $\Pb{\ev A\md D_n}$ is
the same as the quantity $\hat h_n(v_0)$, where $\hat h_n$
is the function $\hat h$ defined in Proposition~\ref{p.free},
but with $A_0'(n),A_1'(n),A_2'(n)$ and
$\eta_n$ replacing $A_0',A_1',A_2'$ and $\eta$
and $D_n$ replaces $D$. 
(It is $D_n$ replacing $D$, not $\tilde D_n$.
The conditions of Proposition~\ref{p.free} hold for either of
these, but the conformal map we consider is defined on $D_n$.)
Proposition~\ref{p.free} therefore gives
\begin{equation}\label{e.pevDm}
\Pb{\ev A\md D_m} = h\Bl(\frac{V_m-W_m}{U_m-W_m}\Br) + O(\delta^3)\,.
\end{equation}
Write $f(U,V,W) := h\bl((V-W)/(U-W)\br)$.
We Taylor-expand the right hand side in~\eref{e.pevDm} to second order in
$W_m$ and to first order in $V_m-V_0$ and $U_m-U_0$.  
Together with~\eref{e.pea}--\eref{e.apx} this gives
\begin{multline*}
0=\Eb{\Ps{\ev A\md D_m}}-\Ps{\ev A} =
\\
\frac 12 \,
\p_W^2\! f
\,
\Es{W_m^2}
+
\p_W\! f
\,
\Es{W_m}
+
\p_V\! f\,
\frac {2\,\Es{t_m}}{V_0}
+
\p_U\! f
\,
\frac {2\,\Es{t_m}}{U_0}
+
O(\delta^3)\,.
\end{multline*}
Here, the derivatives of $f$ are evaluated at $(V_0,W_0,U_0)$.
(Note that $V$ is complex valued, and we interpret
$\p_V f$ as an $\R$-linear map from $\C$ to $\R$.)
If we plug  in
$V_0=i+O(\delta^3)$ and $U_0=1+O(\delta^3)$
 (as we have seen in Section~\ref{s.core}, one can certainly
find $v_0$ and $u_0$ satisfying $\phi(v_0)=i+O(\delta^3)$
and $\phi(u_0)=1+O(\delta^3)$ if $r_0$ is large), then 
after some tedious but straightforward computations the above
equality simplifies to
$$
\Es{W_m^2}
+
2\,\Es{W_m}
-8\, \Es{t_m}=O(\delta^3)\,,
$$
while $V_0=2\,i+O(\delta^3)$ and $U_0=1+O(\delta^3)$ give
$$
3\,\Es{W_m^2}
+
8\,\Es{W_m}
-24\,\Es{t_m}=O(\delta^3)\,.
$$
Combining these two relations together implies~\eref{e.martp}
and~\eref{e.timep} in the case $k=0$.
For $k>0$, the proof is basically the same;
the only essential difference is that one must use $\eta_k$ in place of $\eta$.
\QED

\begin{theorem}[Driving process convergence]\label{t.peanocara}
For every positive $\eps_1$,$\eps_2$, $\eps_3$ and $\TT$, there is some positive
$r_1=r_1(\eps_1,\eps_2,\eps_3,\TT)$ such that the following holds.
Let $D=D(\alpha,\beta,a,b)\in\Doms^*$
satisfy $\inrad(D)>r_1$ and $\harm_D(0,\alpha)\in [\eps_1,1-\eps_1]$.
Let $\gamma$ be corresponding UST Peano path, let
$\phi:D\to\H$ denote the conformal map which takes
$a$ to $0$, $b$ to %
$\infty$ and satisfies $|\phi(0)|=1$,
let $\hat\gamma:=\phi\circ\gamma$, parameterized according
to capacity from $\infty$, and let $\W(t)$ denote the
Loewner driving process for $\hat\gamma$.
Then there is a coupling of standard Brownian motion $B:[0,\TT]\to \R$
and $\W$ such that
$$
\PB{\sup\bl\{|\W(t)-B(8\,t)|:t\in[0,\TT]\br\}>\eps_2}<\eps_3\,.
$$
\end{theorem}

\proof
The proof is almost identical to the proof of 
Theorem~\ref {t.sleconv}, where we used Skorokhod's embedding, but
one has to be a little careful because it may happen that  $0$
is ``swallowed'' before time $\TT$.

Let us first assume that $\phi(0)$ is close to $i$, say 
$|\phi (0) -i|<1/100$,
and that $\TT$ is small enough so that
\begin{equation}\label{e.smallT}
\TT\le 1/100 \hbox { and } \Pb{B[0,\TT]\subset[-1/10,1/10]}>1-\eps_3/3\,,
\end{equation}
where $B$ is standard Brownian motion.
Take $\delta=\delta(\eps_1,\eps_2,\eps_3)>0$ small
and $\eps=1/10$. Define $r_0 (\eps, \delta)$ as in 
Proposition~\ref {p.peanoLoc}.
Let $k\in\N$ be the first integer where $\inrad (D_k) \le r_0$
or $\harm _{D_k}(0,\alpha_k) \notin [\eps, 1- \eps]$
and define $\TT_0:=\min\{\TT,t_k\}$, where $t_n$ is as in the proposition.
Exactly as in the proof 
of Theorem~\ref{t.sleconv}, Proposition~\ref{p.peanoLoc} 
implies that we may couple $\W$ 
with a Brownian motion $B$ in such a way that
$$
\PB{\sup\bl\{|\W(t)-B(8\,t)|:t\in[0,\TT_0]\br\}>\eps_2/3}<\eps_3/3\,,
$$
if $\inrad (D) \ge r_1$ and $r_1$ is large enough.
By our assumptions regarding $\TT$, we have with high probability
that  for all $t\in [0,\TT_0]$, $\W(t)\in[-1/5,1/5]$.
If we choose $r_1$ large enough, this guarantees that
$\Pb{\TT_0\ne \TT}<\eps_3/3$ and proves the Theorem when 
(\ref {e.smallT}) is satisfied and $\phi(0)$ is close to $i$.

Consider now a general $\TT<\infty$.
Let $\TT_1>0$ be some constant satisfying (\ref{e.smallT}),
and let $z_0:=\phi^{-1}(2\,i\,\TT^2/\TT_1^2)$.
{}From the Koebe distortion theorem it follows that
there is a constant $c=c(\TT,\eps_1)$ such that
$\inr{z_0}( D)\ge c\,\inrad(D)$.
(See, e.g., Thm.~1.3 and Cor.~1.4 in  \cite{\PommeBDRY}.)  
Consequently, by choosing $r_1$ appropriately larger,
we may invoke the above argument with the basepoint
moved from $0$ to a vertex near $z_0$ and with a smaller $\eps_2$.
Rescaling now completes the proof of the theorem.
\QED

\subsection{Uniform continuity}\label{s.peanounif}

In order to prove convergence with respect to a stronger topology,
 tightness will be needed, and we therefore derive in the present
subsection some 
regularity estimates for UST Peano curves with
respect to the capacity parameterization.
 Some 
results from \cite {\SchSLE} will be used.
 
Let $D\subset \C$ be a simply connected domain
containing $0$, whose boundary
is a $C^1$ simple closed path.  Let $a$ and $b$ be two
distinct points on $\p D$. 
In this section, we consider for large $R$ the UST Peano curve
from a point near $R\,a$ to a point near $R\,b$ on a grid approximation
of $R\, D$.
One reason not to
consider arbitrary domains 
is that we need to partially adapt
to the framework of~\cite{\SchSLE}
in order to quote results from there.
 Also, it is natural (since the UST
Peano curve is asymptotically space filling) to impose 
regularity conditions on $\partial D$ in order to get uniform %
regularity estimates for the UST Peano curve. 

Let $\alpha_D$ and $\beta_D$ be respectively the clockwise and 
anti-clockwise arcs of $\p D$ from
$b$ to $a$.
Given $R$ large, let
$D^R=D(\alpha^R,\beta^R,a^R,b^R)\in\Doms^*$
be an approximation of $(R\,D,R\,\alpha_D,R\,\beta_D)$,
in the following sense.
Fix some sufficiently large constant $C>0$; for example,
$C=10$ would do.
We require $\alpha^R$ to be a simple path in $\Z^2$ satisfying
$\pdist(\alpha^R,R\,\alpha_D)\le C$ and  require
$\beta^R$ to be a simple path in the dual grid $(1/2+\Z)^2$
satisfying $\pdist(\beta^R,R\,\beta_D)\le C$.
We also require $\beta^R\cap\alpha^R=\emptyset$,
of course, and that each of  $a^R,b^R$ is a Peano vertex
adjacent to an endpoint of $\alpha^R$ and an endpoint of $\beta^R$.

Let $\gamma=\gamma^R$ be the UST Peano path in $D^R$.
Let $\phi:D\to\H$ be the conformal homeomorphism satisfying
$\phi(a)=0$, $\phi(b)=\infty$ and $|\phi(0)|=1$.
Let $\phi_R:D^R\to\H$ be the conformal homeomorphism
satisfying $|\phi_R(0)|=1$, taking $a^R$
 to $0$ and $b^R$ to $\infty$.
Then $\lim_{R\to\infty}R^{-1}\phi_R^{-1}(z)=\phi^{-1}(z)$,
uniformly in $\H$. 
(This follows, e.g., from Cor.~2.4 in~\cite{\PommeBDRY}.)  
Let $\hat\gamma:=\phi_R\circ \gamma$, parameterized according
to capacity from $\infty$.
Let $g_t:\H\setminus\hat \gamma[0,t]\to\H$ be the
conformal map with the usual normalization
$ g_t (z)-z\to 0$ when $|z| \to \infty$.

\begin{proposition}[Uniform continuity estimate]\label{p.contp}
For every $\eps>0$ and $\TT>0$
there are some positive
$R_0=R_0(D,\TT,\eps)$ and $\delta=\delta(D,\TT,\eps)$
such that for all $R > R_0$
$$
\PB{\sup\bl\{|\hat\gamma(t_2)-\hat\gamma(t_1)|:
t_1,t_2\in[0,\TT],\, |t_2-t_1|\le\delta\br\}>\eps}<\eps\,.
$$
\end{proposition}

We first prove a slightly modified version of this proposition.

\begin{lemma}\label{l.contp1}
For $0<t_1<t_2<\infty$ let
$
Y(t_1,t_2):=\diam\bl( g_{t_1}\circ\hat \gamma[t_1,t_2]\br)
$.
For every $\eps>0$ there is
a $\delta=\delta(D,\eps)>0$ and an $R_0=R_0(D,\eps)>0$
such that for all $R\ge R_0$
\begin{equation}\label{e.harest}
\PB{\sup\bl\{ |\hat \gamma(t_2)-\hat \gamma(t_1)|:
0\le t_1\le t_2\le \tau ,\,
Y(t_1,t_2)\le \delta\br\} \ge \eps}
<\eps\,,
\end{equation}
where
$\tau:=\inf\bl\{t\ge 0:|\hat\gamma(t)|=\eps^{-1}\br\}$.
\end{lemma}

The proof will use Theorems~10.7 and~11.1.(ii) of~\cite{\SchSLE}.
As explained there, the proofs of these theorems are now
easier, because we have established the conformal invariance of
the UST; Corollary~\ref{c.ust}.

\proofof{Lemma \ref {l.contp1}}
Let $\delta_R$ be a positive function of $R$ such that
$\lim_{R\to\infty}\delta_R=0$.
It suffices to show that~\eref{e.harest} holds for
all sufficiently large $R$ with
$\delta_R$ in place of $\delta$.
Let $Z$ denote the semi-circle  $2\,\eps^{-1}\,\p\U\cap\H$, say.
For $R$ large, let $t_1$ and
$t_2$ be such that 
$|\hat \gamma(t_2)-\hat \gamma(t_1)|$ is maximal
subject to the constraints
$0\le t_1\le t_2\le \tau$ and $Y(t_1,t_2)\le \delta_R$.
Note that  $\min_{t\le\tau}\dist\bl(g_{t}(Z),g_t\circ \hat\gamma[t,\tau]\br)$
is bounded from below, as
$\hat\gamma[0,\tau]\subset \eps^{-1}\overline\U$.
(Considering the harmonic measure from a point near $\infty$,
one deduces that the diameter of $g_{t}\bl((1/2)Z\br)$ is bounded below.
The extremal distance between $g_t(Z)$ and $g_t\bl((1/2)Z\br)$ is 
the same as the extremal distance between $Z$ and $(1/2)Z$.
This shows that $\dist\bl(g_{t}(Z),g_t\circ\hat\gamma[t,\tau]\br)
\ge \dist\bl(g_{t}(Z),g_t((1/2)Z)\br)$ is bounded from below.)
Since $Y(t_1,t_2)\le\delta_R\to0$ as $R\to\infty$
the extremal length of the collection
of simple arcs in $\H\setminus g_{t_1}\circ\hat \gamma[t_1,t_2]$
which separate $g_{t_1}\circ\hat \gamma[t_1,t_2]$
from $g_{t_1}(Z)$ goes to zero.  By conformal invariance
of extremal length, it follows that the extremal length
of the collection of simple arcs in
$\H\setminus\hat \gamma[t_1,t_2]$ which separate
$\hat \gamma[t_1,t_2]$ from $Z$ in
$\H\setminus\hat \gamma[0,t_1]$ tends to zero as well.
In particular, the shortest such arc for the Euclidean 
metric, say $\eta$, satisfies
$\lim_{R\to\infty}\length(\eta)=0$.

We are going to study separately the three cases where $\eta$
is close to the origin, close to the real line but not to the 
origin, and not close to the real line. In each case, we will
see that the existence of such an $\eta$ is very unlikely.   
Let $\ev A$ be the event $|\hat\gamma(t_1)-\hat\gamma(t_2)|\ge\eps$.
For $s>0$ let $\ev X_0(s)$ be the event $\dist(0,\eta)<s$,
and let $\ev X_1(s)$ be the event $\dist(\R,\eta)<s$.
We will prove
\begin{align}
\forall s_1>0\,\,\exists R_0>0\,\,\forall R>R_0\quad
&\Pb{\ev A\setminus\ev X_1(s_1)}<\eps\,,
\label{e.X2}
\\
\forall s_0>0\,\,
\exists s_1>0\,\,\exists R_0>0\,\,\forall R>R_0\quad
&\Pb{\ev A\cap\ev X_1(s_1)\setminus\ev X_0(s_0)}<\eps\,,
\label{e.X1}
\\
\exists s_0>0\,\,\exists R_0>0\,\,\forall R>R_0\quad
&\Pb{\ev A\cap\ev X_0(s_0)}<\eps\,.
\label{e.X0}
\end{align}
Using these statements, the proof of the lemma is completed
by choosing $s_0$ according to~\eref{e.X0}, then choosing
$s_1$ according to~\eref{e.X1}, and finally choosing
$R_0$ according to~\eref{e.X2},~\eref{e.X1} and~\eref{e.X0}.
  
\begin{figure}
\SetLabels
(.3*.7)$A$\\
\R(.42*.3)$\hat\gamma_1$\\
\B(.3*.5)$\hat\gamma_2$\\
\L(.42*.77)$\hat\gamma_3$\\
\L(.8*.8)$\hat\gamma(t^*)$\\ %
\endSetLabels
\centerline{\AffixLabels%
{\includegraphics*[height=3in,angle=90]{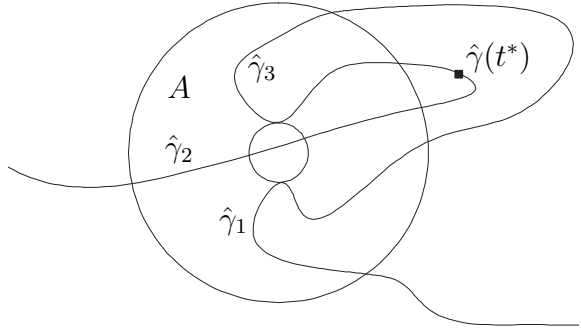}}}
\caption{\label{f.topconfg}The paths $\hat \gamma_1,\hat\gamma_2$ and
$\hat\gamma_3$ and the annulus $A$.}
\end{figure}

We start with~\eref{e.X2}.  Fix some $s_1>0$,
and assume that $\ev A\setminus \ev X_1(s_1)$ holds.
We also assume that $\eps< s_1$.
There is no loss of generality in that assumption, since
$\ev A$ is monotone decreasing in $\eps$.
Since $\lim_{R\to\infty}\length(\eta)= 0$, for large $R$
the two endpoints of $\eta$ must be in $\hat \gamma[0,t_1]$. 
Because $\hat \gamma$ tends to $\infty$
with $t$, it is clear that $\hat \gamma[t_2,\infty)\cap \eta\ne\emptyset$.
In fact, the crossing number of $\hat \gamma[t_2,\infty)$ and $\eta$
must be $\pm1$, since $\hat \gamma$ and $\eta$ are simple curves.
Consider the concentric annulus $A$ whose inner circle 
is the smallest circle surrounding $\eta$ and whose outer
circle has radius $\eps/4$.  
Let $\ball$ denote 
the open disk bounded by the outer circle of $A$,
and note that $\ball\subset\H$, by our assumption $\eps<s_1$.
On the event $\ev A$, %
there is a $t^*\in [t_1,t_2]$ such %
that the distance from $\hat\gamma(t^*)$ to $\eta$ %
is at least $\eps/2$.
In particular $\hat\gamma(t^*)\notin \ball$.  %
Now, $\eta$ separates $\hat\gamma(t^*)$ from $\infty$
in $\H\setminus\hat\gamma[0,t_1]$.
Therefore, if $\ev A\setminus X_1(s)$ holds, %
then $\hat\gamma[0,t_1]\cup\eta$ separates $\hat\gamma(t^*)$ from $\infty$.
Since $\hat\gamma$ is a simple path,
this implies that the arc of $\hat\gamma[0,t_1]$
between the two points $\overline\eta\cap\hat\gamma[0,t_1]$
does not stay in $\ball$.
Hence, $\hat \gamma[0,\infty)\cap \ball$ has three
distinct connected components, say $\hat \gamma_1,\hat \gamma_2,\hat \gamma_3$
each of which intersects the inner circle of $A$,
such that $\hat \gamma_1,\hat \gamma_3\subset\hat \gamma[0,t^*]$ %
and $\hat \gamma_2\subset\hat \gamma[t^*,\infty)$ and %
$\hat \gamma_2$ separates $\hat \gamma_1$ from $\hat \gamma_3$ within $\ball$.
See Figure~\ref{f.topconfg}.

Note that adjacent to one side of $\phi^{-1}_R(\hat\gamma_2)$ lies $T$, the UST,
and $T^\dagger$, the
dual UST, is adjacent to the other side. 
Both are connected, and
they do not intersect $\phi^{-1}_R(\hat \gamma)$.  It follows
that there are paths $\chi_1\subset T$ and $\chi_2\subset T^\dagger$
with endpoints in $\phi_R^{-1}(\p \ball)$
each of which intersects the inner boundary of $\phi_R^{-1}(A)$.
But the diameter of the inner boundary of $R^{-1}\phi_R^{-1}(A)$ goes
to zero as $R\to\infty$ and the distance between the two boundary
components of $R^{-1}\phi_R^{-1}(A)$ does not.
Hence, by~\cite[Theorem 10.7]{\SchSLE},
the probability that such a configuration
appears somewhere goes to zero with $R$.
(Although the result from~\cite{\SchSLE} refers to the
UST in the whole plane, the proof is local, and since 
we are bounded away from the boundary, the result is applicable here.)
This proves~\eref{e.X2}.

Now fix $s_0>0$ and let $s_1>0$ be much smaller.
Assume that $\ev A\cap\ev X_1(s_1)\setminus\ev X_0(s_0)$ holds,
$\eps<s_1$, 
and that $R$ is large.  
Also assume that $\eta$ is closer to $[0,\infty)$ than
to $(-\infty,0]$.  
Note that $\eta$ is then bounded away from $(-\infty,0]$.
Let $A$ be defined as above, and let
$\ball$ be the intersection of $\H$ with the disk bounded
by the outer boundary component of $A$.
We now need to consider two distinct possibilities.
Either both endpoints of $\eta$ are on $\hat\gamma[0,t_1]$,
and then the configuration is topologically as in the argument
for~\eref{e.X2}, or one endpoint of $\eta$ is on $[0,\infty)$.
But it is easy to see that in either case  there is a simple
path in $T^\dagger$ which intersects $\phi_R^{-1}(\eta)$ 
whose endpoints are in $\phi_R^{-1}(\p \ball)$, by an argument very
similar to the one given above.
Now~\cite[Theorem~11.1.(ii)]{\SchSLE} shows that these
events have small probabilities if $s_1$ is small.
The case where $\eta$ is closer to $(-\infty,0]$ is
treated similarly, with the roles of the tree and the dual tree
switched.  Thus~\eref{e.X1} is established.

To prove~\eref{e.X0}, let $s>0$, and
let $v\in\Z^2$ be a vertex closest to $\phi_R^{-1}(i\,s)$.
Let $v^\dagger\in (\Z+1/2)^2$ be a dual vertex adjacent to $v$.
Let $\chi$ be the simple path from $v$ to $\alpha^R$ in $T$,
and let $\chi^\dagger$ be the simple path from
$v^\dagger$ to $\beta^R$ in $T^\dagger$.
We may sample $\chi$ by running simple random walk
from $v$ on $H(D^R)$ stopped on hitting $\alpha^R$,
and loop-erasing it.
It therefore follows by Proposition~\ref{p.free}
that if $s$ is sufficiently small,
then the diameter of $\phi_R(\chi)$ is smaller than
$\eps/10$ with probability at least $1-\eps/10$.
Moreover, there is some $s_0>0$ such that for
all sufficiently large $R$ with probability
at least $1-\eps/10$ the distance from $\phi_R(\chi)$
to $0$ is at least $s_0$,
and the same two estimates will hold for $\chi^\dagger$.
Let $D'$ be the domain bounded by
$[v,v^\dagger]\cup\chi\cup\chi^\dagger\cup\p D^R$
which has the initial point $a^R$ of $\gamma$ on its boundary. 
Note that $\gamma$
crosses the boundary of $D'$ exactly once, through the
segment $[v,v^\dagger]$.  In particular,
if $\diam\bl(\phi_R(D')\br)<\eps$ and $\ev A$ holds, %
then $\eta$ is not contained in $\phi_R(D')$. %
This proves~\eref{e.X0}, and completes the proof of the lemma.
\QED

\proofof{Proposition~\ref{p.contp}}
Theorem~\ref{t.peanocara} implies that
we may find some $r>0$ such that
$\PB{\sup\bl\{|W(t)|:t\in[0,\TT]\br\}\ge r}<\eps/4$
for all sufficiently large $R$. %
Let $\eps':=\min\{\eps,r^{-1}\}$,
and let $\delta'$ denote the $\delta$ obtained by using
Lemma~\ref{l.contp1} with $\eps'$ in place of $\eps$.
Since Brownian motion is a.s.\ continuous, Theorem~\ref{t.peanocara}
implies that there is some $\delta>0$ such that if $R$ is large enough
we have 
$$
\PB{\sup\bl\{|W(t_1)-W(t_2)|:t_1,t_2\in[0,\TT],\,|t_1-t_2|<\delta
\br\}\ge\delta'}<\eps/4\,.
$$ 
Lemma~\ref{l.diambound} 
applied to the path
$t\mapsto g_{t_1}\circ\gamma(t-t_1)-W(t_1)$ now implies
$$
\PB{\sup\bl\{Y(t_1,t_2):0\le t_1\le t_2\le \TT,\,|t_1-t_2|<\delta
\br\}\ge C(\delta^{1/2}+\delta')}<\eps/4\,. %
$$ 
Now the proof is completed by using Lemma~\ref{l.contp1}.
\QED

\subsection{Consequences}
\label {s.peanoconclusions}

In this section we gather some consequences, starting with the
following two theorems.

\begin{theorem}[Chordal \SLEkk8/ traces a path]\label{t.slecont}
Let $\tilde g_t$ denote the chordal \SLEkk8/ process
driven by $B(8\,t)$, where $B(t)$ is standard Brownian motion.
Then a.s.\ for every $t>0$ the map $\tilde g_t^{-1}$ extends
continuously to $\closure\H$ and
$\tilde\gamma(t):=\tilde g_t^{-1}(B(8\,t))$ is a.s.\ continuous.
Moreover, a.s.\ $\tilde g_t^{-1}(\H)$ is the
unbounded component of $\H\setminus \tilde\gamma[0,t]$ for
every $t\ge 0$.
\end{theorem}

\begin{theorem}[Peano path convergence]\label{t.peano}
Let $D\subset \C$ be a  domain containing zero, such that 
$\p D$ 
is a $C^1$-smooth simple closed path.
 Let $\p D=\alpha_D\cup\beta_D$ be a partition of the boundary of $D$
 into two
nontrivial complementary arcs.
For $R>0$ let $(D^R,\alpha^R,\beta^R)$ %
be an approximation of $(R\,D,R\,\alpha,R\,\beta)$, as described in Section~\ref{s.peanounif}.
Let $\gamma=\gamma^R$ denote the UST Peano curve in $D^R$ with %
the corresponding boundary conditions.
Let $\phi_R:D^R\to\H$ denote the conformal map %
which takes the initial point of $\gamma$ to $0$, the
terminal point to $\infty$ and satisfies $|\phi_R(0)|=1$.
Let $\hat\gamma:=\phi_R\circ \gamma$, parameterized by capacity
from $\infty$.
Then the law of $\hat\gamma$  %
 tends weakly to the law of $\tilde\gamma$ from Theorem~\ref{t.slecont}. %
\end{theorem}

Here, we think of $\hat\gamma$ and $\tilde\gamma$ as %
elements of the space of continuous maps from $[0,\infty)$ to $\closure\H$, %
with the topology of locally uniform convergence. %

A consequence of the theorem is that $R^{-1}\gamma$ is %
close to $R^{-1}\phi_R^{-1}\circ\tilde\gamma$.
That is, we may approximate the UST Peano path $\gamma$
by the image of chordal \SLEkk8/ in $D^R$. %

The analogue of Theorem~\ref{t.slecont} was proven in~\cite{\RSsle}
for all $\kappa\ne 8$, but the particular case $\kappa=8$ could not be handled
there.  It is fortunate that the convergence of the UST Peano path
to \SLEkk8/ settles this problem.
By Remark 7.5 from~\cite{\RSsle} it follows that with the notations %
of Theorem~\ref{t.slecont} for %
$\kappa\ge 8$ we have $\tilde g_t^{-1}(\H)=\H\setminus \tilde\gamma[0,t]$ for %
every $t\ge 0$ a.s. %

\proofof{Theorems~\ref{t.slecont} and~\ref{t.peano}}
Let $W(t)=W_R(t)$ denote the chordal Loewner driving process for $\hat\gamma$.
Fix a sequence $R_n \to \infty$.
First, note that 
the family of laws of $\hat\gamma$ is tight, %
because of Proposition~\ref{p.contp} and the Arzela-Ascoli Theorem 
(see, for instance, \cite[Theorem 2.4.10]{\KaratsasShreve}). %
Also, Theorem~\ref{t.peanocara} implies
that the law of $W$ converges weakly to the law of $B(8\,t)$. %
Hence, there is a subsequence of $R_n$ such that the law
of the pair $(\hat\gamma,W)$ converges weakly to some probability
measure $\mu$.
Let $(\gamma^*,W^*)$ be random with law $\mu$. 
Then we may identify $W^*(t)$ with $B(8\,t)$.
By the chordal analogue of Lemma~\ref{l.lconv}, which
is valid with the same proof, it follows that for all
$t>0$,
$\tilde g_t^{-1}(\H)$ is the unbounded component
of $\H\setminus\gamma^*[0,t]$.  Since $\gamma^*$ is
continuous, elementary properties of
conformal maps imply that $\tilde g_t^{-1} $ extends
continuously to $\closure\H$  (e.g., Theorem 2.1 in~\cite{\PommeBDRY}).  
It is easy to verify that a.s.\ for every $t>0$,
$\tilde g_t^{-1}(B(8t))=\gamma^*(t)$, using the
fact that $\gamma^*[t,t']$ is contained in a small
neighborhood of $\gamma^*(t)$ when $t'-t>0$ is small. %
This proves Theorem~\ref{t.slecont}.
Because the law of the limit path $\gamma^*$ does not depend on the subsequence,
the original sequence converges, and so Theorem~\ref{t.peano} is proved as well.
\QED

We now list some easy 
consequences of Theorems~\ref{t.slecont} and \ref{t.peano}.

\begin {corollary}[Radial \SLEkk8/ traces a path]
\label {c.radialslecont}
Let $\tilde g_t$ denote a radial \SLEkk8/ 
process driven by $\WWzeta(t)=\exp \bl( i B(8t)\br)$ where
$B$ is standard Brownian motion. Then, almost surely, for 
every $t>0$, the map $\tilde g_t^{-1}$ extends continuously 
to $\overline \U$. Moreover $\tilde g_t^{-1} \bl(\WWzeta (t)\br)$
is almost surely continuous.
\end {corollary}
  
\proof
This follows readily from Theorem~\ref{t.slecont} and the absolute
continuity relation between radial and chordal \SLEkk8/ derived in 
\cite {\LSWii}, Proposition~4.2.
\QED

\proofof{Theorem \ref{t.intropeano}}
Define $G_\delta= \delta D^{1/\delta}$
where $D^{1/\delta}$ is defined as in Theorem~\ref{t.peano}.
Consider the situation of Theorem~\ref{t.peano}.
As previously remarked, it follows from~\cite[Cor.~2.4]{\PommeBDRY}  
that $\lim_{R\to\infty}R^{-1}\phi_R^{-1}(z)=\phi^{-1}(z)$,
uniformly in $\H$. 
Consequently, Theorem~\ref {t.peano} shows that for all $\TT>0$, the UST Peano curve 
 scaling limit up to capacity $\TT$ from $b$ is equal  to 
$\phi^{-1}\circ \tilde \gamma$ up to time $\TT$.
It therefore suffices to prove that for all $\eps>0$ there is an $\eps'>0$
such that for all sufficiently large $R$ with probability at least
$1-\eps$ the part of $\gamma$ after the first time it hits the $\eps'$-neighborhood
of $b$ stays within the $\eps$-neighborhood of $b$.
This is easily proved by the same argument used to prove~\eref{e.X0} applied
to the reversal of the UST Peano path, which is also a UST Peano path.
\QED

\begin {corollary}[Path reversal] \label {c.sle8rev} 
The law of the chordal \SLEkk8/ curve is invariant under
simultaneously reversing time and inverting in the unit circle,
up to a monotone increasing time-change.
More precisely, if $\tilde \gamma$ is the chordal \SLEkk8/
curve from $0$ to infinity defined in Theorem~\ref {t.slecont},
then a time-change of $\bl( -1 / \tilde \gamma (1/t), t \ge 0 \br)$ has the 
same law as $\tilde \gamma$.
\end {corollary}  
 
\proof
This follows immediately from the fact that the reversed UST Peano curve
is also a UST Peano curve.
\QED

\section{Random walk estimates}\label{s.RW}

The goal of this section is to prove the remaining random walk estimates
and thereby complete the proofs of the theorems.
Basically, we show that under certain boundary conditions, discrete harmonic 
functions converge to continuous harmonic functions 
satisfying corresponding boundary conditions, as the mesh of the 
grid goes to zero. 
This general principle is not new, of course (see, e.g.,~\cite{\Collatz}), 
but it seems that the precise statements which are needed here do not 
appear in the literature. 
In particular, our results make no smoothness assumptions on the
boundary.
It should perhaps be noted that some of the following proofs (and most 
likely the results too) are special to two dimensions.

\subsection {Preliminary lemmas}\label{s.RWprel}

We now state some lemmas on discrete harmonic functions, which will 
be helpful in the proofs of 
Proposition~\ref {p.phit}, Lemma~\ref {l.dG} and Proposition~\ref {p.free}.

\medskip

For  $\delta>0$, define the discrete derivatives
$$ \partial_x^\delta f(v) := 
\delta^{-1} ( f (v+ \delta) - f (v) )
,\qquad
\partial_y^\delta f(v) 
:= \delta^{-1} ( f( v + i \delta) - f(v) ) .
$$
Let $\Doms_\delta:=\{\delta D:D\in\Doms\}$;
that is, domains adapted to the grid $\delta\Z^2$.
Similarly, for $D_\delta=\delta D\in\Doms_\delta$, define
$\Vd {D_\delta}:=D_\delta\cap \delta\Z^2=\delta\V D$ and
$\Vbd {D_\delta}:=\delta\Vb D$.

\begin{lemma}[Discrete derivative estimate]\label{l.dde}
There is a constant $C>0$ such that for every $D\in\Doms$
and every bounded function $h:\V D\cup \Vb D\to\R$ that
is harmonic in $\V D$,
\begin{equation}
\label{e.dde}
\p_x^1h(0)\le C\,\inrad(D)^{-1}\,\|h\|_\infty\,,\qquad
\p_y^1h(0)\le C\,\inrad(D)^{-1}\,\|h\|_\infty\,.
\end{equation}
\end{lemma}

This lemma is proved using Green's functions in \cite [Thm.~1.7.1]{\Lbook};
see also~\cite[Lem.~7.1]{MR99j:52021} for a proof of the analogous
statement in the triangular lattice using the maximum principle.
In Section~\ref{s.other}, we rewrite and adapt the proof from
\cite{\Lbook} to more general walks on planar lattices.
One can also rather easily prove the lemma using coupling.

\begin {lemma}
\label {l.harnack}
For all $\eps>0$ and $k\in\N$ there exists  a
$c= c_k(\eps)>0$ such that the following
always holds.  Let $\delta\in(0,c^{-1})$ and let
$D\in\Doms_\delta$ satisfy $\inrad(D)\ge 1/2$.
Let $\p_{a_1}^\delta,\dots,\p_{a_k}^\delta\in\{\p_x^\delta,\p_y^\delta\}$.
Let $h:\Vd D\cup\Vbd D\to[0,\infty)$ be non-negative and harmonic
in $\Vd D$. 
If $v\in\Vd D$ satisfies $|\psi_D(v)|\le 1-\eps$, then
\begin {equation}
\label {difh}
\bl|\p_{a_1}^\delta\p_{a_2}^\delta\cdots\p_{a_k}^\delta h(v) \br|
\le  c \,h(0)
\,.
\end {equation}
\end {lemma}

Note that the case $k=0$, which is included, is a kind of Harnack inequality.

It is easy to give quantitative
estimates for $c_k(\eps)$, but they will not be needed here.
Only $k \le 3$ will be used in the sequel.

\medskip

In the proof of the lemma, the following
simple conformal geometry consequences of
the Koebe distortion theorem~\cite[Thm.~1.3]{\PommeBDRY} will be needed. 
Let $D,\eps$ and $v$ be as in the statement of the lemma.
First, note that $1/4\le \inrad(D)\,\psi_D'(0)\le 1$
follows from the Koebe 1/4 Theorem and the Schwarz Lemma, 
respectively.
Let $\ell=\ell(\eps)$ be large, set $z_j:= j\,\psi_D(v)/\ell$,
and $w_j:=\psi_D^{-1}(z_j)$, $j=0,1,\dots,\ell$.
The Koebe distortion theorem gives upper and lower bounds 
for $\inrad(D)\,|\psi_D'|$ on the preimage of the line segment $[0,z_\ell]$. 
This implies that there is a constant $c_1=c_1(\eps)>0$
such that $\inr {w_j}(D)\ge c_1\inrad(D)$,
 and that if $\ell=\ell(\eps)$ is large then,
$|w_j-w_{j-1}|\le c_1\,\inrad(D)/20$, $j=1,\dots,\ell$.
In particular, if $v_j$ is the vertex in $\Vd D$ closest to $w_j$,
then, provided that $\delta$ is sufficiently small,
$|v_j-v_{j-1}|\le \inr{v_{j-1}}(D)/10$.

\proof
We start with $k=0$. Suppose first that $|v|\le \inrad(D)/10$.  Let 
$W\subset\Vd D$ be the set of vertices $w$ satisfying $h(w)\ge h(v)$.
Then $W$ contains a path from $v$ to $\p D$.
But the probability $p$ that the path traced by simple random walk from
$0$ before exiting $D$ separates
$v$ from $\p D$ is bounded away from $0$.
On that event, the simple random walk hits $W$ before exiting $D$.
Consequently $h(0)\ge p\, h(v)$, as needed.
For arbitrary $v\in \Vd D$ satisfying $|\psi_D(v)|\le 1-\eps$,
as we have noted, the Koebe distortion theorem implies 
that there is an $\ell=\ell(\eps)$ depending only on $\eps$, and
a sequence $0=v_0,v_1,\dots,v_\ell=v$ in $\Vd D$ with $\ell\le\ell(\eps)$
such that $|v_{j}-v_{j-1}|\le\inr {v_j} (D)/10$
for each $j=1,\dots,\ell$.
Consequently, iterating the above result gives
$h(0)\ge p^{\ell}\,h(v)$, and proves the case $k=0$.

Using the above, we know that $h(w)\le c' h(0)$ on
the set of vertices $w\in \Vd D$ such
that $|w-v|\le\inr v (D)/10$, where $c'=c'(\eps)$
is some constant depending only on $\eps$.
Consequently, the case $k=1$ now follows from Lemma~\ref{l.dde}
applied with $v$ translated to $0$.

For $k>1$, the proof is by induction.
By the above, we may assume $v=0$.
Let $M$ be the maximum of $\bl|\p_{a_k}^\delta h(w)\br|/h(0)$ on the
set $V$ of vertices $w\in\Vd D$ satisfying $|w|\le \inrad(D)/10$.
The above shows that $M$ is bounded by a universal constant.
Since $\p_{a_k}^\delta h$ is discrete-harmonic on $V$, the
proof is completed by applying the inductive hypotheses
to the function $\p_{a_k}^\delta h(w)+M\,h(0)$.
\QED

\begin{lemma}[Continuous harmonic approximation]\label{l.goharm}
For every $\eps>0$ there is some $r_0=r_0(\eps)>0$ such
that the following holds.  If $D\in\Doms$ satisfies
$\inrad(D)\ge r_0$ and $h:\V D\cup \Vb D \to[0,\infty)$ is
discrete-harmonic in $\V D$, then there 
exists a harmonic
function $h^* : D \to [0, \infty)$ such that
\begin {equation}
\label {e.hath}
\bl| h^* (v)-h(v)\br|\le \eps\,h(0)
\end {equation}
holds for every vertex $v\in \V D$ satisfying
$\bl|\psi_D(v)\br|< 1-\eps$.
\end{lemma}

\proof
Suppose that the Lemma is not true.
Then, there exists $\eps>0$ and a sequence
 of pairs $(D_n,h_n)$, where $D_n\in\Doms$ satisfies
$\inrad(D_n)\ge n$ and $h_n>0$ is
discrete harmonic in $\V {D_n}$, satisfies $h_n(0)=1$, but
(\ref {e.hath}) fails for every harmonic function $\hat h$.

Set $\delta=\delta_n := 1/ \inrad (D_n)$.
Our objective is to apply compactness to show that
the maps $h_n\circ\psi_D^{-1}$ converge locally uniformly
in $\U$ as $n\to\infty$ along some subsequence to 
some harmonic $\hat h$, so that (\ref {e.hath}) does hold for some $n$.
We put $h^n ( v) := h_n ( v/\delta_n )$.

First, standard compactness properties of conformal maps say
that one can take a subsequence 
such that the maps $\delta_n \psi_{D_n}^{-1}$
converge locally uniformly in $\U$ to
some conformal map, say $\phi$. 
(This follows, for example, from the Arzela-Ascoli theorem,
together with~\cite[Cor.\ 1.4]{\PommeBDRY} with $z=0$ and  
part two of~\cite[Theorem 1.3]{\PommeBDRY}.)  
If $K\subset\U$ is compact, then Lemma~\ref{l.harnack}
shows that there is a constant $C>0$ such that for all
sufficiently large $n$ in the subsequence,
the discrete derivatives  $|\p^\delta_x h^n|$
and $|\p^\delta_y h^n|$ are bounded by $C$ in $\phi(K)\cap \V{\delta D_n}$.
By a variant of the Arzela-Ascoli Theorem,
it then follows that there is some continuous
$h^*: \phi (\U) \to [0, \infty)$ and
a further subsequence such that for every compact $K\subset\phi(\U)$,
$$%
\sup\Bl\{\bl| h^n (v)- h^*(v)\br|: 
v\in K\cap\delta\Z^2\Br\}\to 0
$$%
along the subsequence.
The same argument may also be applied to prove the convergence
of the discrete derivatives of $h^n$ to arbitrary order,
possibly in a further subsequence.
Obviously,
the discrete derivatives of $h^n$
will converge to the corresponding continuous derivatives
of $h^*$; that is,
\begin{equation}\label{e.dconv}
\sup\Bl\{\bl|
\p^\delta_{a_1}\cdots\p^\delta_{a_k} h^n (v)- 
\p_{a_1}\cdots\p_{a_k} h^*(v)\br|: v\in K\cap\delta\Z^2\Br\}\to 0\,,
\end{equation}
where $\p_{a_j}^\delta\in\{\p_x^\delta,\p_y^\delta\}$
and $\p_{a_j}\in\{\p_x,\p_y\}$ is the corresponding continuous
derivative, $j=1,2,\dots,k$.
The fact that $h^n$ is discrete-harmonic translates to
$(\p^\delta_x)^2 h^n(v-\delta)+ (\p^\delta_y)^2 h^n(v-i\delta)=0$.
Therefore,~\eref{e.dconv} shows that $h^*$ is harmonic. 
This completes the proof.
\QED

\begin{lemma}[Boundary hitting]\label{l.bdhit}
For every $\eps_1,\eps_2>0$ there is a $\delta=\delta(\eps_1,\eps_2)>0$
such that if $D\in\Doms$ and $w\in\V D$
is a vertex satisfying $\bl|\psi_D(w)\br|\ge 1-\delta$, then
the probability that simple random walk started at $w$
will hit
$$
\bl\{v\in \V D:\bl|\psi_D(v)-\psi_D(w)\br|>\eps_1\br\}
$$
 before hitting $\p D$ is at most $\eps_2$.
\end{lemma}

\begin{figure}
\SetLabels
\L(.33*.36)$z$\\
\T(.31*-.01)$w$\\
\L(.33*.2)$A_1$\\
\R(.02*.2)$A_2$\\
\L(.62*.2)$A_3$\\
(.9*.3)$K_1$\\
(.18*.3)$K_3$\\
(.5*.3)$K_2$\\
\endSetLabels
\centerline{\AffixLabels{\includegraphics*[height=2.3in]{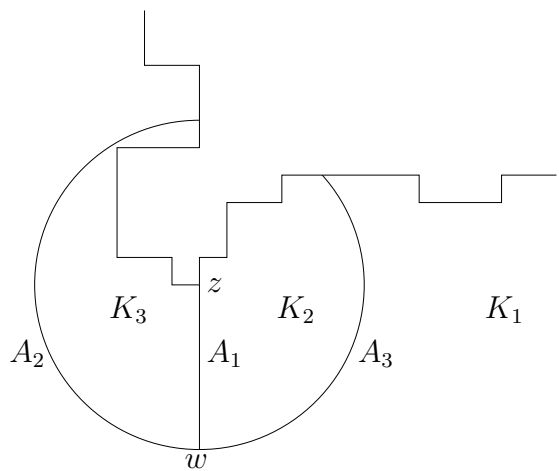}}}
\caption{\label{f.arcs} The arcs $A_j$ and the components $K_j$.}
\end{figure}

\proof
We first prove the lemma in the case where $\eps_2$
is very close to $1$.
Let $\delta>0$ be much smaller than $\eps_1$.
Fix some vertex $w\in\V D$,
and suppose that $\bl|\psi_D(w)\br|\ge 1-\delta$.
Let
$$
\alpha:=\bl\{z\in D:\bl|\psi_D(z)-\psi_D(w)\br|=\eps_1\br\}\,.
$$
Let $z_1$ be a point in $\p D$ closest to $w$
and set $r:=\dist(w,\p D)=|z_1-w|$.
Let $A_1$ be the line segment $[w,z_1]$.
Let $Q$ be the connected component of
$C(z_1,r)\cap D$ which contains $w$, where
$C(z,r)$ denotes the circle of radius $r$ and center
$z$.
Then $Q$ is an arc of a circle.  
Let $A_2$ and $A_3$ denote the two connected components of
$Q\setminus\{w\}$.
See Figure~\ref{f.arcs}. 
For $j=1,2,3$, let $K_j$ be the connected component
of $D\setminus (A_1\cup A_2\cup A_3)$
which does not have $A_j$ as a subset of its boundary.

Because $\delta$ is small compared to $\eps_1$,
the Koebe distortion theorem  
(e.g., Corollaries 1.4 and 1.5 in \cite {\PommeBDRY})  
shows that $\alpha\cap C(w,r/8)=\emptyset$.
For $j=1,2,3$, let $\ev S_j$ be the collection of all
paths which stay in $K_j$ from the first time they hit $C(w,r/8)$
until the first exit from $D$.
Let $B(t)$ denote Brownian motion started from $w$.
It is easy to see that there is a universal constant $c_1>0$
such that $\Ps{B\in\ev S_j}>c_1$ for $j=1,2,3$.
For example, to prove this for $j=3$
observe that the collection of Brownian paths
which first hit $C(w,r/8)$ in $K_3$ and later
hit $A_3$ before $A_1\cup A_2$ has probability
bounded away from zero.

Suppose for the moment that $\alpha$ intersects $A_1$ and $A_2$.
Consider a subarc $\alpha'\subset\alpha$ whose endpoints are
in $A_1$ and $A_2$, which is minimal with respect
to inclusion.  Then $\alpha'\subset\closure{K_3}$
or $\alpha'\subset\closure{K_1\cup K_2}$.
If $\alpha'\subset\closure{K_3}$, then $\alpha'$
separates $C(w,r/8)$ from $\p D$ in $K_3$.
Consequently, on the event $B\in\ev S_j$, $B$ hits
$\alpha$ before hitting $\p D$.  However,
by choosing $\delta$ to be sufficiently small and invoking
conformal invariance of harmonic measure,
we may ensure that the latter event has probability
smaller than $c_1$.  An entirely similar argument
rules out the possibility that $\alpha'\subset \closure{K_1\cup K_2}$.
Similarly, it is not possible that $\alpha$ intersects both
$A_1$ and $A_3$ or that $\alpha$ intersects both
$A_2$ and $A_3$.  Hence, there is some $j\in\{1,2,3\}$ such
that $\alpha\cap K_j=\emptyset$.  Let $j'$ be such a $j$.
 
By the convergence of simple random walk to Brownian motion,
it is clear that there is some universal constant $r_0>0$ such that
if $r>r_0$, then the probability that simple random
walk started from $w$ is in $\ev S_{j'}$ is
at least $c_1/2$.  This establishes the lemma
in the case where $\eps_2\in(1-c_1/2,1]$ and $r>r_0$.
Suppose $r\le r_0$.  Then there are two grid paths
of bounded length starting from $w$ to $\p D$ that
are disjoint except at $w$.   If
$\alpha$ intersects both these paths, then this gives
a lower bound for the continuous harmonic measure
of $\alpha$ from $w$. Consequently, by making
$\delta$ small enough, we can make sure that
this does not happen.  
Thus, again, with probability bounded away
from $0$ the random walk from $w$ hits $\p D$
before $\alpha$, since it may follow any one of
these two paths.
This proves the lemma in the case where $\eps_2\in(c_2,1]$,
where $c_2$ is some universal constant.

The Koebe distortion theorem implies that there is 
a constant $c>0$ such that if $v_1,v_2\in \V D$
are neighbors, then $1-|\psi_D(v_2)|\le c\, \bl(1-|\psi_D(v_1)|\br)$.
(See, e.g., Corollaries 1.4 and 1.5 in \cite {\PommeBDRY}.)  
Consequently, we may iterate the above restricted
case of the lemma and use the Markov property,
thereby proving the lemma for arbitrary $\eps_2>0$.
\QED

\subsection {The hitting probability estimate} \label{s.phitpf}

\proofof{Proposition \ref {p.phit}}
Let $\eps_1>0$ be much smaller than $\eps$.
We consider the 
discrete harmonic function
 $h(w):= H(w,u) / H(0,u)$.  For
$\delta>0$ let 
$$
V(\delta, \eps_1):=\bl\{
z\in\V D:\bl|\psi_D(z)\br|\ge 1-\delta,\,
|\psi_D(z)-\psi_D(u)|>\eps_1\br\}
\,.
$$

Our first goal is to show that 
for every $\eps_1\in(0,1/4)$,
there is some $\delta=\delta(\eps_1)>0$ and some $r_0=r_0(\eps_1)>0$
such that 
\begin{equation}
\label{e.vanish}
\max\bl\{h(z):z\in V(\delta,\eps_1)\}<\eps_1
\end{equation}
provided that $\inrad(D)>r_0$.
This will be achieved by first showing that $h$ is not too 
large on the set 
$$
W:=\bl\{z\in \V D: \eps_1/2\ge |\psi_D(z)-\psi_D(u)|\ge \eps_1/3\br\}
$$
and then letting $\delta$ go to zero and appealing to Lemma \ref {l.bdhit}.

Assume that $\inrad(D)$ is sufficiently large so that
 any nearest neighbor path from $0$ to $u$ in $D$
has a vertex in $W$. 
Let $M$ denote the maximum of $h$ on $W$.
We claim that $M$ is bounded 
by a constant $c=c(\eps_1)$ depending
only on $\eps_1$.
Indeed,
let 
$K$ be the set of all $v\in\V D$ satisfying
$h(v)\ge M/2$ and let $K'$ be the union of all
edges where both endpoints are in $K\cup\{u\}$.  Then
the maximum principle shows that $K'$ is connected and contains 
a simple nearest neighbor path $J$ joining $W$ to $u$ 
whose vertices are in $\{  z \in\V D 
\ : \ | \psi_D (z) - \psi_D (u) | \le \eps_1/ 2 \}\cup\{u\}$.
Note that, in particular, 
$\diam(\psi_D( J))\ge\eps_1/3$.
Consequently, the continuous harmonic measure from $0$ of
$J$ in $D$ is bounded from below by some constant $c_1(\eps_1)>0$. 

 We claim that the discrete harmonic measure $H_D(0,J)$ of $J$ 
at the origin is
also bounded away from $0$ if $r_0$ is large enough.
Indeed, let $D' = D\setminus J$
and let $A$ be the arc on $\p\U$ corresponding to $J$
under the map $\psi_{D'}$.  The length of $A$ is
bounded from below, since it is equal to $2\pi$ times
the harmonic measure of $J$.
Let $A'$ denote the middle subarc of $A$ having
half the length of $A$.  By Lemma~\ref{l.bdhit}
applied to the domain $D'$, it follows that there is a
$c_2=c_2(\eps)\in(0,1/10)$,
such that $H_{D'}(v,J)\ge 1/2$ on vertices $v$ such
that $\psi_{D'}(v)$ is within distance $c_2$ of $A'$.
Using Lemma~\ref{l.goharm} with $\epsilon$ replaced by $c_2/4$,
we find that if $r_0$ is large, there is a
non-negative harmonic function $h^*_J:D'\to[0,\infty)$
such that $\bl|h^*_J(v)-H_{D'}(v,J)\br|\le c_2/4$
for all $v\in \V {D'}$ satisfying $\bl|\psi_{D'}(v)\br| < 1-c_2/4$.
Take $a\in A'$ and $z:=(1-c_2/2)a$.  Then it follows from the Koebe 
distortion theorem (as in the argument towards the end of the
proof of Proposition~\ref{p.bas}) that we may find a vertex 
$v\in\V {D'}$ such that $\bl|\psi_{D'}(v)-z\br|<c_2/4$, assuming
that $r_0$ is large enough.  Thus,
$h^*_J(v)> H_{D'}(v,J)-c_2/4\ge 1/2-c_2/4>1/4$.
By the Harnack principle applied to $h^*_J$, there is a
universal constant $c_3>0$ such that $h^*_J(z)\ge c_3$.
Since this applies to every $z \in (1-c_2/2)A'$,
the mean value property for $h^*_J$ gives
$h^*_J(0)\ge c_3 \length(A')/(2\pi)$.
Since $\bl|h^*_J(0)-H_{D'}(0,J)\br|<c_2/4$, our claim that
$H_D'(0,J)$ is bounded away from $0$ is established.
Since $h$ is positive, harmonic $h(0)=1$ and $h\ge M/2$ on $J$,
this also gives the bound $M\le c(\eps)$.
   
Since $h$ is harmonic, Lemma~\ref{l.bdhit} with 
$\tilde \eps_1 = \eps_1 /2$ and 
$\eps_2:=\tilde\eps_1/c(\eps) $
instead of $\eps_1, \eps_2$ implies that 
if $\delta=\delta(\eps_1)$ is sufficiently
small, and $\inrad(D)$ is large enough
to guarantee that $W$ separates $V(\delta, \eps_1)$
from $u$ (in the graph-connectivity sense),
then $h(z)\le  M \eps_2 < \eps_1$ for all $z \in V (\delta, \eps_1)$:
  ~\eref{e.vanish}
holds.

Now apply Lemma~\ref{l.goharm} again to conclude
that there is a harmonic function $\hat h:\U\to[0,\infty)$
such that 
$$
\bl|\hat h\circ\psi_D(z)-h(z)\br|<\eps_1
$$
for all $z\in \V D$ such that
$|\psi_D(z)|<1-\delta/4$.
Set $\tilde h(z):= \hat h\bl((1-\delta/2)z\br)$.
We know that $\tilde h\ge 0$ in $\p\U$,
$\tilde h(z')\le 2\eps_1$ on the set 
$S:=\bl\{ z'\in\p\U: |z'-\psi_D(u)|\ge 2\eps_1\br\}$.
Consequently, the Poisson representation of $\tilde h$ gives
$$
\tilde h(z) =O(\eps_1)+\int_{\p\U\setminus S} \tilde h(z')\,
\frac{1-\bl|z\br|^2}
{\bl|z-z'\br|^2}
\,|dz'|
\,.
$$
Since $\tilde h(0)=1+O(\eps_1)$ and $\eps_1$ is arbitrary,
the Proposition follows.
\QED

\subsection{Some Green's function estimates}
\label{s.dGpf}

As opposed to Proposition \ref {p.phit}, Lemma~\ref{l.dG} 
requires only crude bounds.  
It is actually possible to prove that $G_0 (0,0) - G_m (0,0)$ is close
to $t_m$, but we do not need this result here. 

\proofof{Lemma \ref {l.dG}}
We start with~\eref{e.lowG}.
Let $S$ be the set of vertices in $\V D$ satisfying~\eref{e.vrange},
and assume $S\ne\emptyset$.
For a random walk starting from a vertex in $S$, there is probability
bounded away from zero that within $\inrad(D)^2$ steps
it will exit $D$.  This gives
\begin{equation}\label{e.dgup}
\sum_{w\in S} G_D(0,w) \le O(1)\,\inrad(D)^2\,.
\end{equation}
On the other hand, with probability bounded away from
zero, the number of steps into vertices in $S$ for the random walk
started at $0$ that is stopped on exiting $D$ is
greater than  $\inrad(D)^2$.  Therefore
\begin{equation}\label{e.dglow}
O(1)\,\sum_{w\in S} G_D(0,w) \ge \inrad(D)^2\,.
\end{equation}
By reversing the walk, we know that $G_D(0,w)=G_D(w,0)$.
Since $G_D(w,0)$ is harmonic on $\V D\setminus\{0\}$,
the Harnack Principle (i.e., $k=0$ in~\eref{difh})
can be used to show that $G_D(w,0)/G_D(w',0)=O(1)$ when $w,w'\in S$.
Combining this with $G_D(0,w)=G_D(w,0)$ and the estimates~\eref{e.dgup},
\eref{e.dglow} gives~\eref{e.lowG}. 

By Lemma~\ref{l.diambound}, we have
\begin{equation}\label{e.diambd}
\diam \bl(\psi_D\circ\gamma[0,m])=O(\delta)\,.
\end{equation}
In the following, we fix $\gamma[0,m]$ (that is, it will
be considered deterministic).
Let $z$ be the vertex where simple random walk from $0$
first exists $D_m$.  By considering what happens to the
random walk after first hitting $z$ we get the identity
$ G_0(0,v)-G_m(0,v) = \Eb{ G_0(z,v)}$
(where $G_0(z,v)=0$ for $z\notin\V D$, by definition).
Consequently,
$$
G_0(0,v)-G_m(0,v) \le
\Pb{z\in\gamma[0,m]}
\,
\max\bl\{ \Eb{ G_0(\gamma_j,v)}:j=1,\dots,m\br\}
\,.
$$
By~\eref{e.diambd}, the continuous harmonic measure from $0$
of $\psi_D\circ\gamma[0,m]$ in $\U$ is $O(\delta)$.
Therefore, the continuous harmonic measure from $0$ of $\gamma[0,m]$ in
$D$ is also $O(\delta)$.
As in the argument given in Section~\ref{s.phitpf},
this implies that if $\inrad(D)$ is large enough,
$\Pb{z\in\gamma[0,m]}=O(\delta)$.

Let $K$ denote the disk $\bl\{w\in D:|w-v|<\inrad(D)/10\br\}$
and fix some $j\in \{1,2,\dots,m\}$.
Since $\psi_D\circ\gamma[0,m]$ is contained
in $\closure\U\setminus \bl(1-O(\delta)\br)\U$.
It follows that the continuous harmonic measure
of $K$ from $\gamma_j$ in $D$ is $O(\delta)$.
If $\psi_D(\gamma_j)$ is sufficiently close to $\p\U$
(how close may depend on $\delta$),
then we can make sure that the corresponding
discrete harmonic measure $H_D(\gamma_j,K)$
is less than $\delta$, by Lemma~\ref{l.bdhit}.
If $\psi_D(\gamma_j)$ is not close to $\p\U$,
then when $\inrad(D)$ is large the bound $H_D(\gamma_j,K)\le O(\delta)$
follows by the convergence of the discrete harmonic measure
to the continuous harmonic measure, as we have seen before.
If $w\in \V D\cap K$ neighbors with a vertex outside of
$K$, then $G_0(w,v)=O(1)$ 
follows from~\eref{e.dglow} by translating $w$ to $0$.
Hence, $G_0(\gamma_j,v)= O(1)\,H_D(\gamma_j,K)=O(\delta)$.
Putting these estimates together completes the proof.
\QED

\subsection{Mixed boundary conditions}
\label{s.pfreepf}

Recall that Proposition \ref {p.free}, which we will now prove,
is not used in the proof of Theorem~\ref{t.lerw}.

\proofof{Proposition \ref {p.free}}
Suppose first that the distance between $A_0$ and $A_1$ is at least $\eps$.
Let $A_0^*$ and $A_1^*$ denote the two connected components
of $\eta\setminus (A_0'\cup A_1')$, such that the sequence
$(A_0',A_0^*,A_1',A_1^*)$ conforms to the counterclockwise
order along $\eta$.  This induces a corresponding partition
$E_2=E_0^*\cup E_1^*$ of $E_2$, according to whether or not the first
point on the edge is in $A_0^*$ or in $A_1^*$.

We need to use the discrete harmonic conjugate function $\hat k$
of $\hat h$.  To be perfectly precise, it is necessary to 
set some combinatorial infrastructure: we first define a (multi-) graph 
$\hat H$ and $\hat k$ will be defined on the planar dual
$\hat H^\dagger$ of $\hat H$.
The vertices of $\hat H$ are $V_H\cup\{v_0,v_1\}$ (where $v_0$ and
$v_1$ are new symbols not appearing in $V_H$).
As edges of $\hat H$ we take all the edges of $H$, and, additionally,
for every $j=0,1$ and every directed edge
$[v,u]$ in $E_j$, there is a corresponding edge $[v,v_j]$ in $\hat H$.
Finally, there is also the edge $[v_0,v_1]$ in $\hat H$.
Consider a planar embedding of $\hat H$ which extends the planar embedding
of $H$, such that $v_0$ and $v_1$ are in the unbounded component
of $\C\setminus H$.
Let $\hat H^\dagger$ denote the planar dual of $\hat H$.
Then there is a unique edge $[v_0^\dagger,v_1^\dagger]$ in $\hat H^\dagger$
which crosses $[v_0,v_1]$.
We choose the labels so that $v_j^\dagger$ naturally corresponds
to $A_j^*$, $j=0,1$. 
Set $\hat h(v_j):=j$, $j=0,1$.  If we consider $\hat h$ as a function
on $\hat H$, then it is discrete harmonic except at $v_0$ and $v_1$.
This easily implies (see, e.g.,~\cite{DehnTiling},
or, more explicitly,~\cite{\BSsqharm})
that there is a discrete harmonic conjugate $\hat k$
defined on the vertices
of $\hat H^\dagger$; that is, for every directed edge
$e=[u,v]$ in $\hat H$ if $\{u,v\}\ne \{v_0,v_1\}$, then the discrete
Cauchy-Riemann equation
$\hat h(v)-\hat h(u)=\hat k(v^\dagger)-\hat k(u^\dagger)$
holds, where $[u^\dagger,v^\dagger]$ is the edge
of $\hat H^\dagger$ intersecting $e$ from right to left.
In fact, $\hat k$ is harmonic in $\hat H^\dagger$ except at $v_0^\dagger$
and $v_1^\dagger$.  The function $\hat k$ is unique, up to an additive
constant.  We choose the additive constant so that
$\hat k(v_0^\dagger)=0$.  Since $\hat h\ge 0$, by considering
the neighbors of $v_0$ and the orientation, it follows
that $\hat k(v_1^\dagger)\ge 0$.

Consider a sequence $D_n$ of such domains satisfying $\inrad(D_n)\ge n$,
with arcs $\eta=\eta_n$ and such harmonic 
functions $\hat h_n,\hat k_n$.
Let $L_n$ denote the maximum value of $\hat k_n$, which
is the value of $\hat k_n$ on $v_1^\dagger$.

Since $\eps>0$ is fixed, we can consider a subsequence 
of $n\to\infty$ such that the arcs $A_0$ and $A_1$
converge to arcs $\tilde A_0$ and $\tilde A_1$ of length at least $\eps$,
and the distance between them is at least $\eps$.
Let $\tilde A_0^*$ and $\tilde A_1^*$ denote the two
components of $\p\U\setminus(\tilde A_0\cup\tilde A_1)$,
so that $\tilde A_0,\tilde A_0^*,\tilde A_1,\tilde A_1^*$ is
the positive order along $\p\U$ of these arcs.

We now separate the argument into two cases according to
whether or not $L_n>1$.  Suppose that $L_n>1$ for infinitely
many $n$ and take a further subsequence of $n$ such that
$L_n>1$ along that subsequence.  Then $\hat k_n/L_n$
and $\hat h_n/L_n$ are both bounded by $1$.
It follows from Lemma~\ref{l.goharm} that
after taking a further subsequence, if necessary,
there are harmonic functions
$h$ and $k$ on $\U$
such that $L_n^{-1}\hat h_n\circ\psi_{D_n}^{-1}\to  h$
and $L_n^{-1}\hat k_n\circ\psi_{D_n}^{-1}\to  k$
uniformly on compact subsets of $\U$
(appropriately interpreted, since $\hat h_n$ and $\hat k_n$
are only defined on vertices and dual vertices, 
not on every point of $D_n$).
Moreover, \eref{e.dconv} shows that
$h$ and $k$ are harmonic conjugates,
because the discrete Cauchy-Riemann equations tend to the
continuous Cauchy-Riemann equations.

By Lemma~\ref{l.bdhit}, it follows that
$k$ is respectively equal to $0$ and $1$ 
in the relative interior of $\tilde A_0^*, \tilde A_1^*$,
and similarly $h$ has boundary values $0$ and $1/ \tilde L$
in $\tilde A_0$ and $\tilde A_1$,
where $\tilde L:=\lim_{n\to\infty} L_n$
(where the limit is along the subsequence, and must exist
and be finite).
By Schwarz reflection, say, this implies that
$h$ and $k$ satisfy Neumann boundary
conditions in $\tilde A_0^* \cup \tilde A_1^*$ and $\tilde A_0 \cup \tilde A_1$,
respectively.
It now easily follows (e.g., from the maximum principle)
that $ h+i k$ is the (unique)
conformal map taking $\U$ to the rectangle
$[0,1/\tilde L]\times [0,1]$ which takes the four
arcs $\tilde A_0, \tilde A_0^*, \tilde A_1, \tilde A_1^*$ to the corresponding sides
of the rectangle.

The argument in the case where $L_n\le 1$ for infinitely many
$n$ proceeds in the same manner, except that one should not
divide $\hat h_n\circ \psi_{D_n}^{-1}$ and 
$\hat k_n\circ \psi_{D_n}^{-1}$ by $L_n$.

It remains to remove the assumption that the distance between
$A_0$ and $A_1$ is at least $\eps$.
Observe that the probabilistic description of $\hat h$ shows that it
is monotone increasing in $A_1'$ and monotone decreasing in $A_0'$.
Take $\eps'>0$ much smaller than $\eps$.
Then $\hat h(0)$ for the given configuration is bounded from
above by the value of $\hat h(0)$ for the configuration where arcs
of length $\eps'$ are removed at the two ends of $A_1$, and
$A_1'$ is adjusted accordingly.  Similarly,
$\hat h(0)$ is bounded from below by the value of $\hat h(0)$
for the configuration where such arcs are removed at the two ends of
$A_0$.  The difference between the value of $h$ for
original versus any of the modified configurations goes to zero as $\eps'\to 0$,
since $h$ depends continuously on $(A_1,A_2)$, as long as
the length of $A_1\cup A_2$ is not zero.
Consequently, we get the Proposition by applying the restricted 
version proved above with $\eps'$ in place of $\eps$ and by ``sandwiching''.
\QED

\section{Other lattices}\label{s.other}

For convenience and simplicity, the proofs up to now have been
written for the loop-erased random walk and UST Peano curve
on the square grid.  
The purpose of the present section
is to briefly indicate how to 
adapt the proofs to more general walks on more general grids.
In order to keep this section short, we will not 
try to consider the most general cases.

Let $L$ be a
(strictly two-dimensional) lattice in $\R^2$;
that is, $L$ is a discrete additive
subgroup of $\R^2$ that is not contained in a line.
Discrete means that there is some neighborhood of $0$
whose intersection with $L$ is $\{0\}$.
Suppose that $G$ is a planar graph whose vertices are the elements
of $L$, and $G$ is invariant under translation by elements of $L$.
That is, if $u,v\in L$ are neighbors in $G$ and $\ell\in L$, then
$\ell+u$ neighbors  $\ell+v$.
It is not hard to verify that there is a linear map taking $L$
to the triangular lattice such that neighbors in $G$ are mapped
to vertices at distance $1$.  In particular, as a graph,
$G$ is isomorphic to the triangular grid or to the square grid.

Let $N$ be the set of neighbors of $0$ in $G$, and let
$N':=\{0\}\cup N$.
Let $X$ be an $N'$-valued random variable,
and let $X_1,X_2,\dots$ be an i.i.d.\ sequence where each
$X_n$ has the same law as $X$.  Consider the random
walk 
$$
S_n:=\sum_{j=1}^n X_j
$$
on $G$.
We are interested in the situation where the scaling limit of
$S_n$ is standard Brownian motion.  For this purpose, we require that
$\Es{X}=0$ and that the covariance matrix of $X$ is the identity
matrix.  (Note that 
if the covariance matrix of $X$ is non-degenerate but not
equal to the identity, we can always apply a linear transformation
to the system to convert to the above situation.  Therefore,
what we say below also applies in that case, provided that we 
appropriately modify the linear complex structure on $\R^2$.)

Note that under these assumptions, the Markov chain corresponding
to the walk $S_n$ does not need to be reversible.  An interesting
particular example the reader may wish to keep in mind is
where $\Pb{X=\exp(2\,\pi\,i\,j/3)}=1/3$ for $j=0,1,2$.

\begin{theorem}\label{t.genlerw}
Theorem~\ref{t.lerw} applies to the loop-erasure of
the random walk $S_n$.
\end{theorem}

\proofof{Theorem~\ref{t.genlerw}}
An inspection of the proof Theorem~\ref{t.lerw}, including
all the necessary lemmas, shows that only the generalization
of the proof of Lemma~\ref{l.dde} to the present framework requires 
special justification, which is given below. \QED

\begin{lemma}\label{l.gdde}
 Let $\tau_r$ denotes the first time $n$ with
$|S_n| \geq r$.   There exists a constant
$C$, depending on $X$ but not on $r$, such that
for all $r \geq C$, $w\in N$ and $y\in L$, 
\[  \Bl|\Prob^0[S_{\tau_r} = y] -
    \Prob^w[S_{\tau_r} = y] \Br| \leq C\; r^{-1} \;
   \Prob^0[S_{\tau_r} = y].\]
\end{lemma}

Here, $\P^w$ denotes the law of the Markov chain
started from $w$; that is, the law of $\bl(S_n+w:n\in\N\br)$
under $\P=\P^0$.  This Lemma is clearly sufficient to
provide the necessary analogue of Lemma~\ref{l.dde} for $S_n$.

\proof
There are various ways to prove the lemma (via coupling for instance).
We give here a proof based on Green's functions, as in~\cite{\Lbook}.
Without loss of generality, we assume that $\P[X=0]>0$ and that
$L$ is the minimal lattice containing $\bl\{w\in N: p_w>0\br\}$.
Then the random walk is irreducible on $L$.  
The discrete Laplacian $\Delta_X$ associated with $X$ is defined by
$$
\Delta_X f(z) := \Eb{f(z+X)}-f(z)\,.
$$
Let  $a$ be the potential kernel for
the random walk,
\[  a(z) := 
\sum_{j=0}^\infty\Bl(\Prob^0[S_j = 0]
          - \Prob^0[S_j = -z]\Br)
. \]
It is known that the series converges, and, in fact
\begin{equation}\label{e.az}
  a(z) = c_1 \log |z | + c_2 + O(|z|^{-1})\,,
\end{equation}
as $|z|\to\infty$, $z\in L$,
(where $c_1, c_2$ depend on the law of $X$).
This is proved in~\cite{\FukaiUchiyama} for the lattice $\Z^2$ 
with arbitrary nondegenerate covariance matrix (with
an appropriate dependence on the matrix), so the above
follows for other $L$ by applying a linear transformation.
Since $\P^0[S_j=-z]=\P^z[S_j=0]$, it follows that
$$
\Delta_X a(z)=\begin{cases} 1&z=0\,,\\
0&z\ne 0\,.
\end{cases}
$$
Let $G_r$ denote the Green's function for the walk
in $L\cap r\U$; that is,
$G_r(z,z'):=\sum_{j\in\N}\P^z[j<\tau_r,\,S_j=z']$.
Note that for all
$z,w\in L\cap r\U$,
\begin{equation}
\label{e.aG}   a(z-w) + G_r(z,w) = \E^z\bl[a(S_{\tau_r}-w)\br] , 
\end{equation}
since for fixed $w$ both sides are $\Delta_X$-harmonic for $z\in L\cap r\U$
and equality holds for $z \in L\setminus r\U$.
Set $M:= \max\{|w|:w\in N\}$ and $Z=\bl\{z\in L: r/2\le|z|<M+r/2\br\}$.
By~\eref{e.az} and~\eref{e.aG} $G_r(z,w) = c_1 \log 2 + O(r^{-1})$
for $z\in Z$ and $w\in N'$.
The same argument applied to the reversed walk $-S_j$,
which has potential kernel $\bar a(z) = a(-z)$
and Green's function $\bar G_r(z,w) = G_r(w,z)$,
gives 
\begin{equation}\label{e.Gappx}
\forall w\in N'\ \forall z\in Z\qquad
 G_r(w,z) = c_1 \log 2 + O(r^{-1}).
\end{equation}
Assuming $r>4\,M$, by considering the last vertex in $Z$ visited
by the walk before time $\tau_r$ we obtain for all $w\in N'$
and all $y\in L$
\begin{equation*}
   \Prob^w[S_{\tau_r} = y] =
   \sum_{z\in Z} 
             G_r(w,z) \, %
\P^z\!\bl[S_{\tau_r}=y,\, \min\{j\ge 1:S_j\in Z\}>\tau_r\br] .
\end{equation*}
Together with~\eref{e.Gappx}, this completes the proof of the lemma.
\QED

Observe that Theorem~\ref{t.lerw} also holds for the simple
random walk on the honeycomb grid, because two steps on the
honeycomb lattice are the same as a single step on a triangular
grid containing every other vertex on the honeycomb grid,
and so Lemma~\ref{l.gdde} may be applied.

\medskip

We now turn our attention towards spanning trees, and the 
generalizations of 
Corollary~\ref {c.ust} and Theorem~\ref{t.intropeano}. 
Suppose that $X$ and $-X$ have the same distribution, so that 
the walk $S$ is reversible. For an edge $e=[x,y]$, define
$p_e = \P [ X = y-x ] = \P [ X = x-y]$.
In this case, it is easy to generalize the definition of UST
to a measure on trees related to the law of $X$. This can 
be done either using Wilson's algorithm,
or equivalently by giving to each tree $T$ a probability 
that is proportional to the product of the transition 
probabilities along the edges of $T$. In other words,
$ 
\P [T] = Z^{-1} \prod_{e \in T}
p_e$,   
where $Z$ is a normalizing constant.
(The equivalence is proved in \cite{\WilsonAlg}; see also~\cite{\Llerw}.)
We call this the UST corresponding to the walk $S$ (even if this  
probability measure is not uniform).
Note that Lemma~\ref{l.pmark}  holds also 
in the present setting because the probability $\P[T]$
is given in terms of a product.

\begin{theorem}\label{t.reversegen}
Assuming that $-X$ has the same distribution
as $X$ (i.e., $S_j$ is reversible), %
Corollary~\ref{c.ust} and Theorem~\ref{t.intropeano} hold for the
UST corresponding to the walk $S$.
\end{theorem}
\proof 
The proof of Corollary~\ref{c.ust} holds in this generality.
In the proof of Theorem~\ref{t.intropeano}, the only significant
changes concern the discrete harmonic conjugate function,
used in the proof of Proposition~\ref{p.free}.
Recall that there as an appropriate definition for
the discrete harmonic conjugate for reversible walks on planar
graphs, where the discrete Cauchy-Riemann equation is
modified (see~\cite{DehnTiling} or~\cite[\S 6.1]{\KenyonDD}). %
If $G$ is graph-isomorphic to the square grid, the same
is true for the dual graph.  If $G$ is graph-isomorphic
to the triangular grid, then the dual is graph-isomorphic
to the honeycomb grid.
As pointed out above, Lemma~\ref{l.gdde} may therefore be applied to the
harmonic conjugate.
The details are left to the reader.
\QED

In the non-reversible setting, instead of a spanning tree,
one should consider a spanning arborescence, which is an
oriented tree with a root and the edges are oriented towards the root.
Fix a finite Markov chain with state space $V$ and a root $o\in V$.
Consider the measure on spanning arborescences of $V$ with root $o$,
where the probability for $T$ is proportional to the product of
the transition probabilities along the directed edges of $T$.
This is the analogue of the UST in the non-reversible setting.
Wilson's algorithm holds in this generality~\cite{\WilsonAlg},
however, the choice of the root $o$ clearly matters.

If we consider a finite piece of the lattice $L$,
and we wire part or all of the boundary, it is natural
to pick the wired vertex as the root.  With this convention,
Corollary~\ref{c.ust} holds for the wired tree.
It would be interesting to see if the free tree with root
chosen at $0\in D$ is invariant under conformal maps preserving $0$,
say (in the non-reversible setting).
Of course, one needs to choose a grid approximation of
$D$ where there is an oriented path from each vertex to
the root $0$.

In the proof of Theorem~\ref{t.intropeano} we have used reversibility
in two places.
The proofs of Theorems 10.7 and 11.1 of~\cite{\SchSLE}, which we quoted,
currently require reversibility.  However, these results were
only used to improve the topology of convergence to \SLE/.
More seriously, Section~\ref{s.pfreepf} uses the conjugate harmonic
function, whose definition in the non-reversible setting is not clear. 
Notwithstanding the obstacles,
it seems likely that these results can be proven
in the non-reversible setting too.

\bibliographystyle{halpha}
\addcontentsline{toc}{section}{Bibliography}
\bibliography{mr,prep,notmr}

\bigskip

\filbreak
\begingroup
\small
\parindent=0pt

\vtop{  %
\hsize=2.3in  %
Greg Lawler\\  %
Department of Mathematics\\  %
310 Malott Hall\\  %
Cornell University\\  %
Ithaca, NY 14853-4201, USA\\  %
{lawler@math.cornell.edu}  %
}  %
\bigskip  %
\vtop{  %
\hsize=2.3in  %
Oded Schramm\\  %
Microsoft Corporation\\  %
One Microsoft Way\\  %
Redmond, WA 98052, USA\\  %
{schramm@microsoft.com}  %
}  %
\bigskip  %
\vtop{  %
\hsize=2.3in  %
Wendelin Werner\\  %
D\'epartement de Math\'ematiques\\  %
B\^at. 425\\  %
Universit\'e Paris-Sud\\  %
91405 ORSAY cedex, France\\  %
{wendelin.werner@math.u-psud.fr}  %
}  %
\endgroup  %

\filbreak
\end{document}